%% file: IOPR-ToricKato-JEP.tex
\documentclass[10pt, a4paper]{amsart}

\usepackage{geometry}
\usepackage{a4wide}
\usepackage[all,cmtip]{xy}

\usepackage{numprint}

\usepackage[english]{babel}
\usepackage{scalerel, mathtools}
\usepackage[backref=page]{hyperref}
\usepackage{microtype}
\usepackage{fancyhdr}
\usepackage{amsmath, amssymb,color}
\usepackage{amsthm}

\usepackage[T1]{fontenc}
\usepackage[utf8]{inputenc}
\usepackage{lmodern}\normalfont
\DeclareFontShape{T1}{lmr}{bx}{sc} { <-> ssub * cmr/bx/sc }{}

\newcommand{\Inv}{\mathop\mathrm{Inv}}
\newcommand{\DD}{\mathbb{D}}

\newcommand{\abs}[1]{\vert #1\vert}

\newcommand{\ov}{\overline}

\newcommand{\un}{\underline}
\newcommand{\BB}{\mathbb{B}}
\newcommand{\Spec}{\mathop\mathrm{Spec}}

\newcommand{\e}{\mathrm{e}}

\newcommand{\del}{\partial}

\newcommand{\Int}{\mathrm{Int}}

\newcommand{\Hom}{\mathrm{Hom}}

\newcommand{\Aut}{\mathrm{Aut}}

\newcommand{\id}{\operatorname{\text{\sf id}}}

\newcommand{\Pic}{\operatorname{Pic}}
\newcommand{\coker}{\operatorname{coker}}
\newcommand{\im}{\operatorname{im}}

\newcommand{\Bl}{\operatorname{Bl}}
\newcommand{\orb}{\operatorname{orb}}
\newcommand{\ord}{\operatorname{ord}}

\newcommand{\Oo}{\mathcal{O}}

\newcommand{\TT}{\mathbb{T}}

\newcommand{\al}{\alpha}
\newcommand{\be}{\beta}

\newcommand{\la}{\lambda}

\newcommand{\GL}{\mathrm{GL}}

\renewcommand{\phi}{\varphi}

\newcommand{\ce}{\mathcal{C}^\infty}
\newcommand{\CC}{\mathbb{C}}

\newcommand{\RR}{\mathbb{R}}
\newcommand{\Ss}{\mathbb{S}}

\newcommand{\ZZ}{\mathbb{Z}}
\newcommand{\NN}{\mathbb{N}}
\newcommand{\PP}{\mathbb{P}}

\newcommand{\QQ}{\mathbb{Q}}

\newcommand{\Ka}{K\"{a}hler}

\newcounter{Mycounter}[section]
\newcounter{lemma}[section]
\newcounter{claim}[section]

\newcounter{fact}[section]

\newcounter{sublemma}[section]

\newcounter{corollary}[section]

\newcounter{theorem}[section]

\newcounter{conjecture}[section]

\newcounter{proposition}[section]

\newcounter{definition}[section]

\newcounter{example}[section]

\newcounter{remark}[section]

\newcounter{problem}[section]

\newcounter{question}[section]

 \stepcounter{thmx}

\makeatletter

\@addtoreset{equation}{section}

\@addtoreset{footnote}{section}

\makeatother

\renewcommand*{\backref}[1]{}
\usepackage{enumerate}
\def\eqref#1{(\ref{#1})}

\def\1{\sqrt{-1}\:}

\newcommand{\cntrct}                
{\hspace{2pt}\raisebox{1pt}{\text{$\lrcorner$}}\hspace{2pt}}

\makeatletter
\newcommand*\rel@kern[1]{\kern#1\dimexpr\macc@kerna}
\newcommand*\widebar[1]{%
  \begingroup
  \def\mathaccent##1##2{%
    \rel@kern{0.8}%
    \overline{\rel@kern{-0.8}\macc@nucleus\rel@kern{0.2}}%
    \rel@kern{-0.2}%
  }%
  \macc@depth\@ne
  \let\math@bgroup\@empty \let\math@egroup\macc@set@skewchar
  \mathsurround\z@ \frozen@everymath{\mathgroup\macc@group\relax}%
  \macc@set@skewchar\relax
  \let\mathaccentV\macc@nested@a
  \macc@nested@a\relax111{#1}%
  \endgroup
}
\makeatother

\makeatletter
\def\cleardoublepage{\clearpage\if@twoside \ifodd\c@page \else\hbox{}\thispagestyle{empty}\newpage
\if@twocolumn\hbox{}\newpage\fi\fi\fi} \makeatother

\usepackage{comment}

\graphicspath{{Images/}}

\newcommand{\executeiffilenewer}[3]{%
	\ifnum\pdfs
	
	trcmp{\pdffilemoddate{#1}}%
	{\pdffilemoddate{#2}}>0%
	{\immediate\write18{#3}}\fi%
}

\usepackage{xparse}			

\newcommand{\nN}{\mathbb{N}}
\newcommand{\nZ}{\mathbb{Z}}
\newcommand{\nR}{\mathbb{R}}
\newcommand{\nC}{\mathbb{C}}

\newcommand{\one}{{1\hspace*{-0.8 mm}\mathrm{l}}}

\newcommand{\annfan}[1]{\hat{#1}^0}

\NewDocumentCommand{\colfan}{ +m +m +m +O{} }
{
	\Sigma_{#1}^{#2 \,{}^\centerdot #3}
}
\NewDocumentCommand{\colvar}{ +m +m +m +O{} }
{
	\wt{X}_{#1}^{#2 \,{}^\centerdot #3}
}
\NewDocumentCommand{\colmap}{ +m +m +m +O{p} }
{
	{#4}_{#1}^{#2 \centerdot #3}
}
\NewDocumentCommand{\idemap}{ +m +m +O{} }
{
	\sigma{#3}_{#1}^{#2}
}

\newcommand{\sminus}{\setminus}
\newcommand{\wt}[1]{\widetilde{#1}}

\newcommand{\transp}[1]{{#1}^T}

\newcommand{\ui}{\mathfrak{i}}

\usepackage{todonotes}

\title{Toric Kato manifolds}
\author[N. Istrati]{Nicolina Istrati}
\address[Nicolina Istrati]{FB 12/Mathematik und Informatik, Philipps-Universität Marburg, Hans-Meerwein-Str. 6, 35032 Marburg, Germany}
\email{nicolina.istrati@uni-marburg.de}

\author[A. Otiman]{Alexandra Otiman}
\address[Alexandra-Iulia Otiman]{Roma Tre University, Department of Mathematics and Physics, Largo San
	Leonardo Murialdo, Rome, Italy AND
	Institute of Mathematics “Simion Stoilow” of the Romanian Academy, 21, Calea Grivitei,
	010702, Bucharest, Romania 
}
\email{aiotiman@mat.uniroma3.it, alexandra.otiman@imar.ro}

\author[M. Pontecorvo]{Massimiliano Pontecorvo}
\address[Massimiliano Pontecorvo]{Roma Tre University, Department of Mathematics and Physics, Largo San
	Leonardo Murialdo, Rome, Italy}
\email{max@mat.uniroma3.it}

\author[M. Ruggiero]{Matteo Ruggiero}
\address[Matteo Ruggiero]{Université Paris Cité, Sorbonne Université, CNRS, Institut de Mathématiques de Jussieu-Paris Rive Gauche, F-75006 Paris, France
}
\email{matteo.ruggiero@imj-prg.fr}

\begin{document}

\begin{abstract}
	 We introduce and study a special class of Kato manifolds, which we call {\it toric Kato manifolds}. Their construction stems from toric geometry, as their universal covers are open subsets of toric algebraic varieties of non-finite type. This generalizes previous constructions of Tsuchihashi and Oda, and in complex dimension 2,  retrieves the properly blown-up Inoue surfaces. We study the topological and analytical properties of toric Kato manifolds and link certain invariants to natural combinatorial data coming from the toric construction. Moreover, we produce families of flat degenerations of any toric Kato manifold,  which serve as an essential tool in computing their Hodge numbers. In the last part, we study the Hermitian geometry of Kato manifolds. We give a characterization result for the existence of locally conformally K\" ahler metrics on any Kato manifold.  Finally, we prove that no Kato manifold carries balanced metrics and that a large class of toric Kato manifolds of complex dimension $\geq 3$ do not support pluriclosed metrics. 
\end{abstract}
\maketitle

\section{Introduction}


Kato manifolds are compact complex manifolds of non-K\" ahler type and were introduced by M. Kato in 1977 \cite{kato} as manifolds containing a global spherical shell. More specifically, he showed that any such manifold, which we shall call henceforth a {\it Kato manifold}, is constructed in the following way. Let $\pi: \hat{\BB} \rightarrow \BB$ be a
modification of the standard unit ball in $\CC^n$ at finitely many points and let $\sigma: \overline{\BB} \hookrightarrow \hat{\BB}$ be a holomorphic embedding. Glue small neighborhoods of the two boundary components of $\hat{\BB} \sminus \sigma(\BB)$ via the local biholomorphism $\sigma \circ \pi$. The resulting manifold $X(\pi,\sigma)$ is a compact complex manifold with infinite cyclic fundamental group. The couple $(\pi, \sigma)$ is referred to as a {\it Kato data}. 

A Kato manifold has an associated germ given by $F:=\pi \circ \sigma: (\CC^n, 0) \rightarrow (\CC^n, 0)$. Many of its properties are encoded in $F$, however the germ alone does not fully characterize the manifold, unless the complex dimension is $2$ and the surface is minimal. This comes from the fact that for $n=2$, any modification is a sequence of blow-ups at points, while in higher dimension modifications can be much more complicated. As a result, Kato surfaces are much studied and fairly well understood (see for instance \cite{naka84}, \cite{dl84}, \cite{dot},\cite{tel}, \cite{fp} etc.), while in higher dimension very little is known.  

One main motivation for studying Kato manifolds comes from the fact that a big part of them carry locally conformally \Ka\ metrics, but no Vaisman metrics \cite{bru}, \cite{iop}. They constitute therefore an immense source of manifolds where one can study the interplay between the existence of special Hermitian metrics and different geometric and cohomological properties.  

In \cite{iop} we considered the simplest class of Kato manifolds, which corresponds to the case when $\pi: \hat{\BB} \rightarrow \BB$ is given by successive blow-ups at special points and $\sigma$ is a standard chart of a blow-up. We studied several of their analytical invariants, but did not have the tools to compute any Hodge number.

In the present paper, we introduce a much more general class,  the {\it toric Kato manifolds}. They are determined by a toric Kato data, namely a smooth toric modification $\pi: \hat{\CC}^n \rightarrow \CC^n$ at $0$ and a chart $\sigma: \CC^n \rightarrow \hat{\CC}^n$ satisfying a natural $(\CC^*)^n$-equivariance property. The corresponding germ has the form
\begin{equation}\label{germ}
F(z)=F_{\un\la,A}(z):=\un\la z^A=(\la_1z_1^{a_{11}}\cdots z_n^{a_{1n}},\ldots ,\la_nz_1^{a_{n1}}\cdots z_n^{a_{nn}})
\end{equation}
where $\un\la=(\la_1,\ldots, \la_n)\in(\CC^*)^n$ and  $A=(a_{jk})\in\GL(n,\ZZ)$. 

Starting from the fan of $\hat\CC^n$, viewed as a toric variety, and the embedding $\sigma$, we construct a natural infinite fan $\Sigma_A$ endowed with a $\ZZ$-action. Then using the germ $F$, we define a $\ZZ$-invariant open set in the toric variety of non-finite type $\wt{X}_{\un\la, A}\subset X(\Sigma_A,\ZZ^n)$. Our starting point is the following description of toric Kato manifolds (see \ref{secondDescr} for a more precise statement).

\begin{theorem}
If $X$ is any toric Kato manifold with germ $F_{\un\la,A}$, then we have a biholomorphism:
\begin{equation*}
X\cong\wt{X}_{\un\la,A}/{\mathbb{Z}}.
\end{equation*}
\end{theorem}

This second point of view of toric Kato manifolds generalizes the toric description of Inoue surfaces \cite{in1}, \cite{in2} given by Oda \cite[Section~14]{oda},  and a construction given by Tsuchihashi in \cite{tsu} of a class of non-K\" ahler manifolds with infinite cyclic fundamental group. Other generalisations of Tsuchihashi's construction, going in different directions, have been considered in the literature.  
For instance, Sankaran \cite{san} generalizes \cite{tsu} to higher rank fundamental group, while Battisti and Oeljeklaus \cite{bo} produce new manifolds by interpolating between the Sankaran construction and the LVMB manifolds of \cite{b}. 


The above description  allows for a better conceptual understanding of toric Kato manifolds, and in particular unveils new ways for computing different invariants using classical techniques from toric geometry. For instance, denoting by $a_j$ the number of $j$-dimensional cones of the fan of $\hat\CC^n$, we have the following  (see \ref{bettith}).

\begin{theorem}
	The $n$-dimensional toric Kato manifold $X$ has the following Betti numbers:
	\begin{gather*}
	b_0(X)=b_1(X)=b_{2n-1}(X)=b_{2n}(X)=1\\
	b_{2j+1}(X)=0, \ \ 1\leq j\leq n-2\\
	b_{2j}(X)=-1+\sum_{s=j}^n(-1)^{s-j}\binom{s}{j}\left(a_{n-s}+\binom{n}{s+1}\right), \ \ 1\leq j\leq n-1.
	\end{gather*}
\end{theorem}

The matrix $A$ in \eqref{germ} has only non-negative coefficients. Consider the maximal sub-matrix $P$ of $A$ which is a permutation matrix. When $P=A$, the resulting toric Kato manifold is a primary Hopf manifold. When $P\in \GL(n-1,\ZZ)$, the universal cover of the corresponding toric Kato manifold $X$ is all of $X(\Sigma_A,\ZZ^n)$, and we call $X$ of \textit{parabolic type}. In all the other cases, $\wt{X}$ is a proper subset of $X(\Sigma_A,\ZZ^n)$ and we call $X$ of \textit{hyperbolic type}. Imitating Nakamura's classification of Inoue surfaces \cite{naka84} in terms of curves, we have the following characterization (see \ref{hiperpara}).

\begin{theorem}
Let $X$ be a toric Kato manifold. 
\begin{enumerate}
\item $X$ is a primary Hopf manifold if and only if any of its $(\CC^*)^n$-invariant curves is elliptic; 
\item $X$ is of hyperbolic type if and only if any $(\CC^*)^n$-invariant curve is rational;
\item $X$ is of parabolic type if and only if $X$ contains a unique $(\CC^*)^n$-invariant elliptic curve, and at least one rational $(\CC^*)^n$-invariant curve. 
\end{enumerate}
\end{theorem}

Concerning analytic invariants, we determine explicitly the canonical line bundle in terms of the maximal toric divisor and show that the Kodaira dimension of any toric Kato manifold is negative (\ref{Kodairadim}). Also, we compute several Hodge numbers (see \ref{holF}, \ref{cohOX} and \ref{cohOX2}).

\begin{theorem}
Let $X$ be a toric Kato manifold with divisor $D$ induced by the exceptional divisor of $\pi$. Then one has  $H^0(X,\Omega_X^p)=0$ for any $p\geq 1$. If moreover $X$ is of hyperbolic type, then one has:
		\begin{gather*}
		h^{0,0}(X)=h^{0,1}(X)=1,\  h^{0,p}(X)=0, \ p\geq 2\\
		h^{1,p}(X)=0, \ p\neq 1\\
		\ h^{1,1}(X)=b_2(X)=\sharp \{ \text{irreducible components of }\, D\}>0.
		\end{gather*}
\end{theorem}

We note that the primary Hopf manifolds have the same Hodge numbers as above, and that our proof also works for parabolic type manifolds with $|\un\la|$ small enough (see \ref{invParabolic}). In particular, for all these cases one has $b_k=\sum_{p+q=k}h^{p,q}$ for $k\leq 2$. We suspect that all Kato manifolds should have the Hodge numbers of the above theorem, however showing this is beyond our means since much of our proof is based on the toric description of our manifolds. We should also note that \cite{san} computes all the Hodge numbers of the toric Kato manifolds with $P(A)=\emptyset$, under the additional assumption that $A$ is diagonalizable over $\CC$ and irreducible over $\QQ$. It would be interesting to know if one could generalize his methods to our general setting. 

A main tool for computing the Hodge numbers of toric Kato manifolds is developed in Section \ref{degen}, where we exhibit two different types of flat toric degenerations of toric Kato manifolds. The first one generalizes a previous construction of Nakamura \cite{naka83} for surfaces. In this case we deform, in a smooth family, any toric Kato manifold to a singular space obtained by identifying two invariant hypersurfaces on a smooth compact toric variety.   The second type of degeneration exists only in the hyperbolic case and extends a construction given by Tsuchihashi in \cite{tsu}. In this case, the singular fiber is again given by the identification of two invariant hypersurfaces in a compact, possibly singular, toric variety. A special feature of this second family is that all the smooth fibers are biholomorphic.

Finally, guided by the general principle of finding special Hermitian metrics on complex manifolds of non-K\" ahler type, we investigate the case of Kato manifolds. 
Special Hermitian metrics usually arise by imposing some power of the fundamental form $\Omega$ to be in the kernel of some differential operator. The existence of such metrics is far from being guaranteed on a  generic compact complex manifold. As a consequence, constant effort has been put to find new specific examples and to unravel the different restrictions that the existence of such metrics might impose.

In this direction, our first result consists in a characterization for the existence of locally conformally \Ka\ (lcK) metrics (see \ref{lcK}).

\begin{theorem} 
Let $(\pi:\hat\BB\rightarrow\BB,\sigma:\overline{\BB}\rightarrow\hat\BB)$ be a Kato data and let $X=X(\pi,\sigma)$ be the corresponding Kato manifold. 
The following are equivalent:
\begin{enumerate}
\item $X$ admits a locally conformally K\" ahler metric;
\item $\wt{X}$ admits a \Ka\ metric;
\item $\hat{\mathbb{\BB}}^n$ admits a K\" ahler metric.
\end{enumerate}
\end{theorem}

We should note here that in the case of a toric Kato manifold $X$, although the compact torus $\TT^n$ acts on the universal cover $\wt{X}$, this action never descends to $X$ if $X$ is not Hopf. Therefore,  $X$ is not a toric manifold in any classical sense. On the other hand, if $X$ admits an lcK metric, then $\wt{X}$ admits a \Ka\ metric with respect to which it becomes a toric \Ka\ manifold. In this manner, lcK toric Kato manifolds give interesting generalizations of the class of toric lcK manifolds. 

Our last result is concerned with the non-existence of other special Hermitian metrics (see Section~\ref{metrici} for the definitions, \ref{balanced} and \ref{pluri}).

\begin{theorem} 
	\begin{enumerate}
		\item A Kato manifold admits no strongly Gauduchon metric, and in particular no balanced metric and no Hermitian symplectic metric.
			\item If $X$ is a Kato manifold of dimension $n\geq 3$ satisfying $H^{1, 2}_{\overline{\partial}}(X)=0$, then it cannot be endowed with a pluriclosed metric. In particular, a toric Kato manifold of hyperbolic type does not admit pluriclosed metrics, unless it is a surface. 
	\end{enumerate}
\end{theorem}

The paper is organized as follows. We begin by some necessary preliminaries on Kato manifolds. The first definition of toric Kato manifolds, together with the elementary but technical properties of toric Kato data, are given in Section~\ref{toric}. In the next section, we describe and prove the second characterization of toric Kato manifolds. Next, in Section~\ref{geomprop}, we study geometrical properties of toric Kato manifolds, such as the invariant complex submanifolds, or the different natural divisors. Section~\ref{betti} is dedicated to the topological invariants and Section~\ref{analytic} to the computation of the Hodge numbers. In Section~\ref{degen}, we describe the degenerations of toric Kato manifolds needed for the next section. In Section~\ref{isom}, we tackle the problem of the classification of toric Kato manifolds.  Finally, Section~\ref{metrici} deals with the existence of special Hermitian metrics on a general Kato manifold. We conclude the section with an explicit family of examples of toric Kato manifolds in any complex dimension $n \geq 4$ that do not admit any lcK metrics.

\qquad

\section{Preliminaries on Kato manifolds}\label{prel}

In the present paper, a Kato manifold will be a compact complex manifold obtained as follows \cite{kato}, \cite{dl84}. We let $\BB:=\{(z_1, \ldots, z_n) \in \CC^n |  |z_1|^2+\ldots+|z_n|^2 < 1\}$ be the standard open ball in $\CC^n$. Consider a modification $\pi: \hat{\BB} \rightarrow \BB$ at $0$ and a holomorphic embedding $\sigma: \overline{\BB} \hookrightarrow \hat{\BB}$. Glue the two boundaries of $\hat{\BB} \sminus \sigma(\BB)$ via the real analytic CR-diffeomorphism $\gamma:=\sigma \circ \pi$. The result is a compact complex manifold $X(\pi, \sigma)$, named here a \textit{Kato manifold}. Any small neighborhood of the image of $\del\hat\BB$ in $X(\pi,\sigma)$ is a global spherical shell. We shall refer to the couple $(\pi, \sigma)$ as {\it Kato data} and to $F:=\pi \circ \sigma : (\BB, 0) \rightarrow (\BB, 0)$ as {\it the corresponding germ}. 

The class of Kato manifolds introduced in \cite{kato} is slightly more general, as $\pi$ is allowed to be a modification at more than one point, so that the resulting manifolds are proper modifications of the ones we described above \cite[Lemme~2.7, Part~I]{dl84}. However, for the present discussion this generality will not make any difference. 

We denote by $q:\wt{X}\rightarrow X$ the universal cover of a Kato manifold $X=X(\pi,\sigma)$. The deck group $\Gamma$ of $q$ is canonically isomorphic to $\ZZ$, and we will denote also by $\gamma$ the positive generator of $\Gamma$, which is indeed induced by the map $\sigma\circ\pi$.

We start by settling several general facts that will be needed in the paper. 

\begin{lemma}\label{interFB}
Let $F:(\BB,0)\rightarrow (\BB,0)$ be a holomorphic germ with $F(\ov\BB)\subset\BB$. Then $\bigcap_{m\in\NN}F^m(\BB)=\{0\}$.
\end{lemma}
\begin{proof}
By hypothesis there exists $0<r<1$ such that $F(\BB)\subset\BB_r$. By the Schwarz lemma \cite[Theorem~6]{shabat}, we have $||F(z)||\leq r||z||$ for any $z\in\BB$, hence $||F^m(z)||\leq r^m$ for any $m>0$ and any $z\in\BB$. This implies that $\lim_{m\to\infty}\sup_{z\in\BB}||F^m(z)||=0$, and hence $\mathrm{diam}\left(\bigcap_{m\in\NN}F^m(\BB)\right)=0$. The conclusion then follows.
\end{proof}

Given a map $F:(\CC^n,0)\rightarrow (\CC^n,0)$, define its stable set:
\begin{equation*}
W^s(F)=\{z\in\CC^n\mid \lim_{m\to\infty}||F^m(z)||=0\}.
\end{equation*}

\begin{lemma}\label{reuniuneFB}
For any holomorphic map $F:(\CC^n,0)\rightarrow (\CC^n,0)$ with $F(\ov \BB)\subset\BB$, we have $W^s(F)=\bigcup_{m\in\ZZ}F^m(\BB)$, where for $m=-k<0$, $F^{-k}(\BB)$ is the preimage of $\BB$ via $F^k$. 
\end{lemma}
\begin{proof}
The previous lemma shows that $\BB\subset W^s(\BB)$. Clearly $F^m(\BB)\subset W^s(F)$ for any $m>0$ and also for any $m<0$. Conversely, if $z\in W^s(F)$ then there exists $m>0$ such that $F^m(z)\in\BB$, which shows the desired equality.
\end{proof}

Suppose that $F:(\CC^n,0)\rightarrow(\CC^n,0)$ is given by a Kato data $(\pi,\sigma)$, $F=\pi\circ\sigma$, such that $\pi$ is a modification at $0$.
For any $k>0$, let $H_k:=F^{-k}(0)$ so that $H_k\subseteq H_{k+1}$ and put $H_\infty:=\bigcup_{k>0} H_k$. 
Let us set $\mathrm{Inv}(F):=\CC^n\sminus H_\infty$ and $W^*(F):=W^s(F)\cap\mathrm{Inv}(F)$.


\begin{proposition}\label{stableF}
Let $F:(\CC^n,0)\rightarrow(\CC^n,0)$ be a holomorphic map with $F(0)=0$ given by a Kato data $(\pi,\sigma)$ via $F=\pi\circ\sigma$, so that $\pi$ is a modification at $0$ with exceptional divisor $E$. Suppose that $H_\infty=H_m$ for some $m>0$. Let $D$ be the divisor on $X(\pi,\sigma)$ induced by $E$ after gluing. Then we have a biholomorphism $X(\pi,\sigma)\sminus D\cong W^*(F)/\langle F\rangle$.
\end{proposition}
\begin{proof}
If we put $\BB^*:=\BB\cap\mathrm{Inv}(F)$, then by \ref{reuniuneFB} we have $W^*(F)=\bigcup_{m\in\ZZ}F^m(\BB^*)$. Furthermore, by \ref{interFB}, we have:
\begin{equation*}
\bigcup_{m\in\ZZ}\left(F^m(\BB^*)\sminus F^{m+1}(\BB^*)\right)=\bigcup_{m\in\ZZ}F^m(\BB^*)\sminus \bigcap_{m\in\ZZ}F^m(\BB^*)=W^*(F).
\end{equation*}
Now it is clear that $D\cap(\hat\BB\sminus \sigma(\BB))=\pi^{-1}(H_\infty)\cap (\hat\BB\sminus \sigma(\BB))$ and that $\pi:\left(\hat\BB\sminus \sigma(\BB)\right)\sminus \pi^{-1}(H_\infty)\rightarrow \BB^*\sminus F(\BB^*)$ establishes a biholomorphism satisfying $\pi\circ\gamma=F\circ \pi$, where $\gamma=\sigma\circ\pi$. This implies in particular that the group $\langle F\rangle$ acts freely and properly on $W^*(F)$ and that $\pi$ induces the desired isomorphism $X(\pi,\sigma)\sminus D\cong W^*(F)/\langle F\rangle$.
\end{proof}

\begin{proposition}\label{embBall}
Let $(\pi:\hat\BB\rightarrow\BB,\sigma:\BB\rightarrow\hat\BB)$ be a Kato data with germ $F=\pi\circ\sigma$ and let $X=X(\pi,\sigma)$ be the corresponding Kato manifold. Then there exists a holomorphic open embedding $\phi:\hat\BB\sminus \{\sigma(0)\}\rightarrow \wt{X}$.
\end{proposition}
\begin{proof}
We recall that $\wt{X}=\bigsqcup_{m\in\ZZ} W_m \Big/\sim$, where $W_m=\hat\BB\sminus \sigma(\BB)$ for any $m\in\ZZ$ and $\sim$ indicates that $W_m$ is glued to $W_{m+1}$ via $\gamma=\sigma\circ\pi$.

Let us put, for any $l\geq 1$, $\wt{X}_l:=\bigsqcup_{1\leq m\leq l} W_m \Big/\sim  \, \,\subset \wt{X}$. Also, following \cite[Section~1]{iop}, denote by $(\pi_l:\hat\BB^{(l)}\rightarrow\BB,\sigma_l:\BB\rightarrow\hat\BB^{(l)})$ the Kato data obtained by composing $(\pi,\sigma)$ with itself $l$ times. Then it is easy to see that $\wt{X}_l=\hat\BB^{(l)}\sminus \sigma_l(\BB)$.

If we denote by $\pi_{(l)}:\hat\BB^{(l)}\rightarrow\hat\BB^{(l-1)}$ the induced map, then using that $\pi_{(l)}\circ\sigma_l=\sigma_{l-1}\circ F$ we find that $\pi^{-1}_{(l)}$ gives an embedding of $\hat\BB^{(l-1)}\sminus \sigma_{l-1}(F(\BB))$ into $\hat\BB^{(l)}\sminus \sigma_l(\BB)$. Therefore we find inductively, for any $l\geq 1$, open embeddings 
\begin{equation*}
\phi_l:Q_l:=\hat\BB\sminus \sigma(F^{l-1}(\BB))\rightarrow\hat\BB^{(l)}\sminus \sigma_l(\BB)\subset \wt{X}
\end{equation*}
which satisfy $\phi_{l+1}|_{Q_l}=\phi_l$. Now since $\bigcap_{l\geq 1} F^l(\BB)=\{0\}$ by \ref{interFB}, this means that $\bigcup_{l\geq 1} Q_l=\hat\BB\sminus \{\sigma(0)\}$, and thus the family $\{\phi_l\}_{l\geq 1}$ naturally defines an open embedding $\phi:\hat\BB\sminus \{\sigma(0)\}\rightarrow \wt{X}$, which concludes the proof.
\end{proof}

\section{Toric Kato manifolds}\label{toric}

We refer the reader to \cite{fulton}, \cite{oda}, \cite{oda2} or any other classical reference for the theory of toric algebraic varieties. Here, we will only fix notation. 

We denote by $N:=\ZZ^n$, $M:=\Hom_\ZZ(N,\ZZ)$, $N_\RR:=N\otimes_\ZZ\RR$, and
\begin{equation*}
\TT:=N\otimes_\ZZ\Ss^1=(\Ss^1)^n\subset T_N:=N\otimes_\ZZ\CC^*=(\CC^*)^n. 
\end{equation*}

We recall that any toric algebraic variety $X$ is a $T_N$-equivariant partial compactification of $T_N$, and its algebraic structure is encoded by a fan $\Sigma$ with support in $N_\RR$. We write $X=X(\Sigma,N)$. Traditionally, $\Sigma$ is a finite collection of rational cones with some compatibility properties, so that $X$ is covered by a finite number of affine open sets $(X_\tau)_{\tau\in\Sigma}$. However, in the present text, the fan $\Sigma$ is allowed to be infinite, and then $X$ is of non-finite type. 

We denote by $|\Sigma|:=\bigcup_{\tau\in\Sigma}\tau$ the support of $\Sigma$. We will let $\Sigma^{(k)}\subset\Sigma$ denote the subset of $k$-dimensional cones of $\Sigma$. Also, for a collection of vectors $v_j\in N_\RR$, $j\in J$, we will denote by $\langle v_j\mid j\in J\rangle:=\sum_{j\in J}\RR_{\geq 0}v_j$ the cone generated by it. 

Given a cone $\tau\in\Sigma^{(k)}$, we recall the notation $\check\tau:=\{l\in M_\RR\mid \langle l,v\rangle\geq 0\ \forall v\in\tau\}$, which is a cone in $M_\RR$, $\tau^\perp:=\{l\in M_\RR\mid \langle l,v\rangle=0\ \forall v\in\tau\}$, which is a face of $\check\tau$, and $S_\tau:=\check\tau\cap M$. Then one has
\begin{equation*}
X_\tau=\Spec\CC[S_\tau]=\Hom_{u.s.g.}(S_\tau,\CC)
\end{equation*}
where $u.s.g.$ stands for unit semi-groups. We also recall \cite[Theorem~4.2]{oda} that there exists a one-to-one correspondence between cones $\tau\in\Sigma$ and $T_N$-orbits on $X(\Sigma,N)$, given by $\orb\tau=\Hom_{gr.}(\tau^{\perp}\cap M,\CC^*)$, so that $\dim\tau=\mathrm{codim}\orb\tau$.

Finally, the map 
\begin{equation*}
\ord: \CC\rightarrow \RR\cup\{\infty\}, \ \ z\mapsto -\log|z|
\end{equation*}
restricts to a group homomorphism from $\CC^*$ onto $\RR$, and induces a surjective map of fiber $\TT$
\begin{equation*}
\ord:T_N\rightarrow N_\RR.
\end{equation*}

\subsection{Toric Kato data}
We let $z=(z_1,\ldots, z_n)$ denote the standard holomorphic coordinates on $\CC^n$. In all that follows, $\CC^n$ is viewed as a toric variety with the standard action of $T_N$. We will denote an element of the complex torus by $\un \lambda=(\la_1,\ldots,\la_n)\in T_N$, and for a toric variety $X$ and $x\in X$, we simply write $\un\la x$ to denote the action of $\un\la$ on $x$.

We will call a proper modification  $\pi:\hat\CC^n\rightarrow\CC^n$ at $0\in\CC^n$ \textit{a toric modification} if $\hat\CC^n$ is a $T_N$-toric variety and $\pi$ is $T_N$-equivariant. 

\begin{definition}
A Kato data $(\pi:\hat\BB\rightarrow\BB,\sigma:\BB\rightarrow\hat\BB)$ is called a \emph{toric Kato data} if the extension $\pi:\hat\CC^n\rightarrow\CC^n$ is a smooth toric modification at $0$ and there exists $\nu\in\Aut_{gr.}(T_N)$ such that $\sigma(\un \lambda z)=\nu(\un \lambda)\sigma(z)$ for any $z\in \BB$ and $\un \lambda\in T_N$ for which this is defined. In this case, we say that $\sigma$ is $\nu$-equivariant. A Kato manifold will be called a \emph{toric Kato manifold} if it admits a toric Kato data. 
\end{definition}
 
Let $\pi:\hat\CC^n\rightarrow\CC^n$ be a toric modification at $0$, let $\Sigma$ be the fan of $\CC^n$ and let $\hat{\Sigma}$ be the fan of $\hat\CC^n$. Note that $\hat{\Sigma}$ is uniquely determined by $\pi|_{\pi^{-1}(\BB)}$. Then $|\hat{\Sigma}|=|\Sigma|=(\RR_{\geq 0})^n=:C_0$. Moreover, since $\pi$ is a biholomorphism outside $0$, it follows that 
\begin{equation}\label{intRays}
|\hat{\Sigma}^{(1)}\sminus \Sigma^{(1)}|\subset\Int|\Sigma|
\end{equation}
i.e. each new ray of $\hat{\Sigma}$ has only positive components. Conversely, any finite refinement with regular cones $\hat{\Sigma}$ of $\Sigma$ satisfying \eqref{intRays} defines a toric modification of $\CC^n$ at $0$.

For every cone $\tau\in\hat{\Sigma}^{(n)}$, there exists a unique $\ZZ$-basis $\un e^\tau:=\{e^\tau_1,\ldots,e^\tau_n\}$ of $N$, such that $\langle e^\tau_j\rangle\prec \tau$, $j=1,\ldots, n$. If $\un f^\tau:=\{f_1^\tau,\ldots,f_n^\tau\}$ is the corresponding dual basis of $M$, then the map
\begin{equation}\label{defnPhiA}
S_\tau= \Bigl\{\sum_{j=1}^nu_jf_j^\tau\mid u_1\ldots, u_n\in\NN\Bigr\}\rightarrow \NN\langle e_1^*,\ldots, e_n^*\rangle,\ \  \sum_{j=1}^nu_jf_j^\tau\mapsto\sum_{j=1}^nu_je_j^*
\end{equation}
gives rise to a holomorphic chart
\begin{gather*}
\phi_\tau:\CC^n=\Spec\CC[e_1^*,\ldots, e_n^*]\rightarrow \Spec\CC[S_\tau]\subset \hat\CC^n.
\end{gather*}
This chart is uniquely determined up to composing with a permutation of the coordinates.  In what follows, we will call such a chart \emph{a toric chart}. 

Express the vectors of $\un e^\tau$ in the standard basis $e_1,\ldots, e_n$ of $N$, and put $A_\tau:=\begin{pmatrix} e_1^\tau &\cdots & e_n^\tau\end{pmatrix}\in\GL(n,\ZZ)$. Then $A_\tau$ completely determines the chart $\phi_\tau$, so we will also use the notation $\phi_\tau=\phi_{A_\tau}$. Note that $A_\tau$ has only non-negative coefficients, since $\tau\subset C_0$. Moreover, since $A_\tau=(a_{kl})_{k,l}$ satisfies
\begin{equation*}
e_j^*=\sum_{k=1}^na_{jk}f_k^\tau, \ \ 1\leq j\leq n
\end{equation*}
it follows that the map $\pi\circ\phi_\tau:\CC^n\rightarrow\CC^n$, which is induced by 
\begin{equation*}
\NN\langle e_1^*,\ldots, e_n^*\rangle\rightarrow \NN\langle e_1^*,\ldots, e_n^*\rangle, \ \ e_j^*\mapsto \sum_{k=1}^na_{jk}e_k^*
\end{equation*}
satisfies
\begin{equation}\label{germeneTau}
\pi\circ\phi_\tau(z)=z^{A_\tau}:=(z_1^{a_{11}}\cdots z_n^{a_{1n}},\ldots ,z_1^{a_{n1}}\cdots z_n^{a_{nn}}), \ \ z\in\CC^n.
\end{equation}
Using the equivariance of $\pi$, we immediately infer the equivariance relation:
\begin{equation}\label{equivariance}
\phi_\tau(\un \lambda z)=\un \lambda^{A_\tau}\phi_\tau(z), \ \ \forall \un \lambda\in T_N.
\end{equation}

Any other toric chart $\phi'_\tau$ for $X_\tau$ is given by $\phi'_{\tau}=\phi_\tau\circ\hat{s}$, where $s\in\mathcal{S}_n$ is a permutation and $\hat{s}(z)=(z_{s(1)},\ldots, z_{s(n)})=z^{P_s}$, with $P_s=\transp{\begin{pmatrix}
e_{s(1)} &\cdots &e_{s(n)}\end{pmatrix}}$.
Then the corresponding matrix for $\phi'_\tau$ is given by $A'_\tau=A_\tau P_s$.

\begin{lemma}\label{toricChart}
Let $(\pi,\sigma)$ be a toric Kato data. Let $\hat{\Sigma}$ be the fan determined by $\pi$, and let $\tau\in\hat{\Sigma}^{(n)}$ be the cone representing the $T_N$-fixed point $\sigma(0)\in\hat\BB$. Then there exists $\un\la_0\in T_N$ and a toric chart $\phi_A:\CC^n\rightarrow X_\tau$ so that $\sigma=\un\la_0\circ\phi_A$. In particular, the germ $F:=\pi\circ\sigma$ satisfies $F(z)=F_{\un\la_0,A}(z):=\un\la_0z^A$ for all $z\in\CC^n$.
\end{lemma}
\begin{proof}
Note that $\Aut_{gr.}(T_N)\cong\Aut_\ZZ(N)\cong\GL(n,\ZZ)$, and a matrix $A$ determines the group automorphism $\un\la\mapsto\un\la^A$. Hence, by hypothesis, there exists $A\in\GL(n,\ZZ)$ so that $\sigma(\un\la z)=\un\la^A\sigma(z)$, for all $\un\la\in T_N$ and $z\in\BB$ for which this is defined. 

Let $B\in\GL(n,\ZZ)$ be a matrix giving a toric chart $\phi_B:\CC^n\rightarrow X_\tau$. Then $f:=\phi_B^{-1}\circ\sigma:\BB\rightarrow\CC^n$ is a biholomorphism onto its image, and it satisfies, for any $\un\la\in T_N\cap\BB$:
\begin{align*}
f(\un\la)=f(\un\la_1(\un\la\un\la_1^{-1}))=\phi_B^{-1}(\un\la^A\un\la_1^{-A}\sigma(\un\la_1))=\un\la^{B^{-1}A}\un\la_1^{-B^{-1}A}f(\un\la_1)
\end{align*}
where $\un\la_1$ is some fixed element from $T_N\cap\BB$.
Thus, putting $\un\la_2:=\un\la_1^{-B^{-1}A}f(\un\la_1)\in T_N$, we find that $f(z)=\un\la_2 z^{B^{-1}A}$ on $\BB$. Now since $f$ has an holomorphic inverse, it is easy to check that $B^{-1}A$ is a permutation matrix $P_s$ for some $s\in\mathcal{S}_n$. Therefore $A=BP_s$ is a toric matrix and $\sigma=\un\la_0\circ\phi_A$ for $\un\la_0:=\un\la_2^B$, as expected. 
\end{proof}

\begin{lemma}\label{powerKD}
If $(\pi,\sigma)$ is a toric Kato data with germ $F_{\un\la_0,A}$ and $d>1$ is a natural number, then there exist a toric Kato data $(\pi^d,\sigma^d)$ naturally associated to $(\pi,\sigma)$, with corresponding germ $F^d_{\un\la_0,A}$. Moreover, $X(\pi^d,\sigma^d)$ is a cyclic unramified covering  of $X(\pi,\sigma)$ of degree $d$.
\end{lemma}
\begin{proof}
The Kato data $(\pi^d,\sigma^d)$ is obtained by composing $(\pi,\sigma)$ with itself $d$ times, as it was described in \cite[Section~1]{iop}, so that the resulting germ is the $d$-th power of the initial one. From the construction, it is clear that the resulting Kato data is again toric, and that the fan corresponding to $\pi^d$ is given by
\begin{equation*}
\hat{\Sigma}^d:=\{A^k\tau\mid 0\leq k\leq d-1, \ \tau\in\hat{\Sigma}_0\}\cup\{A^{d-1}\tau_A\}
\end{equation*} 
where $\hat{\Sigma}_0:=\hat{\Sigma} \sminus \{\tau_A\}$, with $\hat{\Sigma}$ the fan of $\pi$ and $\tau_A\in\hat{\Sigma}^{(n)}$ the cone representing $\sigma(0)$. Furthermore, by \cite[Lemma~1.7]{iop}, $X(\pi^d,\sigma^d)$ is a cyclic unramified covering of $X(\pi,\sigma)$ of degree $d$.
\end{proof}

\subsection{Properties of toric Kato matrices}
\qquad

\begin{definition}
A matrix $A\in\GL(n,\ZZ)$ will be called a \textit{toric Kato matrix} if there exists a toric modification $\pi:\hat\CC^n\rightarrow\CC^n$ at $0$ and a toric chart $\phi_A:\CC^n\rightarrow\hat\CC^n$ with $\phi_A(\un\la z)=\un \la^A\phi_A(z)$ for all $z\in\CC^n$ and $\un\la\in T_N$.
\end{definition}

\begin{proposition}\label{prop:toricKatomatrix}
A matrix $A\in\GL(n,\ZZ)$ is a toric Kato matrix if and only if its columns are either positive or standard vectors.
\end{proposition}
\begin{proof}
The fact that a toric Kato matrix satisfies the desired properties is a direct consequence of \eqref{intRays}. For the other implication, we need to prove that for any $A \in \GL(n, \ZZ)$ with either positive or standard columns, there exists a smooth toric modification $\pi: \hat{\mathbb{C}}^n \rightarrow \mathbb{C}^n$ at $0$, such that the cone $\tau_A$ generated by the columns of $A$ belongs to the corresponding fan $\hat{\Sigma}$ of $\hat{\mathbb{C}}^n$.
This can be achieved by first considering any simplicial fan containing the cone $\tau_A$ and whose rays are generated either by the standard vectors or positive vectors, and then regularizing it via barycentric subdivision (which adds only positive rays, see, e.g., \cite[p. 48]{fulton}).
\end{proof}

\ref{prop:toricKatomatrix} shows that our setting is more general than the one considered in \cite{iop}. In fact, toric modifications of $0 \in \nC^n$ need not dominate the blow-up of the origin, as the following example shows.

\begin{example}
Consider the toric Kato matrix
$A=\begin{pmatrix}
3 & 2 & 3 \\
2 & 1 & 2 \\
2 & 2 & 1
\end{pmatrix}
\in \GL(3,\nZ)$.
Then the cone $\tau_A$ generated by the three columns of $A$ intersects in its interior the cone generated by $e_1=(1,0,0)$ and $e_1+e_2+e_3=(1,1,1)$. Hence no toric modification $\pi:\hat{\nC}^3 \to \nC^3$ over $0$ for which the fan $\hat{\Sigma}$ of $\hat{\nC}^3$ contains the cone $\tau_A$ dominates the blow-up of the origin.

In Figure \ref{fig:example4} we give two different explicit realizations of toric modifications $\pi$ containing $\tau_A$.


Let us explain how the fan on the right can be described in terms of blow-ups and contractions. We start from the standard fan of $\nC^3$, where the rays are generated by the canonical basis $(1,0,0)$, $(0,1,0)$ and $(0,0,1)$. When performing a blow-up, we only indicate the generator $v$ of the ray associated to the exceptional divisor obtained via blow-up. This determines uniquely the blow-up, and the new fan $\hat{\Sigma}$ is obtained from the previous one $\Sigma$ by replacing any cone of $\Sigma$ containing $v$ by all the cones given as the sum of any of its proper faces and $v$. 

Firstly, we perform three point blow-ups corresponding to the rays generated by $(1,1,1)$, $(2,2,1)$ and $(2,1,2)$. 
Secondly, we perform the blow-up of four equivariant lines, corresponding to the rays generated by $(3,2,1)$, $(2,1,1)$, $(3,2,2)$ and $(5,3,3)$ respectively. Next, we notice that the last operation, corresponding to introducing the divisor $E$ given by the ray through $(5,3,3)$, is on the one hand the blow-up of the line corresponding to the cone $\langle (3,2,2),(1,1,1)\rangle$, but also the blow-up of the line corresponding to the cone $\langle (2,1,2),(3,2,1)\rangle$. Finally, this allows us to contract $E$ by erasing the ray generated by $(5,3,3)$ together with the faces connecting this ray to $(2,1,1)$ and $(3,2,2)$ and the induced changes on the other cones. The composition of the last blow-up and the contraction corresponds to a flip in the language of Minimal Model Program.

\begin{tiny}
\begin{figure}[h]
\centering
\begin{minipage}[htbp]{\columnwidth}
\def\svgwidth{0.48\columnwidth}
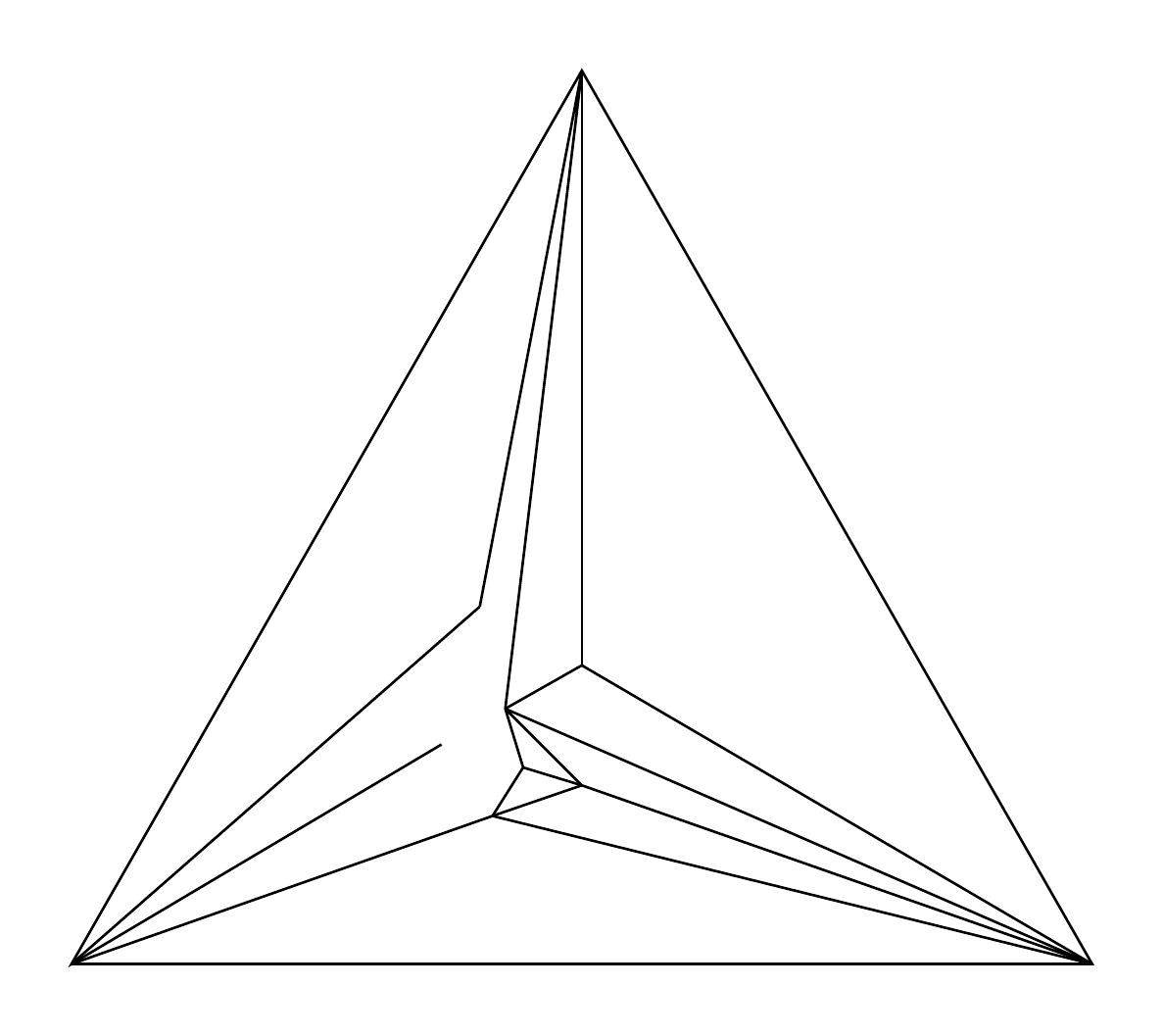
\hfill
\def\svgwidth{0.48\columnwidth}
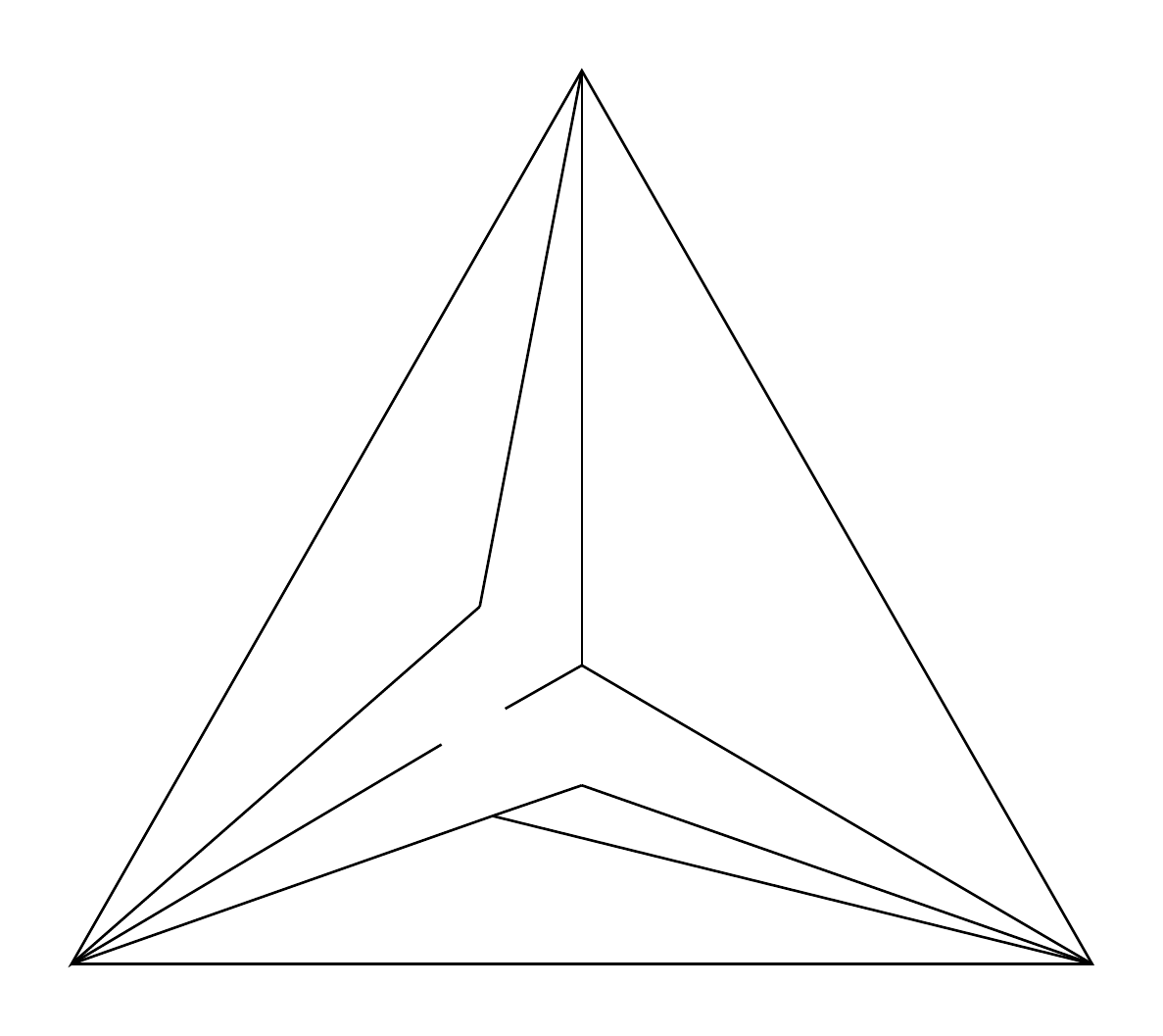
\end{minipage}
\caption{Toric Kato data with $\pi$ not dominating the blow-up of the origin.}
\label{fig:example4}
\end{figure}
\end{tiny}
\end{example}

Let $(\pi,\sigma=\un\lambda_0\circ\phi_A)$ be a toric Kato data, and let $\tau_A\in\hat{\Sigma}^{(n)}$ be the  cone generated by the columns of $A$, or equivalently, the cone corresponding to the fixed point $\sigma(0)$. Then $A$ is a positive matrix if and only if $\tau_A\subset\Int C_0$. In general, we have:

\begin{lemma}\label{permA}
Let $A=\begin{pmatrix}
A_1 &A_2 &\cdots &A_n\end{pmatrix}$ be a toric Kato matrix.
There exists a unique maximal subset $P(A)\subset\{1,\ldots, n\}$ satisfying that there exists a permutation $s:P(A)\rightarrow P(A)$ so that $A_j=e_{s(j)}$ for any $j\in P(A)$. Furthermore, there exists $m_0\geq 1$ so that for any $m\geq m_0$ and for any $j\in\{1,\ldots, n\}\sminus P(A)$, the $j$-th column of $A^m$ has strictly positive components.
\end{lemma}
\begin{proof}
Let us define 
\begin{equation*}
S_1(A):=\{j\in\{1,\ldots, n\}\mid A_j\in\{e_1,\ldots, e_n\}\}
\end{equation*} 
and let  $k:S_1(A)\rightarrow \{1,\ldots, n\}$ be the injective map satisfying  $A_j=e_{k(j)}$ for all $j\in S_1(A)$. Moreover, define inductively, for all $m\geq 2$, $S_m(A):=S_{m-1}(A)\cap k^{-1}(S_{m-1}(A))$. Then one has:
\begin{equation*}
\ldots\subseteq S_m(A)\subseteq S_{m-1}(A)\subseteq\ldots\subseteq S_1(A)
\end{equation*}
and so there exists $m_0\geq 1$ with $S_m(A)=S_{m_0}(A)$ for all $m\geq m_0$. One clearly has then $P(A)=S_{m_0}(A)$ and $s=k|_{P(A)}$.

On the other hand, it is easy to check that $S_m(A)=S_1(A^m)$ for all $m\geq 1$. Moreover, by \eqref{intRays}, a $j$-th column of $A^m$ has only positive components exactly when $j\notin S_1(A^m)$, from which the conclusion follows.
\end{proof}

For any matrix $A=(a_{kl})_{1\leq k,l\leq n}\in\GL(n,\ZZ)$ and $J\subset\{1,\ldots, n\}$, we denote by $J^c=\{1,\ldots, n\}\sminus J$ and we put $A_J:=(a_{kl})_{k,l \in J}$, which is the matrix obtained from $A$ by erasing the lines and columns prescribed by $J^c$. We also use the notation $\tau_J\in\Sigma^{(|J|)}$ for the cone generated by $e_j$ with $j\in J$ and we put $N_J:=\mathrm{span}_{\ZZ}\{e_j|j\in J\}\subset N$, respectively $T_J:=N_J\otimes_\ZZ\CC^*$, $\CC_J=N_J\otimes \CC$ etc. Finally, we denote by $p_J$ the natural projection from $N$ to $N_J$, as well as its $\RR$ or $\CC$-linear extension,   and for $v\in \CC^n$, we also write $v_J:=p_J(v)$.

\begin{lemma}\label{subKatomatrix}
Let $A=(a_{kl})_{1\leq k,l\leq n}$ be a toric Kato matrix and let $J\subseteq P(A)$ be any subset with $|J|<n-1$ which is fixed by the permutation $s$ given in \ref{permA}. Then $A_{J^c}$ is again a toric Kato matrix. 
\end{lemma}
\begin{proof}

Let $\pi:\hat\CC^n\rightarrow\CC^n$ be a toric modification at $0$ such that $A$ gives a toric chart in $\hat\CC^n$.  Let 
$\hat\CC_{J^c}:=\ov{\orb \tau_J}\subset\hat\CC^n$ be the strict transform of $\CC_{J^c}=\ov{\orb \tau_J}\subset\CC^n$ via $\pi$. Note that since $A$ preserves $N_{J}$, it follows that also the morphism of $T_N$ given by $\nu_A(\un\la)=\un\la^A$  fixes $T_{J}$. Since $\hat\CC_{J^c}$ is uniquely determined as the subset of $\hat\CC^n$ on which $T_{J}$ acts trivially, by the $\nu_A$-equivariance of  $\phi_A$ it follows that $\phi_A(\CC_{J^c})\subset\hat\CC_{J^c}$. Moreover, by the definition \eqref{defnPhiA} of $\phi_A$ we find that $\phi_A|_{\CC_{J^c}}=\phi_{A_{J^c}}$. Hence $A_{J^c}$ gives a toric chart in the toric modification $\pi|_{\hat\CC_{J^c}}:\hat\CC_{J^c}\rightarrow \CC_{J^c}$.
\end{proof}

\begin{lemma}\label{PerronNP}
Let $A=(a_{kl})_{1\leq k,l\leq n}\in\GL(n,\ZZ)$ be a toric Kato matrix with $|P(A)|<n-1$. Then $A$ has a simple real eigenvalue $\al>1$ so that for any other $\al\neq\be\in\Spec(A)$, we have $|\be|<\al$. Moreover, $A$ admits a Perron eigenvector $f_A$ with $Af_A=\al f_A$ so that $\langle e_j^*, f_A\rangle>0$ for all $1\leq j\leq n$. 
\end{lemma}
\begin{proof}
Let $s$ be the permutation of $P(A)$ given by \ref{permA} and let $d=\ord s$. Then for $m=m_0d$ big enough, $(A^m)_{P(A)}=\mathrm{Id}$ and  $B:=(A^m)_{P(A)^c}$ is a positive matrix. After replacing $A$ by some power, we suppose that this happens for $m=1$.

By the Perron-Frobenius theorem, $B$ has a simple real eigenvalue $\al$ and a Perron eigenvector $f_B$ with positive components. Since $|\det(B)|=1$ and $|P(A)^c|>1$, we have $\al>1$. 

For any $\be\in\CC$ and $k\times k$ matrix $M$, denote by $V_M(\be)\subset \CC^k$ the generalized eigenspace of $M$ for the eigenvalue $\be$. Put $P=P(A)$ and write $p:=p_{P^c}:
\CC^n\rightarrow \CC_{P^c}$ for the natural projection. Then $p$ is compatible with the splitting
\begin{equation*}
\CC^n=\bigoplus_{\be\in\CC}V_A(\be)\rightarrow \bigoplus_{\be\in\CC}V_B(\be)=\CC_{P^c}
\end{equation*}
and $\Spec(A)=\Spec(B)\cup\Spec(A_{P})$, where $\Spec(A_{P})=\{1\}$. In particular, we infer that $\al\in\Spec(A)$ is a simple eigenvalue and there exists an eigenvector $f_A\in\CC^n$ with $p(f_A)=f_B$. For any $k\in P$, we find:
\begin{equation*}
\langle e_k^*,f_A\rangle=\frac{1}{\al-1}\sum_{j\in P^c}a_{kj}\langle e_j^*, f_B\rangle>0
\end{equation*}
from which the conclusion follows.
\end{proof}

For a toric Kato matrix $A$ with $|P(A)|<n-1$, let $f_A$ be a Perron vector for $A$ with positive components. Then $\transp{A}$ acting on $(\CC^n)^*$ also has a Perron vector $f_A^*$ with $\langle f_A^*, f_A\rangle=1$. We note that since $\pm 1$ are the only possible rational eigenvalues of $A$, $f^*_A$ cannot be chosen rational. We also observe that since $\CC_{P(A)}$ is spanned by eigenvectors of $A$ with eigenvalue different from $\al$, for any $k\in P(A)$ we have $\langle f_A^*,e_k\rangle=0$, while for $k\in P(A)^c$, we have $\langle f_A^*, e_k\rangle >0$.
 Define the half-space:
\begin{equation*}
H(A):=\{v\in N_\RR\mid \langle f_A^*,v\rangle>0\}\subset N_\RR.
\end{equation*}
Also, for a toric Kato matrix with $P(A)^c=\{j\}$, put
\begin{equation*}
H(A):=\{v\in N_\RR\mid\langle e_j^*,v\rangle >0\}\subset N_\RR.
\end{equation*}

\qquad

\subsection{Toric Kato germs}

Note that for any toric Kato germ $F=F_{\un\la_0,A}$, if $m_0>0$ is such that $A^{m_0}$ has strictly positive components on the columns $j\in P(A)^c$, by \ref{permA} we have 
\begin{gather*}
H_\infty=\bigcup_{m>0}F^{-m}(0)=F^{-m_0}(0)=\bigcup_{j\in P(A)^c} \CC_{\{j\}^c}\\
\Inv(F)=\{z\in\CC^n\mid z_j\neq 0, \ j\in P(A)^c\}=\CC_{P(A)}\times T_{P(A)^c}\supseteq T_N.
\end{gather*}
Let us put $W_T(F):=W^s(F)\cap T_N\subseteq W^*(F)$. Then we have the following:

\begin{lemma}\label{splitWs}
Let  $F=F_{\un\la_0,A}$ be a toric Kato germ. Let $B:=A_{P(A)^c}$, $\un\la':=(\un\la_0)_{P(A)^c}$  and let $F'=F_{\un\la',B}$ be the induced toric Kato germ acting on $T_{P(A)^c}$. Then we have a biholomorphism:
\begin{equation*}
W_T(F)=T_{P(A)}\times W_T(F').
\end{equation*}
In particular, if $|P(A)|=n-1$ then $W_T(F)=T_N$.
\end{lemma}
\begin{proof}
As $W_T(F)=W_T(F^m)$ for any $m\in\NN^*$, we can suppose without loss of generality that $A_{P(A)}=\id$. Let us denote by $p:=p_{P(A)^c}$ and by $p':=p_{P(A)}$ the two projections. Since $p\circ F=F'\circ p$, it is clear that $p(W_T(F))\subseteq W_T(F')$. Denote by $\un\la'':=(\un\la_0)_{P(A)}$. Then we have
\begin{equation}\label{FmQ}
(F^m(z))_{P(A)}=(\un\la'')^mz_{P(A)}z_{P(A)^c}^{Q_m}, \ \ m\in\NN
\end{equation}
where $Q_m$ is an integer valued matrix determined by $A$, whose coefficients are each at least $m$.  From this, it is clear that the fiber of $p:W_T(F)\rightarrow p(W_T(F))\subset W_T(F')$ is $T_{P(A)}$. It remains to show that this restriction of $p$ is surjective. 

Let $w\in W_T(F')$ and let $0<\epsilon<1$. Since $w\in p(W_T(F))\Leftrightarrow (F')^m(w)\in p(W_T(F))$ $\forall m\in\NN$, we can suppose that $w$ satisfies $|\la_j|\prod_{k\in P(A)^c}|w_k|<\epsilon$ for each $j\in P(A)$. Then  \eqref{FmQ} implies that for any $z\in T_N$ with $p(z)=w$ and each $j\in P(A)$, $m\in\NN$, we have $|F^m(z)_j|<|z_j|\epsilon^m$. Thus every such $z$ is in $W_T(F)$ and so $w\in p(W_T(F))$. 

Finally, if $P(A)^c=\{j\}$, then $F'(w)=\un\la' w$ with $|\un\la'|<1$, from which it follows that $W_T(F')=\CC^*$ and so $W_T(F)=T_N$. 
\end{proof}

\section{A different construction of toric Kato manifolds}\label{construct}

In this section, we wish to give a different realization of toric Kato manifolds using toric geometry. In particular, we will see that the universal cover of a toric Kato manifold is an open subset of a toric variety of non-finite type. The construction generalizes the one given by Tsuchihashi in \cite{tsu}, as well as  the known toric constructions for Kato surfaces \cite{oda}.

In all that follows, we fix a smooth toric modification $\pi:\hat\CC^n\rightarrow\CC^n$ at $0$ with $\pi\neq\id$ and denote by $\hat{\Sigma}$ the fan of $\hat\CC^n$. We fix a toric chart $\phi_A:\CC^n\rightarrow X_{\tau_A}\subset\hat\CC^n$ for some toric Kato matrix $A\in\GL(n,\ZZ)$ and $\tau_A\in\hat{\Sigma}^{(n)}$. Finally, we fix $\un\la_0\in T_N$ so that $\sigma:=\un\la_0\circ\phi_A$ satisfies $\ov{\sigma(\BB)}\subset\pi^{-1}(\BB)=:\hat\BB$. In particular, $(\pi|_{\hat\BB},\sigma|_{\BB})$ is a toric Kato data. 

Let $\hat{\Sigma}_0:=\hat{\Sigma}\sminus\{\tau_A\}$, and define the infinite fan
\begin{equation}\label{eqn:SigmaA}
\Sigma_A:=\{A^m\tau\mid \tau\in\hat{\Sigma}_0, \ m\in\ZZ\}
\end{equation} 
on which the group $\Gamma_A:=\{A^m\mid m\in\ZZ\}$ acts naturally. Let $X(\Sigma_A,N)$ be the toric variety associated to the fan $\Sigma_A$,  let $\wt{D}_T:=\sum_{\nu\in\Sigma_A^{(1)}}\ov{\orb\nu}$ be the associated toric divisor and let
\begin{equation*} 
D_+:=\sum_{\substack{\nu\in\Sigma_A^{(1)}\\ \nu\not\prec \tau_{P(A)}}}\ov{\orb(\nu)}
\end{equation*}
where we recall that $\tau_{P(A)}=\langle e_k\mid k\in P(A)\rangle$. Define the open set
\begin{equation*}
\wt{X}_{\un\la_0,A}:=\Int(\ov{W_T(F_{\un\la_0,A})}^{X(\Sigma_A,N)})\subset X(\Sigma_A,N).
\end{equation*}
Recall that $W_T(F_{\un\la_0,A})\subset T_N\subset X(\Sigma_A,N)$ is $F_{\un\la_0,A}$-invariant. Also note that the map $F_A$ extends to an automorphism of the toric variety $X(\Sigma_A,N)$ as the map induced by $N\rightarrow N$, $v\mapsto Av$. Similarly, $\un\la_0\in T_N$ extends to an automorphism of $X(\Sigma_A,N)$. Therefore, we have a natural action of $U:=\langle F_{\un\la_0,A}\rangle$ on $\wt{X}_{\la_0,A}$.

Before giving the main theorem of the section, we show the following technical lemma needed in the proof.

\begin{lemma}\label{frontieraOmega}
	Let $V\subset N_\RR$ be an open subset which satisfies $\forall a\in\RR$ with $a\geq 1$, $\forall v\in V$, $av\in V$,  and let $\Omega:=\ord^{-1}(V)\subset T_N$. Let $\Sigma$ be a fan in $N$ and let $X(\Sigma,N)$ be the associated toric variety. Denote by $\ov\Omega^\Sigma$ the closure of $\Omega$ in $X(\Sigma, N)$, and by $\ov V^N$ the closure of $V$ in $N_\RR$. Then for any cone $\tau\in\Sigma^{(k)}$ with $1\leq k\leq n$, one has:
	\begin{enumerate}[(i)]
		\item if $\tau\cap V\neq\emptyset$, then $\orb\tau\subset\ov\Omega^\Sigma$;
		\item if furthermore $V$ is a cone and $(\tau\sminus\{0\})\cap\ov V^N=\emptyset$, then $\orb\tau\cap\ov\Omega^\Sigma =\emptyset$.
	\end{enumerate}
	In particular, if $|\Sigma|\subseteq V\cup\{0\}$, then $\Omega\cup D_T$ is an open subset of $X(\Sigma, N)$, where $D_T=\bigcup_{\tau\in\Sigma^{(1)}}\ov{\orb\tau}$ is the support of the maximal toric divisor of $X(\Sigma,N)$.
\end{lemma}
\begin{proof}
	Let us fix such a cone $\tau\in\Sigma^{(k)}$. In what follows, we let $\Int (S_\tau):=S_\tau\sminus(\tau^\perp\cap M)=\{l\in S_\tau\mid \exists v\in\tau \ \langle l,v\rangle >0\}$.

	Recall that one has a natural map given by extension by $0$:
	\begin{align*}
	j:\orb\tau=\Hom_{gr.}(\tau^{\perp}\cap M,\CC^*)&\rightarrow X_\tau=\Hom_{u.s.g.}(S_\tau,\CC)\\
	j(u)(l)&=\left\{\begin{array}{ll} u(l), & l\in\tau^\perp\cap S_\tau\\
	0, &  \text{otherwise.}\end{array}\right.
	\end{align*}
	Applying the map $\ord$ to the above, $j$ induces a natural injection:
	\begin{equation*}
	j_\RR:\orb\tau_\RR:=\Hom_{gr.}(\tau^\perp\cap M,\RR)\rightarrow X_{\tau,\RR}:=\Hom_{u.s.g.}(S_{\tau},\RR\cup\{\infty\})
	\end{equation*}
	given this time by extension by $\infty$. Moreover, its image is given by 
	\begin{equation}\label{imj}
	j_\RR(\orb\tau_\RR)=\{u\in\Hom_{u.s.g.}(S_{\tau},\RR\cup\{\infty\})\mid \ u^{-1}(\infty)=\Int(S_{\tau})\}.
	\end{equation}
	
	On the other hand, the natural inclusion $\Omega\subset T_N\subset X_\tau$ reads, after applying the map $\ord$, as follows:
	\begin{equation*}
	k:V\rightarrow X_{\tau,\RR}, \ \ k(v)(l)=\langle l,v\rangle, \ \ l\in S_{\tau}.
	\end{equation*}
	
	Since $\orb\tau\subset X_\tau$, one has $\orb\tau\cap \ov\Omega^\Sigma=\orb\tau\cap\ov\Omega^\tau$, where the latter closure is taken in $X_\tau$. Furthermore, since $\ord$ is surjective, this intersection is nonempty precisely when $j_\RR(\orb\tau_\RR)\cap \ov{k(V)}\neq \emptyset$.

	Suppose first that $\tau\cap V\neq\emptyset$ and let $u\in j_{\RR}(\orb\tau_\RR)$. In order to show $(i)$, it suffices to prove that $u\in\ov{k(V)}$. There exists $q\in N_\RR$ so that $u|_{\tau^\perp\cap S_\tau}=\langle\cdot, q\rangle|_{\tau^\perp\cap S_\tau}$. Since $V$ is open, there exists $v\in V$ which is in the relative interior of $\tau$, which is equivalent to $u_0:=\langle\cdot, v\rangle|_{\Int S_\tau}>0$. Furthermore, for $m\geq 1$, since $mv\in V$ and using the properties of $V$, there exists $c_m>m$ big enough so that $v_m:=q+c_mmv=c_m(\frac{1}{c_m}q+mv)\in V$. We find:
	\begin{align*}
	k(v_m)|_{\tau^\perp\cap S_\tau}&=\langle\cdot,q\rangle|_{\tau^\perp\cap S_\tau}=u|_{\tau^\perp\cap S_\tau}\\
	k(v_m)|_{\Int S_\tau}&>\langle \cdot, q\rangle|_{\Int S_\tau}+m^2u_0\rightarrow\infty
	\end{align*}
	hence $\lim_{m\to\infty}k(v_m)=u\in\ov{k(V)}$.
	
	For the second part, suppose that $V$ is a cone and that there exists $u=\lim_{m\to\infty}k(v_m)\in j_\RR(\orb\tau_\RR)$ with $v_m\in\ov V^N$ for each $m\geq 1$. Let us fix the standard scalar product on $N_\RR$ and the orthogonal splitting $N_\RR=N_0\oplus \RR \tau$, so that the map $d:N_0\rightarrow\Hom(\tau^\perp,\RR)$, $q\mapsto\langle\cdot,q\rangle$ becomes a continuous isomorphism. In particular, it follows that $u|_{\tau^\perp\cap S_\tau}=d(q)|_{S_\tau}$ for a unique $q\in N_0$. With respect to the splitting, we write $v_m=v_m^0+v_m^\tau$ with  $v_m^0\in N_0$ and $v_m^\tau\in\RR\tau$. Since $d(v_m^0)|_{\tau^\perp\cap S_\tau}=k(v_m)|_{\tau^\perp\cap S_\tau}\rightarrow d(q)|_{S_\tau}$ as $m\rightarrow\infty$, it follows that $\exists\lim_{m\to\infty}v_m^0=q$. On the other hand, we have:
	\begin{equation*}
	\infty=\lim_{m\to\infty}k(v_m)|_{\Int S_\tau}-\langle\cdot,q\rangle|_{\Int S_\tau}=\lim_{m\to\infty}k(v_m-v_m^0)|_{\Int S_\tau}=\lim_{m\to\infty}\langle\cdot,v_m^\tau\rangle|_{\Int S_\tau}.
	\end{equation*}
	Hence there exists $m_0\geq 1$ so that $\langle\cdot,v_m^\tau\rangle|_{\Int S_\tau}>0$ for all $m\geq m_0$. In particular, $v_m^\tau\in\tau\sminus\{0\}$ for all $m\geq m_0$ and $a_m:=||v_m^\tau||\rightarrow \infty$. Since the sequence $(\frac{1}{a_m}v_m^\tau)_{m\geq m_0}\subset \tau$ is bounded, up to passing to a subsequence, it converges to an element $v_\infty^\tau\in\tau\sminus\{0\}$ with $||v_\infty^\tau||=1$. Now, since $\ov{V}^N$ is a cone,  we infer that $w_m:=\frac{1}{a_m}v_m\in \ov V^N$ and $\lim_{m\to\infty}w_m=v_\infty^\tau\in\ov V^N\cap(\tau\sminus\{0\})$, which concludes $(ii)$.

	Finally, if $|\Sigma|\subseteq V\cup\{0\}$, it follows by $(i)$ that $D_T=\bigcup_{0\neq\tau\in\Sigma}\orb\tau\subset\ov\Omega^\Sigma$, hence $X(\Sigma,N)\sminus(\Omega\cup D_T)=T_N\sminus\Omega$ is closed in $X(\Sigma,N)$,  which concludes the proof.
\end{proof}

We are now ready to prove:

\begin{theorem}\label{secondDescr}
Let $(\pi,\sigma)$ be a toric Kato data with corresponding germ $F_{\un\la_0,A}$. Then the group $U$ acts freely and properly discontinuously on $\wt{X}_{\un\la_0,A}$ and we have a biholomorphism
\begin{equation*}
X(\pi,\sigma)\cong\wt{X}_{\un\la_0,A}/U.
\end{equation*}
Furthermore, if $D$ denotes the divisor of $X(\pi,\sigma)$ induced by the exceptional divisor of $\pi$ and $\wt{D}$ is the preimage of $D$ in the universal cover, then $\wt{D}=D_+$. Finally, we have the following:
\begin{enumerate}[i)]
\item If $P(A)=\emptyset$, then $\wt{X}_{\un\la_0,A}=W_T(F_{\un\la_0,A})\cup \wt{D}_T$. If furthermore $\un\la_0=(1,\ldots, 1)$, then $X(\pi,\sigma)$ is one of the manifolds constructed in \cite{tsu};
\item If $|P(A)|=n-1$, then $\wt{X}_{\un\la_0,A}=X(\Sigma_A,N)$.
\end{enumerate}
\end{theorem}
\begin{proof}
Let us put $\Omega:=W_T(F_{\un\la_0,A})\subseteq T_N$ and $F:=F_{\un\la_0,A}$.  Then $H_\infty=\bigcup_{j\in P(A)^c}\CC_{\{j\}^c}\subset\ov{W_T(F)}^{\CC^n}=\ov{W^s(F)}$. Define the fan:
\begin{equation*}
\Sigma_{A}^+:=\{A^m\tau\mid\tau\in\hat{\Sigma}_0, m\in\NN\}
\end{equation*}
so that $|\Sigma_A^+|\subset C_0=|\Sigma|$ and $X(\Sigma_A^+,N)$ is an open subset of $X(\Sigma_A,N)$. The identity map of $N$ induces a (non-proper) map of toric varieties $\phi:X(\Sigma^+_A, N)\rightarrow \CC^n=X(\Sigma,N)$ and we have 
\begin{equation*}
D_+^+:=\bigcup_{\substack{\nu\in\Sigma_A^{+(1)}\\ \nu\not\prec \tau_{P(A)}}}\ov{\orb\nu}=\phi^{-1}(H_\infty)\subset 
\phi^{-1}(\ov{\Omega}^{\CC^n})=\ov{\Omega}^{X(\Sigma^+_A,N)}
\end{equation*}
from which it follows that $D^+_+\subset\wt{X}_{\un\la_0,A}$. Since for each irreducible component $Q$ of $D_+$ there exists $m\in\ZZ$ so that $F^m(Q)$ is an irreducible component of $D_+^+$, it follows then that also $D_+\subset\wt{X}_{\un\la_0,A}$. Therefore, if $P(A)=\emptyset$ then $D_+=\wt{D}_T$ and we have $\wt{X}_{\un\la_0,A}=W_T(F_{\un\la_0,A})\cup \wt{D}_T$. If $|P(A)|= n-1$ then $W_T(F)=T_N$ by \ref{splitWs} so $\wt{X}_{\un\la_0,A}=X(\Sigma_A,A)$.

Let $\BB^*:=\BB\cap T_N$, $\Omega_0:=\BB^*\sminus F(\BB^*)\subset T_N$ and $W_0:=\ov{\Omega_0}^{X(\Sigma_A,N)}$. Since by the proof of \ref{stableF}, $\Omega_0$ is a fundamental domain for the action of $U$ on $\Omega$, it follows that $W_0$ is a fundamental domain for the action of $U$ on $\wt{X}_{\un\la_0,A}$.

Let $V:=\ord(\BB^*)\subset N_\RR$ and note that it is open and satisfies the condition $v\in V, a\geq 1\Rightarrow a v\in V$. Let $C_0\subset N_\RR$ be the standard cone, so that $V\subset C_0$. By \ref{frontieraOmega} (ii), we have that $\orb\tau\cap\ov{\ord^{-1}(C_0)}^{X(\Sigma_A,N)}=\emptyset$ for any $\tau\in\Sigma_A$ with $\tau\not\subset C_0$, which implies that
\begin{equation*}
\wt{D}_T\cap W_0\subset\bigcup_{\substack{\nu\in\Sigma_A^{(1)}\\ \nu\subset C_0}}\ov{\orb\nu}.
\end{equation*}
Furthermore, since any ray $\nu\in\Sigma_A^{(1)}$ with $\nu\subset\Int C_0$ satisfies $\nu\cap V\neq\emptyset$, \ref{frontieraOmega} (i) implies that $\ov{\orb\nu}\subset\ov{\BB^*}^{X(\Sigma_A,N)}$. In particular, applying $A$, we infer that
\begin{equation*}
\bigcup_{\substack{\nu\in\Sigma_A^{(1)}\\ \nu\subset \Int AC_0}}\ov{\orb\nu}\subset \ov{F_A(\BB^*)}^{X(\Sigma_A,N)}
\end{equation*}
and hence:
\begin{equation*}
\un\la_0\cdot\bigcup_{\substack{\nu\in\Sigma_A^{(1)}\\ \nu\subset \Int AC_0}}\ov{\orb\nu}=\bigcup_{\substack{\nu\in\Sigma_A^{(1)}\\ \nu\subset \Int AC_0}}\ov{\orb\nu}\subset\un\la_0\cdot\ov{F_A(\BB^*)}^{X(\Sigma_A,N)}=\ov{F(\BB^*)}^{X(\Sigma_A,N)}.
\end{equation*}

Therefore $\wt{D}_T\cap W_0=D^{\hat{\Sigma}_0}\cap W_0$, where $D^{\hat{\Sigma}_0}=\bigcup_{\nu\in\hat{\Sigma}_0^{(1)}}\ov{\orb\nu}$, and so $W_0=\ov{\Omega_0}^{X(\hat{\Sigma}_0,N)}$. But $X(\hat{\Sigma}_0,N)=\hat\CC^n\sminus\{\sigma(0)\}$, hence $W_0=\pi^{-1}(\ov\BB)\sminus\sigma(\mathbb \BB)$, which is precisely a fundamental domain for the deck group action on $\widetilde{X(\pi,\sigma)}$. Clearly this identification of $\wt{X}_{\un\la_0,A}$ and of $\widetilde{X(\pi,\sigma)}$ is equivariant, which shows thus that the action of $U$ is free and proper, and that we have the desired isomorphism.

Finally, it is clear from the isomorphism that  $\wt{D}=D_+$.
\end{proof}

For later use, we also show:

\begin{lemma}\label{supportSigma}
The support of $\Sigma_A$ is given by
\begin{equation*}
|\Sigma_A|=\left\{\begin{array}{ll} H(A)\cup\tau_{P(A)}\sminus\left(\RR_{>0}f_A+\tau_{P(A)}\right) , & |P(A)|<n-1\\
H(A)\cup\tau_{P(A)}, &  |P(A)|=n-1\end{array}\right.
\end{equation*}
where for the case $|P(A)|<n-1$, $f_A$ is a Perron vector for $A$ with positive components.
\end{lemma}
\begin{proof}
Since $|\hat{\Sigma}_0|=C_0\sminus\Int(AC_0)$, we find
\begin{equation*}
|\Sigma_A|=\bigcup_{m\in\ZZ}A^m|\hat{\Sigma}_0|=\bigcup_{m\in\ZZ}A^mC_0\sminus\bigcap_{m\in\ZZ}\Int(A^{m}C_0).
\end{equation*}
As in the proof of \ref{permA}, we can suppose without loss of generality that $A_{P(A)^c} =\id$ and that $B:=A_{P(A)^c}$ is a positive matrix.

We first show that 
\begin{equation}\label{eqn:unionfans}
\bigcup_{m\in\ZZ}A^mC_0=H(A)\cup\tau_{P(A)}.
\end{equation}
Indeed, the inclusions $\bigcup_{m\in \mathbb{Z}} A^mC_0\subseteq H(A) \cup \tau_{P(A)}$ and $\tau_{P(A)} \subseteq \bigcup_{m\in \mathbb{Z}} A^mC_0$ are straightforward. Let now $v \in H(A)$.  We shall prove there exists $l>0$ such that $A^lv \in C_0$. For the case  $|P(A)|<n-1$, $f_B=(f_A)_{P(A)^c}$ is the Perron vector of $B$. Let $f_B^*$ be a Perron vector for $\transp{B}$ with $\langle f_B^*,f_B\rangle=1$. For the case $P(A)^c=\{j\}$, we put $f_B=e_j$ and $f^*_B=e_j^*$. If $\al>0$ is the dominant eigenvalue of $A$, or also of $B$, then we have: 
\begin{equation}\label{limB}
\lim_{m \to \infty}\frac{1}{\alpha^{m}}B^m=f_B \cdot \transp{(f^*_B)}
\end{equation}  
which is obvious  in the case $|P(A)|=n-1$ and implied by the Perron-Frobenius theorem for the case $|P(A)|<n-1$. In particular, we have:
\begin{equation*}
\lim_{m \to \infty} \left( \frac{1}{\al^m}A^mv\right)_{P(A)^c}=\langle f_B^*, v\rangle f_B. 
\end{equation*} 

Since $\langle f_B^*,v\rangle> 0$, then there is some $l_0\in\NN$ so that for any $l > l_0$, $(A^lv)_{P(A)^c}$ has positive components. For any $m\in\NN$, put $Q_m:=(a^{(m)}_{jk})_{\substack{j\in P(A)\\ k\in P(A)^c}}$, where $a^{(m)}_{\cdot\cdot}$ denote the components of $A^m$. Then we have:
\begin{equation*}
(A^{l+m}v)_{P(A)} =(A^lv)_{P(A)}+ Q_m(A^lv)_{P(A)^c} 
\end{equation*}
and since  $Q_m$ is a positive matrix, we obtain that for big enough $m$, $A^{l+m}v\in C_0$.

In order to compute the intersection, we distinguish between the two cases. Suppose first that $|P(A)|<n-1$. Then we wish to show that
\begin{equation}\label{eqn:coneinter}
\bigcap_{m\in\ZZ}\mathrm{Int}(A^mC_0)=\RR_{>0}f_A + \tau_{P(A)}.
\end{equation}
It is clear that $\RR_{>0}f_A+\tau_{P(A)}$ is invariant by $A$ and is contained in $\mathrm{Int} C_0$. Moreover, for any vector $v \in C_0 \sminus \tau_{P(A)} \subset H(A)$, we have that
$$
\lim_{m \to +\infty} \frac{A^m v}{\alpha^m} = f_A \langle f_A^*, v \rangle = kf_A
$$
for a suitable $k > 0$. It follows that the cones $A^mC_0=\langle A^me_1,\ldots, A^me_n\rangle$ tend to $\RR_{>0}f_A + \tau_{P(A)}$ when $m \to +\infty$, and so we have:
\begin{equation*}
\RR_{>0}f_A+\tau_{P(A)}\subseteq\bigcap_{m\in\ZZ}A^m\mathrm{Int}C_0\subseteq \bigcap_{m\geq 0}A^mC_0=\RR_{>0}f_A+\tau_{P(A)}
\end{equation*}
which proves \eqref{eqn:coneinter}. 

Suppose now that $P(A)^c=\{j\}$. Then we have:
\begin{equation*}
 \bigcap_{m\in\ZZ}A^mC_0=\tau_{P(A)}.
\end{equation*}
Indeed, clearly $\tau_{P(A)}\subset A^m C_0$ for any $m\in \ZZ$. Conversely, let $v\in\bigcap_{m\in\ZZ}A^mC_0$, so that for each $m\in\ZZ$ there exist $v_m\in \RR_{P(A)}$ and $a_m\in\RR_{\geq 0}$ with 
\begin{equation*}
v=A^m(v_m+a_me_j)=v_m+a_mQ_me_j+a_me_j.
\end{equation*}
In particular, we find that $v_j=a_me_j$ so that $a_m=a_0$ for any $m\in\ZZ$. If $a_0>0$, then $||a_mQ_me_j||$ tends to $\infty$ with $m$, so $||v||$ is unbounded, which is absurd. Thus $a_0=0$ and so $v\in\tau_{P(A)}$. 

Finally, we have
\begin{equation*}
\bigcap_{m\in\ZZ}\mathrm{Int}(A^mC_0)\subseteq \Int C_0\cap \bigcap_{m\in\ZZ}A^mC_0=\emptyset.
\end{equation*}
which concludes the proof.
\end{proof}

\begin{example}
Consider the toric Kato data $(\pi, \sigma)$, where $\pi$ is associated to the fan depicted in the left side of Figure \ref{fig:example1}, and $\sigma = \un\la_0 \phi_A$ with $A=\begin{pmatrix}
1 & 1 & 1 \\
0 & 2 & 1 \\
0 & 1 & 1
\end{pmatrix}
$.
In this case, $P(A)=\{1\}$, and a direct computation shows that the eigenvalues of $A$ are $1, \xi^2, \xi^{-2}$, where $\xi=\frac{1+\sqrt{5}}{2}$ is the golden ratio.
The right side of Figure \ref{fig:example1} shows the fan $\hat{\Sigma}_2$ associated to the toric modification $\pi^2$ given by \ref{powerKD}.

A Perron vector $f_A$ is given by $(\xi, \xi, 1)$, while the dual vector $f_A^*$ is a positive multiple of $(0, \xi, 1)$.
In this case $H(A)=\{(x,y,z) \in \nR^3 |\ \xi y + z > 0 \}$, and $\abs{\Sigma_A}=H(A) \cup \nR_{\geq 0} (1,0,0) \sminus \{(s+t\xi,t\xi,t)\ |\ s \in \nR_{\geq 0}, t \in \nR_{>0}\}$. 

\begin{tiny}
\begin{figure}[h]
\centering
\begin{minipage}[htbp]{\columnwidth}
\def\svgwidth{0.48\columnwidth}
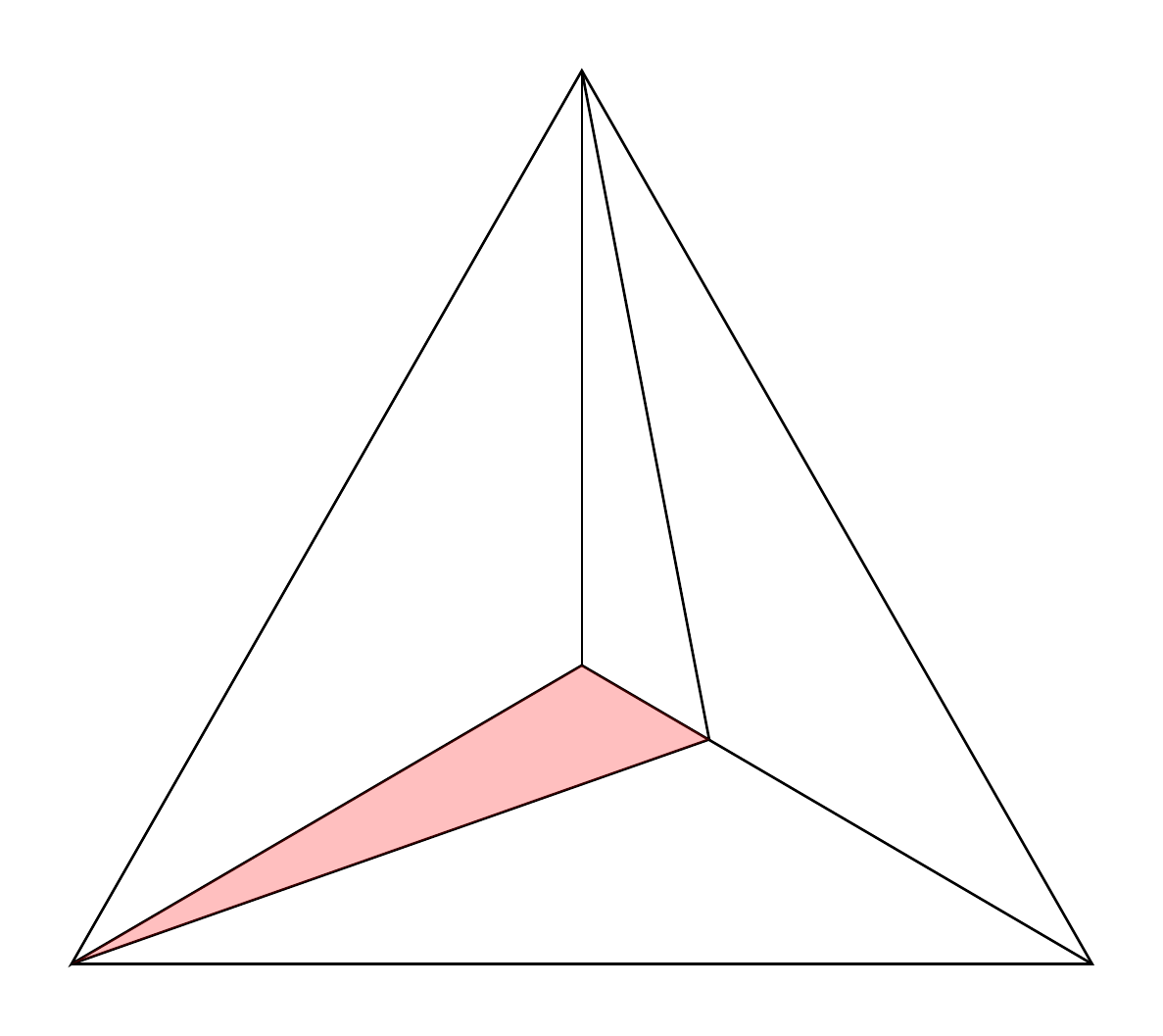
\hfill
\def\svgwidth{0.48\columnwidth}
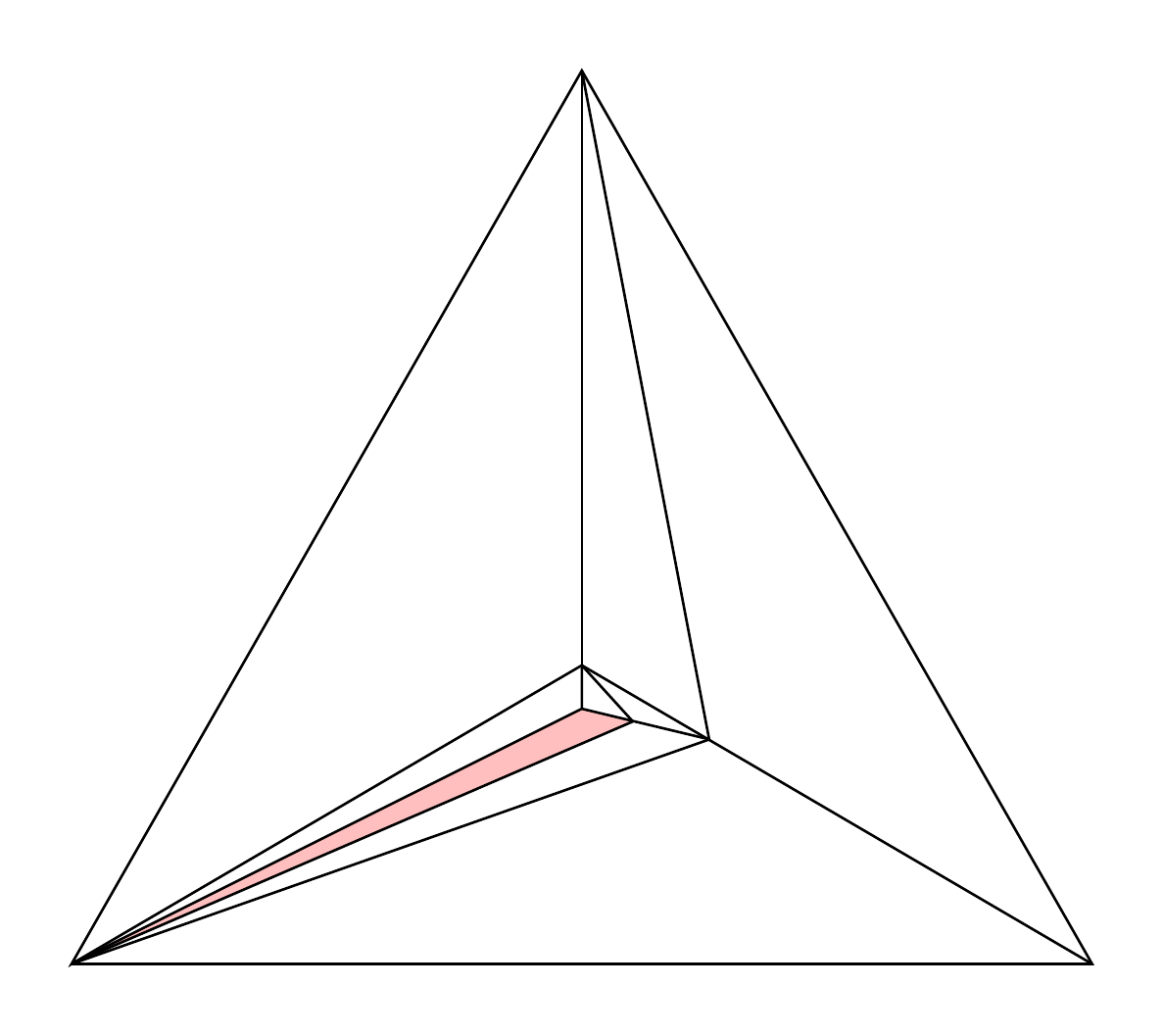
\end{minipage}
\caption{Toric Kato data with $P(A)=\{1\}$.}
\label{fig:example1}
\end{figure}
\end{tiny}
\end{example}

\section{Geometric properties of toric Kato manifolds}\label{geomprop}

In complex dimension $n=2$, the germ $F_{\un\la_0,A}$ determines the toric Kato surface (if we ask it to be minimal). Toric Kato surfaces with $P(A)=\emptyset$ are precisely the properly blown-up Inoue-Hirzebruch surfaces, in the sense of \cite{fp}. As we will see later,  $X(F_{\un\la_0,A})\cong X(F_A)$ for any compatible $\un\la_0\in T_N$. In the case $|P(A)|=1$, toric Kato surfaces are precisely the properly blown-up parabolic Inoue surfaces, and for fixed $A$, they are parametrized by $\DD^*$. Finally, the case $|P(A)|=2$ gives finite quotients of diagonal Hopf surfaces. In higher dimension $n$, it is still true that the case $|P(A)|=n$ gives finite quotients of diagonal Hopf manifolds, while the dichotomy persists between the cases $|P(A)|=n-1$ and $|P(A)|<n-1$. For this reason, we propose the following definition:

\begin{definition}
We call a toric Kato manifold $X$ of matrix $A\in\GL(n,\ZZ)$  of \textit{parabolic type} if $|P(A)|=n-1$. If $|P(A)|<n-1$, then we call $X$ of \textit{hyperbolic type}. 
\end{definition}

In analogy to Nakamura's result in dimension $n=2$ \cite{naka84}, we can detect the type of a toric Kato manifold by looking at its curves. This shows in particular that the above definition makes sense, and does not depend on the chosen Kato data. In order to make this precise, we first introduce one more definition. 

We will call a compact complex subspace $Y$ of a toric Kato manifold $X$ $T$-invariant if its preimage in the universal cover is invariant under the $\TT^n$- action, or equivalently, under the (local) $T_N$-action. All such irreducible subspaces correspond to cones in $\Sigma_A$, and are of two types: either they are immersed toric Kato manifolds, if the corresponding cone is $A$-invariant, or they are complete toric algebraic varieties, if the corresponding cone is not $A$-invariant. In particular, we have the following description of $T$-invariant curves:

\begin{theorem}\label{hiperpara}
Let $X$ be a toric Kato manifold. Then we have the following characterization in terms of $T$-invariant curves:
\begin{enumerate}
\item $X$ is a finite quotient of a diagonal Hopf manifold if and only if all its $T$-invariant curves are elliptic; 
\item $X$ is of hyperbolic type if and only if all $T$-invariant curves are rational;
\item $X$ is of parabolic type if and only if $X$ contains a unique $T$-invariant elliptic curve $E$, and at least one rational $T$-invariant curve. In this case, $E$ is smooth and all other $T$-invariant curves are rational.  Moreover, the intersection number $E\cdot D=0$.
\end{enumerate}
\end{theorem}
\begin{proof}
Let $A\in\GL(n,\ZZ)$ denote a Kato matrix of $X$ and suppose, without loss of generality, that $A_{P(A)}=\id$. Any $T$-invariant curve $C$ of $X$ is given by a cone $\tau\in\Sigma_A^{(n-1)}$, either as $\orb\tau/\Gamma$, in case $\tau$ is $\Gamma$-invariant, or as the image of $\ov{\orb\tau}\cap\wt{X}$. Moreover, $C$ is an elliptic curve in the first case, and toric, thus rational, in the second case.

If $|P(A)|=n$, then $\Sigma_A=\Sigma\sminus\{C_0\}$, so it is clear that $X$ contains $n$ $T$-invariant elliptic curves and no $T$-invariant rational curve. If $|P(A)|<n-1$, then there is no $A$-invariant cone $\tau\in\Sigma_A^{(n-1)}$, so all $T$-invariant curves of $X$ are rational. 

Suppose now that $P(A)^c=\{j\}$. Then the only $A$-invariant cone $\tau\in\Sigma_A^{(n-1)}$ is $\tau_{P(A)}$, so $E:=\orb\tau_{P(A)}/\Gamma\cong\CC^*/\la_j$ is the only elliptic curve in this last case, and is clearly smooth. Moreover, it is clear that there exits at least one $\tau_{P(A)}\neq\tau\in\Sigma_A^{(n-1)}$, giving rise to a  $T$-invariant rational curve.

Recall that by \ref{secondDescr},  $D$ is the image in $X$ of 
\begin{equation*}
D_+=\sum_{\substack{\nu\in\Sigma_A^{(1)}\\ \nu\not\subset \tau_{P(A)}}}\ov{\orb\nu}\cap \wt{X}.
\end{equation*}
Thus, in order to see that $E\cdot D=0$, it is enough to see that $\orb\tau_{P(A)}$ does not intersect $\ov{\orb\nu}$ for any $\nu\in\Sigma_A^{(1)}$ with $\nu\not\subset \tau_{P(A)}$. But this is equivalent to say that any such $\nu$ is not a face of $\tau_{P(A)}$, which is obvious.
\end{proof}

Next we want to give an explicit description of the divisors defined by 
\begin{align*}
\wt{D}^X_T:&=\wt{D}_T\cap \wt{X}=\wt{X}_{\un\la_0,A}\sminus W_T(F_{\un\la_0,A}), \\ 
D_T&=\wt{D}^X_T/_\Gamma.
\end{align*}
It is clear that $\wt{D}^X_T=D_++\sum_{k\in P(A)}X_k$, where $X_k:=\ov{\orb\langle e_k\rangle}^{X(\Sigma_A,N)}\cap\wt{X}_{\un\la_0,A}$. We will call $D_T$ the \textit{toric divisor} of $X$.  It satisfies $q^* D_T=\wt{D}_T^X$, where $q:\wt{X}\rightarrow X$ is the universal cover map. 

Suppose first that the permutation $s$ given by \ref{permA} corresponding to $A$ is the identity. Fix $k\in P(A)$, let $\hat\CC_{\{k\}^c}=\ov{\orb\langle e_k\rangle}\subset\hat\CC^n$ and let $\pi_k:=\pi|_{\hat\CC_{\{k\}^c}}:\hat\CC_{\{k\}^c}\rightarrow \CC_{\{k\}^c}$. Then the proof of \ref{subKatomatrix} shows that the toric chart $\phi_A$ of $\hat\CC^n$ induces a toric chart $\phi_{A_{\{k\}^c}}$ of $\hat\CC_{\{k\}^c}$. Let $\hat{\Sigma}_{k}$ be the fan of $\hat\CC_{\{k\}^c}$ and let $\Sigma_{A_{\{k\}^c}}:=\{A_{\{k\}^c}^m\tau\mid \tau\in\hat{\Sigma}_{k}\sminus\{\tau_{A_{\{k\}^c}}\}, m\in\ZZ\}$. Then we find that 
\begin{equation*}
\ov{\orb\langle e_k\rangle}^{X(\Sigma_A,N)}\cong X(\Sigma_{A_{\{k\}^c}},N_{\{k\}^c})
\end{equation*}
and it is easy to check that $F_{\un\la_0,A}|_{\ov{\orb\langle e_k\rangle}}=F_{(\un\la_0)_{\{k\}^c},A_{\{k\}^c}}$. In particular, we have:
\begin{equation*}
X_k\cong \ov{p_{\{k\}^c}(W_T(F_{\un\la_0,A}))}^{\ov{\orb\langle e_k\rangle}}=\ov{W_T(F_{(\un\la_0)_{\{k\}^c},A_{\{k\}^c}})}^{X(\Sigma_{A_{\{k\}^c}},N)}.
\end{equation*}

We have thus found:
\begin{lemma}
If $A$ is a toric Kato matrix with $s=\id$, and $k\in P(A)$, then $X_k$ is $\Gamma$-equivariantly biholomorphic to $\wt{X}_{(\un\la_0)_{\{k\}^c},A_{\{k\}^c}}$, where $\Gamma$ acts by $F_{\un\la_0,A}$ on the first manifold and by $F_{(\un\la_0)_{\{k\}^c},A_{\{k\}^c}}$ on the second one. In particular, $X_k/\Gamma\cong X(\pi_k,\sigma_k)$, where $\sigma_k=(\un\la_0)_{\{k\}^c}\phi_{A_{\{k\}^c}}$. 
\end{lemma}

In the general case, i.e. when $s\neq\id$, $X_k$ need not be invariant by $F:=F_{\un\la_0,A}$, however $F$ still acts on $\sum_{k\in P(A)} X_k$. In this case, there exists $d\in\NN$ with $s^d=\id$. Then it is easy to see that $W_T(F)=W_T(F^d)$ and that $F^d=F_{\un\la',A^d}$ where $\un\la'=\un\la_0^{\id+\ldots+A^{d-1}}$. Consider the toric Kato data $(\pi^d,\sigma^d)$ given by \ref{powerKD}, with corresponding germ $F^d$. Then clearly $\Sigma_{A^d}=\Sigma_A$. The above argument thus implies that $X_k\cong\wt{X}_{(\un\la')_{\{k\}^c},A^d_{\{k\}^c}}$. To each cycle $J\subset P(A)$ of $s$ then corresponds a compact hypersurface which autointersects transversely
\begin{equation*}
Y_{J}=\left(\bigcup_{k\in J} X_k\right)/\Gamma\subset X(\pi,\sigma)
\end{equation*}
such that the divisor  $D_J:=Y_J$ on $X(\pi,\sigma)$ satisfies $q^*D_J=\sum_{k\in J}X_k$.
Finally, we have found:

\begin{proposition}\label{divD'}
Let $(\pi,\sigma)$ be a toric Kato data with germ $F=F_{\un\la_0,A}$ and let $s$ be the permutation corresponding to $A$. Let $C:=\{J\subset P(A)\mid J$ is a cycle of $s\}$. Then we have:
\begin{equation*}
D_T=D+\sum_{J\in C}Y_J.
\end{equation*}
Moreover, for each $J\in C$, $Y_J$ is a compact (possibly singular) hypersurface whose normalization is a toric Kato manifold. 
\end{proposition}

We also want to establish when the divisors $D,\wt{D}, D_T, \wt{D}^X_T$ are connected. In complex dimension $2$, the situation is well known by the works of Nakamura \cite{naka84} and Dloussky \cite{dl84}, and can be easily read off from the fan $\Sigma_0$. Namely, in the case of parabolic Inoue surfaces, i.e. when $|P(A)|=1$, if $P(A)=\{j\}$, then $\ov{\orb\langle e_j\rangle}$ is fixed by $A$ and disconnected from $\wt{D}$, which itself is connected. It follows that both $\wt{D}^X_T$ and $D_T$ have two connected components, respectively $\wt{D}$ and $\ov{\orb\langle e_j\rangle}$ for $\wt{D}_T$, and $D$ and the elliptic curve $\ov{\orb\langle e_j\rangle}/\Gamma$ for $D_T$. In the case of Inoue-Hirzebruch surfaces, $\wt{D}^X_T=\wt{D}$ has two connected components, given by the two connected components of $|\Sigma_0|\sminus\{0\}$. When $\det A=1$, $A$ fixes the components, so that $D$ has also two connected components. When $\det A=-1$, $A$ interchanges the components, and in this case $D$ is connected. In higher dimension, things are much simpler:

\begin{lemma}\label{conndiv}
Let $X=X(\pi,\sigma)$ be a toric Kato manifold of complex dimension $n\geq 3$. Then the divisors $\wt{D}^X_T,\wt{D}, D_T$ and $D$ are connected. 
\end{lemma}
\begin{proof}
Given a divisor $F$ on a toric variety $X(\Sigma,N)$, $F=\sum_{j\in I}\ov{\orb\nu_j}$, let $G_{F,\Sigma}=(V,E)$ be the graph whose vertexes are given by $V=\{\nu_j\mid j\in I\}$ and edges given by $E=\{(\nu_i,\nu_j)\mid i,j\in I, \nu_i\oplus\nu_j\in\Sigma^{(2)}\}$. It is clear that $F$ is connected if and only if the graph $G_{F,\Sigma}$ is connected. 

Let $A$ be a Kato matrix of $X$. Consider first $\hat\CC^n=X(\hat{\Sigma},N)$, the exceptional divisor $E$ of $\pi$, and the divisors $D_1=E+\sum_{j=1}^n\ov{\orb\langle e_j\rangle}$ and $D_2=E+\sum_{j\notin P(A)}\ov{\orb\langle e_j\rangle}$. Since $E$ is connected and intersects each component $\ov{\orb\langle e_j\rangle}$, $1\leq j\leq n$, also $D_1$ and $D_2$ are connected, so their graphs $G_{D_1,\hat{\Sigma}}$, $G_{D_2,\hat{\Sigma}}$ are connected. Next, consider $\hat\CC^n\sminus\{\sigma(0)\}=X(\hat{\Sigma}_0,N)\subset\hat\CC^n$ and $D'_j=D_j\cap X(\hat{\Sigma}_0,N)$, $j=1,2$. As $n\geq 3$, we have $G_{D_j,\hat{\Sigma}}=G_{D'_j,\hat{\Sigma}_0}$, so $D'_j$ are also connected, $j=1,2$. From this we infer that $\wt{D}_T=\sum_{m\in\ZZ}\gamma^m D'_1$ and $\wt{D}=\sum_{m\in\ZZ}\gamma^m D'_2$ are connected on $X(\hat{\Sigma},N)$, hence also their respective restrictions to $\wt{X}$ are connected and their images in $X$ are connected. This concludes the proof.  
\end{proof}

We also have the following geometric description of Kato manifolds with $P(A)\neq\emptyset$:

\begin{proposition}
Let $X_A:=X(\pi,\sigma)$ be a toric Kato manifold with $\sigma=\un\la_0\sigma_A$. Let $B=A_{P(A)^c}$ and let $(\pi',\sigma'=\un\la'\sigma_B)$ be the naturally induced Kato data from $(\pi,\sigma)$ and $B$. Let  $j:X_B:=X(\pi',\sigma')\rightarrow X_A$ be the natural embedding. We denote by $D_A$ the divisor of $X_A$ induced by the exceptional divisor of $\pi$, and by $D_B=j^*D_A$ the corresponding divisor of $X_B$. Then $X_A\sminus D_A$ has a natural structure  of a holomorphic vector bundle $V$ of rank $|P(A)|$ over  $X_B\sminus D_B$. Moreover, up to passing to a finite unramified cover of $X_A$, we have $V=\oplus_{k\in P(A)} L_k$, where each $L_k\in\Pic(X_B\sminus D_B)$. In particular, if $P(A)^c=\{j\}$, then $X_B=X_B\sminus D_B$ is the unique $T$-invariant elliptic curve $E=\CC^*/\la_j$ of $X_A$, and for each $k\in P(A)$, $L_k$ has negative Chern class determined by $A$.
\end{proposition}
\begin{proof}
Putting
\begin{equation*}
\Sigma':=\Sigma_A\sminus \bigcup_{\substack{\nu\in\Sigma_A^{(1)}\\ \nu\not\subset\tau_{P(A)}}}St_\nu=\{\tau\in\Sigma_A\mid\tau\prec\tau_{P(A)}\}
\end{equation*}
we have that 
\begin{equation*}
X(\Sigma_A, N)\sminus D_{A,+}=X(\Sigma',N)\cong X(\Sigma', N_{P(A)})\times  X(\{0\},N_{P(A)^c}).
\end{equation*}
Moreover, $X(\{0\},N_{P(A)^c})=X(\Sigma_B,N_{P(A)^c})\sminus  D_{B,+}$, and the natural projection $u_\ZZ:N\rightarrow N_{P(A)^c}$ induces a toric morphism $u:X(\Sigma',N)\rightarrow X(\{0\},N_{P(A)^c})$ satisfying $u(\gamma z)=\gamma'(u(z))$ for any $z\in X(\Sigma',N)$, where $\gamma$ is the positive generator of $\Gamma$ acting on $X(\Sigma_A,N)$ and $\gamma'$ is the positive generator of $\Gamma$ acting on $X(\Sigma_B,N_{P(A)^c})$. Since moreover $W_T(F_A)=T_{P(A)}\times W_T(F_B)$ by \ref{splitWs}, it follows that we have an induced vector bundle:
\begin{equation*}
u|_V:V:=X_A\sminus D_A\rightarrow X_B\sminus D_B.
\end{equation*}

With respect to the standard basis $e_1,\ldots e_n$ of $N=\ZZ^n$, we have a natural splitting 
\begin{equation*}
X(\Sigma',N_{P(A)})=\oplus_{k\in P(A)}\CC e_k.
\end{equation*} 
Up to taking a finite cover of $X_A$, which corresponds to taking a positive power of $A$, suppose that $A_{P(A)}=\id$. Then the action of $\Gamma$ preserves each $\CC e_k\times W_T(F_B)$ and $\gamma$ acts on it by 
\begin{equation*}
\gamma(z_k,w)=(\la_k w^{l_k}z_k,\gamma'(w)), \ \ (z_k,w)\in\CC e_k\times W_T(F_B),
\end{equation*} 
where $l_k\in\NN^{|P(A)^c|}$ is the line $(A_{ks})_{s\in P(A)^c}$. It follows that $\CC e_k$ descends to a line bundle $L_k$ over $X_B\sminus D_B$ and that $V=\oplus_{k\in P(A)} L_k$.

 In particular, when $P(A)^c=\{j\}$, each line bundle $L_k$ is given by the multiplier $\epsilon_{k,\gamma}(t)=\la_k t^{A_{kj}}$, $t\in\CC^*=W_T(F_B)$. From this we infer that
\begin{equation*}
c_1(L_k)=-A_{kj}\in H^2(E,\ZZ)=\ZZ.
\end{equation*}

\end{proof}

Finally, we also have:

\begin{proposition}\label{Kodairadim}
Let $X$ be a toric Kato manifold with toric divisor $D_T$ and matrix $A$, and let $L_{\det A}$ be the flat line bundle on $X$ determined by $\rho:\Gamma\rightarrow \CC^*$, $\rho(\gamma)=(-1)^{\det A}$. Then we have
\begin{equation*}
K_X\otimes L_{\det A}=\Oo_X(-D_T).
\end{equation*}
In particular, $X$ has negative Kodaira dimension.
\end{proposition}
\begin{proof} 
Let $e^*_1\ldots e^*_n$ be the standard basis of $M$ and consider the meromorphic $n$-form on $\wt{X}$: 
	\begin{equation*}
	\omega := \frac{d \wt{e}^*_1}{\wt{e}^*_1} \wedge \ldots \wedge \frac{d \wt{e}^*_n}{\wt{e}^*_n}.
	\end{equation*}
Then \cite[Proposition at page~85]{fulton} shows that $\mathrm{div}(\omega)=-\wt{D}_T$. Since moreover $\gamma^*\omega=(-1)^{\det A}\omega$, the conclusion follows.  	
\end{proof}

\section{Betti numbers}\label{betti}

In this section, we fix a toric Kato manifold $X=X(\pi,\sigma)$ of matrix $A$, with $\pi:\hat\CC^n\rightarrow \CC^n$. Let $D$ be the divisor of $X$ induced by the exceptional divisor of $\pi$. We also let $\hat{\Sigma}$ denote the fan of $\hat\CC^n$, we put $\hat{\CC\PP^n}=\hat\CC^n\cup H$, where $H$ is the divisor at infinity, and we let $\Sigma_\PP$ be the corresponding fan of $\hat{\CC\PP^n}$. Let us put $a_j:=|\hat{\Sigma}^{(j)}|$, for each $0\leq j\leq n$. Then we can express the Betti numbers of $X$ purely in terms of the combinatorial data $a_0,\ldots, a_n$.

\begin{theorem}\label{bettith}
The $n$-dimensional toric Kato manifold $X$ has the following Betti numbers:
\begin{gather*}
b_0(X)=b_1(X)=b_{2n-1}(X)=b_{2n}(X)=1\\
b_{2j+1}(X)=0, \ \ 1\leq j\leq n-2\\
b_{2j}(X)=-1+\sum_{s=j}^n(-1)^{s-j}\binom{s}{j}\left(a_{n-s}+\binom{n}{s+1}\right), \ \ 1\leq j\leq n-1.
\end{gather*}
In particular, we have 
\begin{equation*}
b_2(X)=\sharp D, \ \ \chi(X)=a_n-1
\end{equation*} 
where $\sharp$ denotes the number of irreducible components, and $\chi(X)$ denotes the Euler characteristic of $X$.
\end{theorem}
\begin{proof}
Let us suppose that $n>2$, since the statement for $n=2$ is well known. In this proof, we use the notation $H^\bullet(M)$ to denote the cohomology of $M$ with coefficients in $\CC$.

We will first show that $b_j(X)=b_j(\hat\CC^n)$ for $2\leq j\leq 2n-2$. By \cite{kato}, $X=X(\pi,\sigma)$ is diffeomorphic to a Kato manifold $X'=X(\pi,\sigma')$, so that $\sigma'(0)$ no longer meets the exceptional divisor of $\pi$. In particular, there exists a Hopf manifold $Y\cong \Ss^{2n-1}\times\Ss^1$, $p\in Y$ and a modification $\mu:X'\rightarrow Y$ at $p$ induced by $\pi$. Let $B\subset Y$ be a ball centered at $p$, so that $V:=\mu^{-1}(B)\subset X'$ is a neighborhood of the exceptional divisor $F$ of $\mu$ which is homotopy equivalent to $\hat\CC^n$, and let $U:=X'\sminus F\cong Y\sminus \{p\}$. Using the fact that $U\cap V\cong \Ss^{2n-1}$ and that $b_j(U)=b_j(Y)$ for $0\leq j\leq 2n-1$, the Mayer-Vietoris sequence in cohomology for $(U,V)$ gives $H^j(X)=H^j(X')=H^j(\hat\CC^n)$ for $2\leq j\leq 2n-2$. 

Next, we wish to compare the cohomology of $\hat\CC^n$ with that of $\hat{\CC\PP^n}$. Let $U$ be a neighborhood of $H$ in $\hat{\CC\PP^n}$ biholomorphic to a $\DD$-bundle over $H$ inside $\mathcal N_{\hat{\CC\PP^n}}H$, let $V:=\hat\CC^n\subset\hat{\CC\PP^n}$ and let $U^*:=U\cap V$, which is a $\DD^*$-bundle over $H$. We have $H^\bullet(U)=H^\bullet(H)$, while $H^j(U^*)=H^{j}(H)\oplus H^{j-1}(H)=\CC$ for $j\leq 2n-2$. Using this, the Mayer-Vietoris sequence in cohomology then gives:
\begin{equation*}
\xymatrix{ \ar[r] &H^j(\hat{\CC\PP^n})\ar[r] & H^j(H)\oplus H^j(\hat\CC^n)\ar[r]^-{r^*_U-r^*_V} &\CC\ar[r] &H^{j+1}(\hat{\CC\PP^n})\ar[r] &
}
\end{equation*}
where $r_U$, $r_V$ are the corresponding restriction maps, so that $r^*_V=0$ and $r^*_U$ is injective. From this we infer 
\begin{equation}\label{BettiCP}
b_j(\hat\CC^n)=b_j(\hat{\CC\PP^n}), \ \ j \text{ odd }, \ \ b_j(\hat\CC^n)=b_j(\hat{\CC\PP^n})-1, \ \ j \text{ even}.
\end{equation}

Finally, since $\hat{\CC\PP^n}=X(\Sigma_\PP,N)$ is a complete toric algebraic variety, it satisfies \cite[Theorem~10.8]{da}
\begin{equation*}
b_{2j+1}(\hat{\CC\PP^n})=0, \ \ b_{2j}(\hat{\CC\PP^n})=\sum_{s=j}^n(-1)^{s-j}\binom{s}{j}\hat{a}_{n-s}, \ \ j\in\NN
\end{equation*}
where $\hat{a}_j:=|\Sigma_\PP^{(j)}|=a_j+\binom{n}{j-1}$, $j\in \{1,\ldots, n\}$, $\hat{a}_0=a_0$. Together with  \eqref{BettiCP}, this gives the desired formula for $b_\bullet(X)$. 

In particular, since $a_1=\sharp D + n$,  we find
\begin{equation*}
b_2(X)=b_{2n-2}(X)=-1+(a_1+1)-(n\cdot a_0)=\sharp D. 
\end{equation*}
Furthermore, we have:
\begin{gather*}
\chi(\hat{\CC\PP^n})=\sum_{j=0}^n\sum_{s=j}^n\hat{a}_{n-s}(-1)^{s-j}\binom{s}{j}=\sum_{s=0}^n\hat{a}_{n-s}\sum_{j=0}^s(-1)^{s-j}\binom{s}{j}=\hat{a}_n=a_n+n\\
\chi(X)=\sum_{j=1}^{n-1}b_{2j}(X)=\sum_{j=1}^{n-1}(b_{2j}(\hat{\CC\PP^n})-1)=a_n+n-2-(n-1)=a_n-1
\end{gather*}
which concludes the proof.
\end{proof}

\section{Toric degenerations}\label{degen}

\subsection{Nakamura degenerations}

In this section, we describe how the Nakamura degenerations \cite{naka83} of toric Kato surfaces, generalizing constructions of Miyake-Oda \cite{oda}, can be generalized to any dimension.  Let us fix a toric Kato manifold $X=X(\pi,\sigma)$ with germ $F=\un\la F_A$. Recall that $\wt{X}$ can be seen as a partial compactification of $\Omega_F=W_T(F)$ inside the toric manifold $X(\Sigma_A,N)$.  

Also recall that we denoted by $\hat{\Sigma}_0$ the fan of the toric manifold $\hat\CC^n\sminus \{\sigma(0)\}$, which is then a fundamental domain of the action of $\Gamma$ on $\Sigma_A$. Here $\Gamma$ is the deck group of $\wt{X}\rightarrow X$, generated by $\gamma$. Let $\hat{\CC\PP^n}$ be $\hat\CC^n$ with the hyperplane at infinity $H$, and let $\widetilde{\CC\PP^n}$ be the blowup of $\hat{\CC\PP^n}$ along $\sigma(0)$, with exceptional divisor $E$. Then $\widetilde{\CC\PP^n}=X(\Sigma_c,N)$ is a toric manifold, with fan $\Sigma_c$ given by the completion of $\hat\Sigma_0$ with the $n$-dimensional cones:
\begin{gather*}
\tau_j=\langle Ac, Ae_1,\ldots, \widehat{Ae_j}, \ldots Ae_n\rangle , \ \ 1\leq j\leq n\\
\tau'_j=\langle -c, e_1,\ldots, \widehat{e_j},\ldots, e_n\rangle, \ \ 1\leq j\leq n\\
c:=\sum_{j=1}^ne_j.
\end{gather*}

We now wish to construct a manifold $\wt{\mathcal{X}}_N$ as a subset in a toric variety, together with an action of $\Gamma$, giving rise to the Nakamura degeneration. Consider then $\wt{N}=N\oplus \ZZ e$, and define an action $\wt{A}$ on $\wt{N}$ by putting $\wt{A}|_N=A$ and $\wt{A}e=e+Ac$. Put $\nu_k:=\wt{A}^ke$, $k\in\ZZ$, and define regular fans $\wt{\Sigma}_0,\wt{\Sigma}$ in $\wt{N}$ as follows:
\begin{gather*}
\wt{\Sigma}_0:=\{\wt{\tau}:=\tau\oplus\langle e\rangle \mid \tau\in\hat\Sigma_0\}\cup\{\wt{\tau}':=\tau\oplus\langle e,\nu_{-1}\rangle\mid \tau\in\hat\Sigma_0, \tau\subset \del C_0\}\\
\wt{\Sigma}:=\bigcup_{k\in\ZZ}\wt{A}^k\wt{\Sigma}_0.
\end{gather*}

Let us also put $\wt{\un\la}=(\un\la,1)\in T_{\wt{N}}$ and let $\wt{\gamma}$ be the automorphism of $X(\wt{\Sigma},\wt{N})$ given by $\wt{\un\la}\wt{A}$, where we also denote by $\wt{A}$ the natural action induced by $\wt{A}$ on the toric variety. Note that for $z\in T_{\wt{N}}$, $\wt{\gamma}(z)=\wt{F}(z):=\wt{\un\la} z^{\wt{A}}$.  

Let $p_\ZZ:\wt{N}\rightarrow \ZZ e$ be the natural projection, and let $\Sigma_1=\{0,\langle e\rangle \}$ be the fan of $\CC$ in $\ZZ e$. Then $p_\ZZ$ is a map of fans, and so defines naturally a map of toric varieties $p:X(\wt{\Sigma},\wt{N})\rightarrow X(\Sigma_1,\ZZ e)$. We have 
\begin{gather*}
p^{-1}(\CC^*)=p^{-1}(\orb 0)=\bigcup_{\tau\in \wt{\Sigma}, p_\ZZ(\tau)=0}\orb\tau=X(\wt{\Sigma}\cap N, \wt{N})=X(\Sigma_A,N)\times\CC^*
\end{gather*}
and for $(x,t)\in X(\Sigma_A,N)\times\CC^*$ we have $\wt{\gamma}(x,t)=(\un t^A\gamma(x),t)$, where $\un t:=(t,\ldots, t)\in T_N$. Moreover
\begin{gather*}
p^{-1}(0)=p^{-1}(\ov{\orb e})=\bigcup_{\tau\in\wt{\Sigma}, p_\ZZ(\tau)=\langle e\rangle} \overline{\orb\tau}=\bigcup_{m\in\ZZ}\overline{\orb\nu_m}
\end{gather*}
and $\wt{\gamma}(\ov{\orb\nu_m})=\ov{\orb \nu_{m+1}}$. Let us note that $\ov{\orb e}$ in $X(\tilde\Sigma,\tilde N)$ is the toric variety $X(\widehat{St}_{e}, N_{e})$, where $N_{e}:=\wt{N}/\ZZ e$ and
\begin{equation*}
St_{e}=\{\tau\in\wt{\Sigma} \mid \langle e\rangle\prec \tau\}, \ \ \widehat{St}_e=\{\tau \mod \ZZ e\mid \tau\in St_e\}.
\end{equation*}
Under the natural isomorphism $N_e\cong N$ given by the projection, we have $\widehat{St}_e\cong \Sigma_c$. Thus $\ov{\orb \nu_0}$ is equivariantly isomorphic to $\widetilde{\CC\PP^n}$, so that $\ov{\orb\langle \nu_{-1},\nu_0\rangle}=E$, $\ov{\orb \langle \nu_0,\nu_1\rangle}=H$ and $\wt{\gamma}$ sends $E$ biholomorphically to $H$.  

Let $\Delta:=\{t\in\CC|F(t\BB)\subset\BB\}\subset\CC$, which is an open disk containing $1$. Put 
\begin{equation*}
\Omega:=\{(z,t)\in T_N\times\Delta^*\mid z\in W_T(F(t\cdot))\}\subset T_{\wt{N}}
\end{equation*}
and define the smooth manifold 
\begin{equation*}
\wt{\mathcal{X}}_N:=\Int(\ov{\Omega}^{X(\wt{\Sigma},\wt{N})})\subset X(\wt{\Sigma},\wt{N})
\end{equation*}
together with a free and proper action of $\Gamma=\langle\wt{\gamma}\rangle$ and a $\Gamma$-invariant map $p:\wt{\mathcal{X}}_N\rightarrow\Delta$ so that $p^{-1}(t)=\widetilde{X(\pi,\sigma(t\cdot))}$ for $t\neq 0$ and $p^{-1}(0)$ is given by $\ZZ$ copies of $\widetilde{\CC\PP^n}$ glued along $H$ and $E$ via $\wt{\gamma}$. Thus we obtain a flat proper holomorphic map
\begin{equation*}
p_{N}:\mathcal{X}_N:=\wt{\mathcal{X}}_N/\Gamma\rightarrow \Delta
\end{equation*}
with $p_{N}^{-1}(t)=X(\pi,\sigma(t\cdot))$ for $t\in\Delta\sminus \{0\}$ and $p^{-1}_{N}(0)=\widetilde{\CC\PP^n}/_{E\sim_{\wt{\gamma}} H}$. This is precisely the Nakamura degeneration.

Define the toric divisor on $\wt{\mathcal{X}}_N$:
\begin{equation*}
\wt{\mathcal D}=\sum_{\nu\in\wt{\Sigma}^{(1)}}\ov{\orb\nu}\cap\wt{\mathcal X}_N.
\end{equation*}
It is clear that $\wt{\mathcal D}|_{\Delta^*}=\wt{D}^X_T\times\Delta^*$, so that $\wt{\mathcal{D}}\cap p^{-1}(t)=\wt{D}^X_T$ is the maximal toric divisor of $\widetilde {X(\pi,\sigma(t\cdot))}$ for each $t\in\Delta^*$. 
Since moreover $\wt{\mathcal D}$ is $\Gamma$-invariant, we obtain a divisor $\mathcal D:=\wt{\mathcal D}/\Gamma$ on $\mathcal X_N$ so that for each $t\in\Delta^*$, $\mathcal D\cap p_N^{-1}(t)=D_T$ and in particular  $K_{\mathcal X_t}\otimes L_{\det A}=\mathcal O(-\mathcal D)_t$. We have thus obtained:

\begin{theorem}\label{Nakadeg}
Let $X(\pi,\sigma)$ be an $n$-dimensional toric Kato manifold, let $\widetilde{\CC\PP^n}=X(\Sigma_c,N)$ be defined as above and let $S=\widetilde{\CC\PP^n}/_{E\sim_{\wt{\gamma}} H}$, where $\wt{\gamma}$ was defined above. 
There exists a  flat holomorphic proper family $p_N:(\mathcal{X}_N,\mathcal D)\rightarrow\Delta$ over some disk $\Delta\subset\CC$ of radius $R>1$, so that:
\begin{itemize}
\item $\mathcal X_N$ is smooth
\item $\mathcal D$ is an effective divisors on $\mathcal X_N$
\item for each $t\in\Delta^*$, $p_N^{-1}(t)\cong X(\pi,\sigma(t\cdot))$ and $\mathcal{D}\cap p_N^{-1}(t)$ is the toric divisor of $p_N^{-1}(t)$
\item $p_N^{-1}(0)=S$.
\end{itemize}
\end{theorem}

Note that by the discussion of Section~\ref{isom}, the smooth fibers of the family $\mathcal{X}_N$ are all isomorphic if $c\in\im(\id-A)$. 

\subsection{Isotrivial degenerations} Nakamura degenerations are defined for all Kato manifolds, however generally the smooth fibers are not all isomorphic. On the other hand, for most toric Kato manifolds we can define other degenerations, whose smooth fibers will all be isomorphic, but whose central fiber is less concrete and generally more singular. Moreover, the total space of these degenerations need not be smooth. These follow the same ideas as \cite[Section~2]{tsu}.

In order to do so, we first need the following:
\begin{lemma}\label{intv}
Let $A$ be a toric Kato matrix. Then $\im(\id-A)\cap \Int C_0\cap N\neq\emptyset$ if and only if $|P(A)|<n-1$.
\end{lemma}
\begin{proof}
If $|P(A)|=n-1$, so that $j\notin P(A)$, then it is clear that $\im(\id-A)\subset\ker e_j^*$, thus $\im(\id-A)\cap \Int C_0=\emptyset$.

Suppose now that $|P(A)|<n-1$, and let $j\neq k\in P(A)^c$. Suppose that $l\geq 1$ is such that $A_{P(A)^c}$ is positive, and take
\begin{equation*}
v:=(\id+A+\ldots+A^{l-1})(e_j+e_k).
\end{equation*} 
Then $w:=(A-\id)v\in N\cap \im(\id-A)$. Furthermore, we have $w=(A^l-\id)(e_j+e_k)$, and by the choice of $l$ it follows also that $w\in\Int C_0$. 
\end{proof}

\qquad

Let us thus suppose that $X=X(\pi,\sigma)$ is an $n$-dimensional toric Kato manifold with $\sigma=\un\la\sigma_A$ and $|P(A)|<n-1$. Let $F=\un\la F_A$ be its corresponding germ. We will again define a toric manifold $X(\wt{\Sigma},\wt{N})$ together with a toric map to $\CC$.

In order to do so, choose, via \ref{intv}, a primitive element $u\in \im(\id-A)\cap \Int C_0\cap N$, and let $v\in N$ with $Au=Av-v\in\Int AC_0$. Consider again $\wt{N}=N\oplus \ZZ e$, and define an action of  $\wt{A}'$ on $\wt{N}$ by putting $\wt{A}'|_N=A$ and $\wt{A}'e=e$. Put $\nu_k:=(\wt{A}')^k(e+v)$, $k\in\ZZ$, and define fans $\wt{\Sigma}_0,\wt{\Sigma}_u$ in $\wt{N}$ as before. Putting again $\wt{\un\la}=(\un\la, 1)\in T_{\wt{N}}$, we have the natural action of $\wt{\gamma}'$ given by $\wt{\un\la}\wt{A}'$, which restricted to $T_{\wt{N}}$ reads $\wt{\gamma}'(z)=\wt{F}'(z):=\wt{\un\la}. z^{\wt{A}'}$.  

Let again $p_\ZZ:\wt{N}\rightarrow \ZZ e$ be the natural projection, and let $\Sigma_1=\{0,\langle e\rangle \}$ be the fan of $\CC$ in $\ZZ e$. Then $p_\ZZ$ induces a map of toric varieties $p:X(\wt{\Sigma}_u,\wt{N})\rightarrow X(\Sigma_1,\ZZ e)$, with:
\begin{gather*}
p^{-1}(\CC^*)=X(\Sigma_A,N)\times\CC^*\\
p^{-1}(0)=\bigcup_{m\in\ZZ}\overline{\orb\nu_m}.
\end{gather*}
For $(x,t)\in X(\Sigma_A,N)\times\CC^*$ we have $\wt{\gamma}'(x,t)=(\gamma(x),t)$. 

The toric variety $\ov{\orb \nu_0}=X(\widehat{St}_{\nu_0}, N_{\nu_0})$ is isomorphic to $X(\Sigma_u,N)$, via the projection $q:\wt{N}\rightarrow N$, $q|_N=\id$, $qe=-v$ which maps the fan $\widehat{St}_{\nu_0}$ isomorphically to $\Sigma_u$ given by:
\begin{equation*}
\Sigma_u:=\Sigma_0\cup\{\tau'=\langle Au\rangle\oplus A\tau, \ \tau''=\langle -u\rangle\oplus\tau\mid \tau\in\Sigma_0, \tau\subset\del C_0\}.
\end{equation*}
Under this isomorphism, $\wt{\gamma}'$ sends $E':=\ov{\orb \langle -u\rangle}$ biholomorphically to $H':=\ov{\orb\langle Au\rangle}$.

Taking now $\Omega_{\wt{F}'}:=W_T(F_{\un\la,A})\times\CC^*$ and 
\begin{equation*}
\wt{\mathcal X}_u:=\Int(\ov{\Omega_{\wt{F}'}}^{X(\wt{\Sigma}_u,\wt{N})})\subset X(\wt{\Sigma}_u,\wt{N})
\end{equation*}
with the $\Gamma=\langle\wt{\gamma}'\rangle$-invariant projection $p:\wt{\mathcal X}_u\rightarrow \CC$, we obtain a proper holomorphic family
\begin{equation*}
p_{u}:\mathcal{X}_u:=\wt{\mathcal{X}}_u/\Gamma\rightarrow \CC
\end{equation*}
with $p_{u}^{-1}(t)=X$ for $t\in\CC^*$ and $p^{-1}_{u}(0)=X(\Sigma_u,N)/_{E'\sim_{\wt{\gamma}'} H'}$. Moreover, $p_u$ is a flat map by \cite[6.1.5]{ega}  since it is equidimensional and $\mathcal{X}_u$ is Cohen-Macaulay.

As before, define the divisors
\begin{gather*}
\wt{\mathcal D}=\sum_{\nu\in\wt{\Sigma}^{(1)}}\ov{\orb\nu}\cap\wt{\mathcal X}_N\in\mathrm{Div}\wt{\mathcal X},\ \  \mathcal D=\wt{\mathcal D}/\Gamma\in\mathrm{Div}{\mathcal X}
\end{gather*}
Then again one obtains that $\mathcal D\cap p_u^{-1}(t)$ is the toric divisor of $X$, for each $t\in\CC^*$.

Note that the central fiber of the family $\mathcal{X}_u$ depends on $u$, and that any primitive element $u\in\im(\id-A)\cap \Int C_0\cap N$ gives rise to such a family. However,  we cannot in general ensure that $\mathcal{X}_u$ is a smooth manifold, or equivalently that the normalization of the central fiber is smooth, as the following result shows:

\begin{proposition}
The toric variety $X(\Sigma_u,N)$ is smooth if and only if $u=c=\sum_{j=1}^ne_j$. In this case, the family $(\mathcal{X}_u,\mathcal D)$ is equivalent to the Nakamura family $(\mathcal X_N,\mathcal D_N)$, in the sense that there exists a biholomorphism $\phi_X:\mathcal{X}_N\rightarrow\mathcal{X}_c|_\Delta=p_c^{-1}(\Delta)$ with $p_c\circ\phi_X=p_N$, $\phi_X^*\mathcal D|_\Delta=\mathcal D_N$.
\end{proposition}
\begin{proof}
It is clear that if $u=c$, then $X(\Sigma_u,N)\cong \widetilde{\CC\PP^n}$, so it is smooth. Conversely, suppose that $X(\Sigma_u,N)$ is smooth. Since for any $j\in\{1,\ldots, n\}$, $-u,e_1,\ldots, \widehat{e_j}, \ldots, e_n$ is a basis of $N$, and since $u\in\Int C_0$, this implies that $u=c$.

Suppose now that $u=c$. Define $\phi_\ZZ:\wt{N}\rightarrow \wt{N}$ by $\phi_\ZZ|_N=\id$, $\phi_\ZZ e=e+v$, where $v$ is chosen to satisfy $Ac=Av-v$. We have then $\phi_\ZZ(\wt{\Sigma})=\wt{\Sigma}_c$ and $\phi_\ZZ \wt{A} x=\wt{A}'\phi_\ZZ x$ for any $x\in \wt{N}$. This implies that $\phi_\ZZ$ induces a $\Gamma$-equivariant biholomorphism $\phi:X(\wt{\Sigma},\wt{N})\rightarrow X(\wt{\Sigma}_c,\wt{N})$ commuting with the projections to $\CC$ and satisfying $\phi\circ \wt{F}=\wt{F}'\circ\phi$. In particular, it is clear that $\phi(\Omega)=\Omega_{\wt{F}'}\cap p^{-1}(\Delta)$ and that $\phi^*\wt{\mathcal D}=\wt{\mathcal D_N}$, so $\phi$ also induces a biholomorphism between $\wt{\mathcal{X}}_N$ and $\wt{\mathcal{X}}_c|_{\Delta}$, and so a biholomorhism $\phi_X:\mathcal{X}_N\rightarrow\mathcal{X}_c|_\Delta$ with $p_c\circ\phi_X=p_N$ and $\phi_X^*\mathcal D|_\Delta=\mathcal D_N$.
\end{proof}

The conclusion of the discussion is the following:
\begin{theorem}\label{Tsudeg}
Let $X$ be a toric Kato manifold of hyperbolic type of matrix $A$, and let $u\in\im(\id-A)\cap \Int C_0\cap N$ be a primitive element. Then there exists  a  flat holomorphic proper family $p_u:(\mathcal{X}_u,\mathcal D)\rightarrow\CC$ so that:
\begin{itemize}
\item $\mathcal D$ is an effective divisor on $\mathcal X_u$
\item for each $t\in\CC^*$, $p_u^{-1}(t)\cong X$ and $\mathcal{D}\cap p_u^{-1}(t)$ is the toric divisor of $X$
\item $p_u^{-1}(0)=X(\Sigma_u,N)/_{E'\sim_{\wt{\gamma}'}H'}$. 
\end{itemize}
Moreover, $\mathcal X_u$ is smooth if and only if $u=c$, in which case the family $(\mathcal X_u,\mathcal D)$ is equivalent to the Nakamura family $(\mathcal X_N,\mathcal D)$ of \ref{Nakadeg}.
\end{theorem}

\section{Analytic invariants}\label{analytic}

In this section, we compute different analytical invariants of toric Kato manifolds, such as the Hodge numbers $h^{0,q}$, $h^{1,q}$ and $h^{p,0}$, as well as the cohomology of the logarithmic sheaves $\Theta_X(-\log D_T)$ and $\Omega^1(\log D_T)$. To this aim, we wil make use of the degenerations of Secton~\ref{degen} and we will extend the method of Tsuchihashi \cite{tsu}, who computed these invariants  when $P(A)=\emptyset$ and $\underline\lambda=(1,\ldots, 1)$. Using a different method, Sankaran computed in \cite{san} all the Hodge numbers of toric Kato manifolds with $\underline\lambda=(1,\ldots, 1)$, under the additional assumption that $A$ is diagonalisable over $\CC$ and irreducible over $\QQ$, so in particular $P(A)=\emptyset$, obtaining 
$$h^{0,1}=h^{n,n-1}=1, \quad h^{p,p}=b_{2p}, \quad h^{p,q}=0 \text{ otherwise}.$$

We should note that the number $|P(A)|$ is an invariant of toric Kato manifolds - it can be retrieved, for instance, as the number of non-compact irreducible toric divisors of the universal cover, cf. Section~\ref{geomprop}. Moreover, even when $P(A)=\emptyset$ but $1\in\Spec(A)$ (see \ref{ex:nontrivialiso3d}), one cannot simply reduce to the case $\underline\lambda=(1,\ldots, 1)$ (see \ref{cor:condiso}). This shows that our context is much more general than both Sankaran's rank one actions or Tsuchihashi's. 
Sankaran's proof strongly depends on the above mentioned hypotheses, which are quite restrictive from the point of view of toric Kato manifolds. Therefore it would be interesting to know if one can adapt his strategy  to our more general situation.

\begin{proposition}\label{holoFcn}
Let $\wt{X}$ be the universal cover of a toric Kato manifold $X$. 
If $X$ is of hyperbolic type, then $H^0(\wt{X},\mathcal O_{\wt{X}})=\CC$. If $X$ is of parabolic type, then $H^0(\wt{X},\mathcal O_{\wt{X}})=\Oo_\CC(\CC)$.
\end{proposition}
\begin{proof}
Let $A$ be the matrix of $X$. Let $f\in H^0(\wt{X},\mathcal O_{\wt{X}})$. Since $\Omega$ is a $\TT$-invariant domain of $T_N$, we can express $f$ as a series on $\Omega$:
\begin{equation*}
f(z)=\sum_{m\in M}c_m z^m, \ \ z\in\Omega.
\end{equation*}
Since the right hand side extends naturally to a meromorphic function on $\wt{X}$, $f$ writes on $\wt{X}$ as 
\begin{equation}
f=\sum_{m\in M} c_m\wt{m}
\end{equation}
where for each $m\in M$,  $\wt{m}$ is the natural meromorphic function induced by $m$ on $X(\Sigma_A,N)$. 

Now since for any $\tau\in\Sigma^{(1)}_A$, $f$ has no poles on $\ov{\orb\tau}$, it follows that for any $m\in M$ with $c_m\neq 0$ we have
\begin{equation*}
\langle m, v\rangle \geq 0, \ \forall v\in\tau.
\end{equation*}
Since $C:=\ov{|\Sigma_A|}=\ov{H(A)}$ is a convex cone by \ref{supportSigma}, we readily infer:
\begin{equation*}
f=\sum_{m\in \check{C}\cap M}c_m\wt{m}.
\end{equation*}

Now from \ref{supportSigma} we have $\check C\cap M=\RR_{\geq 0}f_A^*\cap M=\{0\}$ if $|P(A)|<n-1$ and $\check C\cap M=\NN e_j^*$ if $P(A)^c=\{j\}$. From this it follows that in the parabolic case, $f$ is an entire function in the $z_j$-variable. 
\end{proof}

\begin{remark}\label{sectionFlat}
For a parabolic toric Kato manifold $X$ with $P(A)^c=\{j\}$, define 
\begin{gather}
\nonumber \rho:\Gamma\cong\ZZ\langle F_{\un\la_0,A}\rangle\rightarrow \CC^*,\ \  \rho(F_{\un\la_0,A})=\la_j\\
\mathcal L=\widetilde{X(\pi,\sigma)}\times_\rho \CC
\end{gather}
where $\la_j\in\DD$ is the $j$-th component of $\un\la_0$. Then it follows that the holomorphic function on the universal cover defined by $f(z)=z_j$ on $T_N$ gives rise to a holomorphic section of $\mathcal L$. 
\end{remark}

\begin{corollary}
If $X$ is a toric Kato manifold which is not Hopf and $L\in\Pic_0(X)$ is non-trivial, then $H^0(X, L)\neq 0$ if and only if $X$ is of parabolic type and $L$ is a positive power of the line bundle defined in \ref{sectionFlat}. In this case, $H^0(X, L)\cong\CC$.
\end{corollary}

Next, let $\Theta_X(-\log D_T)$ denote the locally free subsheaf of the tangent sheaf given by germs of holomorphic vector fields tangent to $D_T$, let $\Omega^1_X(\log D_T)$ be the dual sheaf and let $\Omega^p(\log D_T):=\bigwedge^p(\Omega^1(\log D_T))$ for $p\geq 2$.

\begin{proposition}\label{holF}
If $X$ is any toric Kato manifold, then $H^0(X,\Omega_X^p)=0$ for any $p\geq 1$.
\end{proposition}
\begin{proof}
For Hopf manifolds the result is well known, so we only treat the other cases. We start by noting that the map 
\begin{equation*}
M\otimes\Oo_{\wt{X}}\rightarrow\Omega^1_{\wt{X}}(\log\wt{D}_T), \ \ m\otimes f\rightarrow f\cdot\frac{d\wt{m}}{\wt{m}}
\end{equation*}
where $\wt{m}$ is the natural meromorphic function on $\wt{X}$ associated to $m$, is an isomorphism of $\Oo_{\wt{X}}$-modules \cite[Proposition at page~87]{fulton}. In particular, for any $p\geq 1$, $\Omega_{\wt{X}}^p(\log\wt{D}_T)\cong \bigwedge^p M\otimes\Oo_{\wt{X}}$.

Let $f_1,\ldots, f_n$ be a $\ZZ$-basis of $N$ so that $\nu_j:=\langle f_j\rangle\in\Sigma_A^{(1)}$ for each $1\leq j \leq n$ and let $m_1,\ldots, m_n$ be the dual basis of $M$. For each $J=(1\leq j_1<\ldots<j_p\leq n)$, $1\leq p\leq n$, denote by 
\begin{equation*}
\omega_J:=\frac{d\wt{m}_{j_1}}{\wt{m}_{j_1}}\wedge\ldots\wedge \frac{d\wt{m}_{j_p}}{\wt{m}_{j_p}}\in H^0(\wt{X},\Omega^p_{\wt{X}}(\log D_T)).
\end{equation*}

 Using \ref{holoFcn}, we find that
\begin{equation*}
H^0(\wt{X},\Omega^p_{\wt{X}}(\log \wt{D}_T))=\left\{\begin{array}{ll}  \mathrm{span}_{\CC}\langle \omega_J\mid |J|=p\rangle, & |P(A)|<n-1\\
\mathrm{span}_{\Oo_\CC(\CC)}\langle \omega_J\mid |J|=p\rangle, 
&P(A)^c=\{j\}\end{array}\right.
\end{equation*}
where $\Oo_\CC(\CC)$ denotes the ring of entire functions in $z_j=\wt{e}_j^*$.

Since $\Omega^p_{\wt{X}}$ is a subsheaf of $\Omega^p_{\wt{X}}(\log\wt{D}_T)$ and since each $\omega_J$  has poles along $\sum_{j\in J}\ov{\orb\nu_j}$, it is clear that in the case $|P(A)|<n-1$ we have $H^0(X,\Omega^p_X)=H^0(\wt{X},\Omega^p_{\wt{X}})=0$ for any $p\geq 1$. In the case $P(A)^c=\{j\}$, the same argument shows that $H^0(X,\Omega^p_X)=H^0(\wt{X},\Omega^p_{\wt{X}})=0$ for $p\geq 2$. For $p=1$, we find that 
\begin{equation*}
H^0(\wt{X},\Omega^1_{\wt{X}})=\{ \omega = f(z_j)dz_j\mid f\in\Oo_\CC(\CC)\}.
\end{equation*}
Since $\gamma=F_{\un\la,A}$ acts on such an element $\omega=\sum_{k=0}^\infty c_kz_j^k dz_j$ by $\gamma^*\omega=\sum_{k=0}^\infty\la_j^{k+1}c_kz_j^kdz_j$ and as $|\la_j|<1$, we find
\begin{equation*}
H^0(X,\Omega^1_X)=H^0(\wt{X},\Omega^1_{\wt{X}})^\Gamma=0
\end{equation*}
which concludes the proof.
\end{proof}

\begin{theorem}\label{cohOX}
Let $X$ be a toric Kato manifold of hyperbolic type. Then we have:
\begin{equation*}
H^0(X,\Oo_X)=H^1(X,\Oo_X)=\CC, \ \ H^p(X,\Oo_X)=0, \ p\geq 2.
\end{equation*}
\end{theorem}
\begin{proof}
Let $A$ be a toric Kato matrix for $X$. Since $|P(A)|<n-1$, there exists $u\in\im(\id-A)\cap \Int C_0\cap N$  a primitive element. Consider the family $p_u:\mathcal X_u\rightarrow \CC$ of \ref{Tsudeg}, of central fiber $S$ and generic fiber $X$. Then by upper-semicontinuity \cite[Theorem~4.12]{bs} we have:
\begin{equation}\label{semicO}
\dim H^p(X,\Oo_X)\leq \dim H^p(S,\Oo_S), \ \ p\in\NN.
\end{equation}

Let us first compute $H^\bullet(S,\Oo_S)$. Let $n:\hat{S}=X(\Sigma_u,N)\rightarrow S$ be the normalization map, and let $Y=n(E')=n(H')$ be the double locus of $S$, with inclusion map $j:Y\rightarrow S$. Then we have an exact sequence of sheaves on $S$:
\begin{equation}\label{normalseq}
\xymatrix{
0\ar[r] & \Oo_S\ar[r] &n_*\Oo_{\hat{S}} \ar[r] & j_*\Oo_Y\ar[r] &0
}
\end{equation}
and the corresponding long exact sequence in cohomology. Since $Y$ and $\hat{S}$ are compact toric varieties, we have $H^p(Y,\Oo_Y)=H^p(\hat{S},\Oo_{\hat{S}})=0$ for any $p\geq 1$ cf. \cite[Corollary~7.4]{da}. This implies that $H^p(S,\Oo_S)=0$ for $p\geq 2$ and gives rise to the exact sequence:
\begin{equation*}
\xymatrix{
0\ar[r] &H^0(S,\Oo_S)\ar[r] &\CC\ar[r]^0 &\CC\ar[r] &H^1(S,\Oo_S)\ar[r] &0.}
\end{equation*}
We infer that $H^0(S,\Oo_S)=H^1(S,\Oo_S)=\CC$.

The cohomology of $\Oo_S$ together with \eqref{semicO} immediately gives the desired result for $H^p(X,\Oo_X)$ for $p\neq 1$. On the other hand, since $h^{1,0}(X)=0$ by \ref{holF}, using again \eqref{semicO} we find:
\begin{equation*}
1=b_1(X)\leq h^{1,0}(X)+h^{0,1}(X)=h^{0,1}(X)\leq h^{0,1}(S)=1
\end{equation*}
hence we also have $H^1(X,\Oo_X)=\CC$. This concludes the proof.
\end{proof}

\begin{corollary}
Let $X$ be a toric Kato manifold of hyperbolic type. Then we have an isomorphism $\Pic_0(X)\cong H^1(X,\CC^*)\cong \CC^*$ given by
\begin{equation*}
\la\in\CC^*\mapsto L_\la:=\wt{X}\times_{\rho_\la}\CC, \ \ \rho_\la\in\Hom(\Gamma,\CC^*), \ \rho_\la(n)=\la^n.
\end{equation*}
Moreover, we have a short exact sequence:
\begin{equation*}
\xymatrix{
0\ar[r] &\CC^*\ar[r] &\Pic(X)\ar[r]^-{c_1} &H^2(X,\ZZ)\ar[r]&0.
}
\end{equation*}
\end{corollary}
\begin{proof}
Since by \ref{cohOX} we have $H^1(X,\Oo_X)=\CC=H^1(X,\CC)$, we find that $\Pic_0(X):=H^1(X,\Oo)/H^1(X,\ZZ)\cong H^1(X,\CC^*)=\CC^*$. Since moreover $H^2(X,\Oo_X)=0$, the exponential sequence gives rise to the exact sequence:
\begin{equation*}
\xymatrix{ 0\ar[r] &\CC^*\ar[r] & \Pic(X)\ar[r]^-{c_1} &H^2(X,\ZZ)\ar[r] &H^2(X,\Oo_X)=0}
\end{equation*}
from which we conclude.
\end{proof}

\begin{theorem}\label{cohDT}
Let $X$ be a toric Kato manifold with toric divisor $D_T$ and matrix $A$. Then we have:
\begin{gather*}
H^0(X,\Theta_X(-\log D_T))\cong \ker(A-\id)\\
H^0(X,\Omega_X^1(\log D_T))\cong \ker(\transp{A}-\id). 
\end{gather*}
If moreover $X$ is of hyperbolic type, then we also have:
\begin{gather*}
H^1(X,\Theta_X(-\log D_T))\cong \coker(A-\id), \ \ H^p(X,\Theta_X(-\log D_T))=0, \ p\geq 2\\
H^1(X,\Omega_X^1(\log D_T))\cong \coker(\transp{A}-\id), \ \ H^p(X,\Omega_X^1(\log D_T))=0, \ p\geq 2.
\end{gather*}
\end{theorem}
\begin{proof}
We only show the statements for the sheaf $\Theta_X(-\log D_T)$, since for $\Omega_X^1(\log D_T)$ the proof is dual. We have 
\begin{equation*}
H^0(X,\Theta_X(-\log D_T))=H^0(\wt{X},\Theta_{\wt{X}}(-\log \wt{D}_T))^{\Gamma}.
\end{equation*}
Since $\wt{X}$ is an open set of the toric manifold $X(\Sigma_A,N)$, we find $\Theta_{\wt{X}}(-\log \wt{D}_T)\cong \Oo_{\wt{X}}\otimes N$ and so, by \ref{holoFcn}, the conclusion follows for the case $|P(A)|<n-1$.

If $P(A)^c=\{j\}$, then by \ref{holoFcn}, any element $s\in H^0(\wt{X},\Oo_{\wt{X}}\otimes N)$ can be expressed as:
\begin{equation*}
s=\sum_{m\in\NN}a_mz_j^m, \ \ a_m\in N\otimes\CC.
\end{equation*}
The $\Gamma$-invariance of $s$ reads:
\begin{equation*}
Aa_m=\la_j^{-m}a_m, \ \ \forall m\in\NN.
\end{equation*}
Using that $|\la_j|<1$ and that for any $\mu\in\Spec(A)$, $|\mu|=1$, we find: 
\begin{equation*}
H^0(X,\Theta_X(-\log D_T))\cong\ker(A-\id).
\end{equation*}

Suppose next that $|P(A)|<n-1$, so that there exists a primitive element $u\in\im(\id-A)\cap \Int C_0\cap N$, and consider the family $p_u:(\mathcal X_u, \mathcal D)\rightarrow \CC$ of \ref{Tsudeg}, of central fiber $S$ and generic fiber $X$. Consider the sheaf $\mathcal F:=\Theta_{\mathcal X}(-\log\mathcal D)$, and denote by $\mathcal F_t$ its restriction to any fiber $p_u^{-1}(t)$, $t\in\CC$.   
Then for any $t\in\CC^*$, we have $\mathcal F_t\cong\Theta_X(-\log D_T)\oplus \Oo_X$, so by upper-semicontinuity \cite[Theorem~4.12]{bs} we find:
\begin{equation}\label{semicDT}
\dim H^p(X,\Theta_X(-\log D_T))+\dim H^p(X,\Oo_X)\leq \dim H^p(S,\mathcal F_0), \ \ p\in\NN.
\end{equation}

We must thus compute $H^\bullet(S,\mathcal F_0)$. Consider again the exact sequence \eqref{normalseq} and tensorize it with the flat sheaf $\mathcal F_0$, obtaining:
\begin{equation*}
\xymatrix{
0\ar[r]&\mathcal F_0\ar[r] &n_*\Oo_{\hat{S}}\otimes\wt{N} \ar[r] &j_*\Oo_Y\otimes \wt{N}\ar[r]&0.}
\end{equation*}
As before, passing to the long exact sequence in cohomology we find that  $H^p(S,\mathcal F_0)=0$ for $p\geq 2$ and so $H^p(X,\Theta_X(-\log D_T))=0$ for $p\geq 2$. Moreover, we have the short exact sequence:
\begin{equation*}
\xymatrix{
0\ar[r] &H^0(S,\mathcal F_0)\ar[r] &\wt{N}\otimes\CC\ar[r]^{\wt{A}'-\id} &\wt{N}\otimes\CC\ar[r] &H^1(S,\mathcal F_0)\ar[r] &0}
\end{equation*}
which implies $H^1(S,\mathcal F_0)=\coker(\wt{A}'-\id)=\coker (A-\id)\oplus\CC$. Thus, from \eqref{semicDT} and \ref{cohOX} we infer the inequality:
\begin{equation}\label{semicDT2}
\dim H^1(X,\Theta_X(-\log D_T))\leq\dim\coker(A-\id).
\end{equation}

At the same time, $H^1(X,\Theta_X(-\log D_T))$ can be computed, as $\Gamma$-equivariant cohomology, by a spectral sequence \cite[Subsection~5.2]{gr}. Namely, if $q:\wt{X}\rightarrow X$ denotes the covering map, then $\Theta_{\wt{X}}(-\log\wt{D}_T)$ is a $\Gamma$-sheaf on $\wt{X}$ and $\Theta_X(-\log D_T)=q_*^\Gamma\Theta_{\wt{X}}(-\log\wt{D}_T)$. Therefore, we have a spectral sequence 
\begin{equation*}
E^{p,q}_2=H^p(\Gamma,H^q(\wt{X},\Theta_{\wt{X}}(-\log \wt{D}_T)))\Rightarrow H^{p+q}(X,\Theta_X(-\log D_T)).
\end{equation*}
Moreover, since $\Gamma\cong\ZZ$, we have $E_2^{p,\bullet}=0$ for any $p\geq 2$. Hence this spectral sequence degenerates at the second page, implying that $H^1(X,\Theta_X(-\log D_T))=E_2^{1,0}\oplus E_2^{0,1}$. In particular, we have
\begin{equation}\label{revsemic}
\dim H^1(X,\Theta_X(-\log D_T))\geq \dim H^1(\Gamma,H^0(\wt{X},\Theta_{\wt{X}}(-\log\wt{D}_T))).
\end{equation}
But using \ref{holoFcn}, the fact that $\Ss^1$ is the  Eilenberg–MacLane space of $\Gamma$ and Poincaré duality for local systems, we find:
\begin{align*}
H^1(\Gamma,H^0(\wt{X},\Theta_{\wt{X}}(-\log\wt{D}_T)))&=H^1(\Gamma,N\otimes\CC)\cong H^0(\Gamma,N^*\otimes\CC)^*\\
&=(\ker(\transp{A}-\id))^*=\coker (A-\id) .
\end{align*}
This, together with \eqref{revsemic} and \eqref{semicDT2} gives us the desired conclusion.
\end{proof}

\begin{theorem}\label{cohOX2}
Let $X$ be a toric Kato manifold of hyperbolic type, with divisor $D$ induced by the exceptional divisor of $\pi$. Then we have:
\begin{gather*}
h^{1,p}(X)=0, \ p\neq 1\\
 \ h^{1,1}(X)=\sharp D>0.
 \end{gather*}
\end{theorem}
\begin{proof}
Let $A$ be a Kato matrix for $X$, let $s$ be its corresponding permutation and let $C$ denote the set of cycles of $s$. Let us write:
\begin{equation*}
D_T=\sum_{j=1}^k D_j=D+\sum_{J\in C}Y_J,
\end{equation*} 
where each $D_j$ is an immersed irreducible hypersurface. Let us put $t:=\sharp C$, so that $k=\sharp D+t$. Let $\hat{D}_j$ be the normalization of $D_j$ for each $j=1,\ldots, k$ and let $\al:\hat{D}_T=\bigsqcup_{j=1}^k\hat{D}_j\rightarrow X$ be the natural map given by inclusion and normalization. Then we have an exact sequence:
\begin{equation}\label{Residue}
\xymatrix{
0\ar[r] &\Omega_X^1\ar[r]&\Omega^1_X(\log D_T)\ar[r]^-P &\al_*\Oo_{\hat{D}_T}\ar[r] &0
}
\end{equation}
where $P$ is the Poincaré residue map. 

Each component of $\hat{D}_T$ is either a compact toric variety, or a toric Kato manifold of hyperbolic type, and we have precisely $t$ components of the last type. Using \ref{cohOX} and \ref{cohDT}, we infer that $H^p(X,\Omega^1_X)=0$ for $p\geq 3$ and we find the following exact sequence:
\begin{gather}\label{ExactSeq}
\nonumber \xymatrix{0\ar[r]&H^0(X,\Omega^1_X(\log D_T))\ar[r]&\CC^k\ar[r]&H^1(X,\Omega^1_X)\ar[r]&&}\\
\xymatrix{
&&\ar[r]&H^1(X,\Omega^1_X(\log D_T))\ar[r]^-{P_1}&\CC^{t}\ar[r]&H^2(X,\Omega^1_X)\ar[r]&0.
}
\end{gather}

Now the point is to show that the map $P_1$ induced by $P$ is surjective. This will readily imply that $h^{1,2}(X)=0$. Furthermore, as the Euler characteristic of the exact sequence is $0$, and as $\dim H^0(X,\Omega^1_X(\log D_T))=\dim H^1(X,\Omega^1_X(\log D_T))$ by \ref{cohDT}, we also infer that $h^{1,1}=k-t=\sharp D$. 

Let us thus show that $P_1$ is surjective. We recall that 
\begin{equation*}
q^*\sum_{J\in C}Y_J=\sum_{j\in P(A)}X_j
\end{equation*} 
where $q$ is the covering map, and for each $j\in P(A)$, $X_j=\ov{\orb\langle e_j\rangle}\cap \wt{X}$. The map $P_1$ is induced by the $\Gamma$-equivariant map of $\Gamma$-sheaves on $\wt{X}$:
\begin{gather*}
\wt{P}_1:\Omega^1_{\wt{X}}(\log \wt{D}_T)\cong M\otimes\Oo_{\wt{X}}\rightarrow \oplus_{j\in P(A)}\Oo_{X_j}\\
m\otimes f\mapsto (\langle m,e_j\rangle\cdot f|_{X_j})_{j\in P(A)}
\end{gather*}
so that 
\begin{equation*}
P_1=\wt{P}_1^*:H^1_\Gamma(\wt{X},\Omega^1_{\wt{X}}(\log \wt{D}_T))=H^1(X,\Omega^1_X(\log D_T))\rightarrow \oplus_{J\in C}H^1(\hat Y_J,\Oo_{\hat Y_J}). 
\end{equation*}

On the other hand, using again the spectral sequence of \cite{gr} which degenerates at the second page, together with \ref{cohOX} and \ref{cohDT}, we find
\begin{gather*}
\coker(\transp{A}-\id)=H^1(X,\Omega^1_X(\log D_T))=H^1(\Gamma, H^0(\wt{X},\Omega^1_{\wt{X}}(\log \wt{D}_T)))=H^1(\Gamma,M\otimes\CC),\\
\CC^t=\oplus_{J\in C}H^1(\hat Y_J,\Oo_{\hat Y_J})=H^1(\Gamma, \oplus_{j\in P(A)}H^0(X_j,\Oo_{X_j}))=H^1(\Gamma, \CC^{|P(A)|})
\end{gather*}
where $\Gamma$ acts on $\CC^{|P(A)|}$ via the permutation $s$. 
Hence $P_1$ is the map induced in group cohomology by $Q:M\otimes\CC\rightarrow \CC^{|P(A)|}$, $m\mapsto (\langle m, e_j\rangle)_{j\in P(A)}$. Clearly $Q$ is surjective, and let $K=\ker Q$, so that we have an exact sequence
\begin{equation*}
\xymatrix{ 0\ar[r] &K\ar[r] &M\otimes\CC\ar[r]^-Q &\CC^{|P(A)|}\ar[r] &0.
}
\end{equation*}
This further induces an exact sequence
\begin{equation*}
\xymatrix{
H^1(\Gamma, M\otimes \CC)\ar[r]^-{P_1} & H^1(\Gamma, \CC^{|P(A)|})\ar[r] & H^2(\Gamma, K).}
\end{equation*}
But $H^2(\Gamma, K)=0$ since $\Gamma=\ZZ$, so indeed $P_1$ is a surjective map. This concludes the proof of the theorem.
\end{proof}

\begin{remark}\label{invParabolic}
We note here that for toric Kato manifolds $X(\pi,\sigma_A(t\cdot))$ of parabolic type, the conclusions of \ref{cohOX}, \ref{cohDT} and \ref{cohOX2} also hold if $|t|$ is small enough. The arguments are exactly the same, except that one has to replace the families $p_u$ by the Nakamura family $p_N$ in \ref{cohOX} and \ref{cohDT}. Since in the Nakamura family the smooth fibers need not be isomorphic, the semi-continuity argument then works only for parabolic Kato manifolds close to the singular fiber. 
\end{remark}

\section{Isomorphism classes of toric Kato manifolds}\label{isom}

The aim of this section is to describe isomorphism classes of toric Kato manifolds, according to the combinatorial data provided by the toric modification $\pi:\hat{\CC}^n \to \CC^n$ (specifically the fan $\hat{\Sigma}$ of $\hat{\CC}^n$), and the Kato germ $F=F_{\un\lambda, A}$.

\subsection{Collapsing models}

To this aim, we need first to introduce some more constructions related to toric Kato manifolds, inspired by the collapsing models of \cite{dl84}, and generalizing constructions introduced above.

Denote by $\Sigma$ the fan of $\nC^n$, by $\hat\Sigma$ the fan of $\hat\CC^n$ and set $\annfan{\Sigma} = \hat{\Sigma} \sminus A\Sigma$, where $A\Sigma=\{A\tau\ |\ \tau \in \Sigma\}$. For any $\ell \in \nZ \cup \{-\infty\}$ and $m \in \nZ \cup \{+\infty\}$ with $\ell \leq m$ we set
$$
\colfan{A}{\ell}{m} := \bigcup_{\ell \leq k < m} A^k \annfan{\Sigma} \cup A^m \Sigma,
$$
where we set $A^{\pm \infty}\tau = \emptyset$ for any cone $\tau$.

\begin{lemma}\label{lem:colfan}
The set $\colfan{A}{\ell}{m}$ defines a regular fan, which is finite if and only if $\ell, m \in \nZ$.
When $\ell=-\infty$, the support is given by $\abs{\colfan{A}{-\infty}{m}}=H(A) \cup \tau_{P(A)}$ for any $m \in \nZ$. 
In general, we have $\abs{\colfan{A}{\ell}{m}} \supseteq \abs{\colfan{A}{\ell}{m'}}$
for all $\ell \in \nZ \cup\{-\infty\}$, $m \leq m' \in \nZ \cup \{+\infty\}$, $\ell \leq m$, with strict inclusion if and only if $m \in \nZ$ and $m'=+\infty$.
\end{lemma}
\begin{proof}
The first part of the statement is trivial, the rest follows directly from \eqref{eqn:unionfans}.
\end{proof}

To simplify notations, we will omit $\ell$ and $m$ when they are infinite.
Notice that $\colfan{A}{0}{d}=\hat{\Sigma}^d$ is the fan associated to $\pi^d$ as defined in the proof of \ref{powerKD};
$\colfan{A}{}{}$ corresponds to the infinite fan $\Sigma_A$ defined by \eqref{eqn:SigmaA}, while $\colfan{A}{0}{}$ corresponds to the infinite fan $\Sigma_A^+$ defined in the proof of \ref{secondDescr}.

To the fan $\colfan{A}{\ell}{m}$ is associated the toric variety (possibly of non-finite type) $X(\colfan{A}{\ell}{m},N)$.
When $m \in \nZ$, we denote by $0_m$ the point associated to the cone $A^m \tau_0$, where $\tau_0 \in \Sigma$ is the cone associated to the origin of $\nC^n$.
Notice that when $\ell=0$ and $m \in \nZ_{>0}$, the point $0_m$ is exactly the point $\sigma^m(0)$, where $\sigma^m$ is the second entry of the Kato data given by \ref{powerKD}.

All these toric varieties are related one to the other via some natural morphisms.

\begin{definition}
Let $\ell \in \nZ \cup \{-\infty\}$, and $m \leq m' \in \nZ \cup \{+\infty\}$, with $\ell \leq m$.
The toric morphism $\colmap{\ell}{m}{m'}:X(\colfan{A}{\ell}{m'},N) \to X(\colfan{A}{\ell}{m},N)$, induced by the identity of $N$, is called the \emph{collapsing map} from $X(\colfan{A}{\ell}{m'},N)$ to $X(\colfan{A}{\ell}{m},N)$.
\end{definition}

The $F$-invariant open set $W_T(F)$ sits inside the torus $T_N$, which is an open dense subset of $X(\colfan{A}{\ell}{m},N)$ for any $\ell \leq m$.
Similarly to Section~\ref{construct}, we set
\begin{equation*}
\colvar{F}{\ell}{m}:=\Int(\ov{W_T(F)}^{X(\colfan{A}{\ell}{m},N)})\subset X(\colfan{A}{\ell}{m},N).
\end{equation*}
The variety $\colvar{F}{\ell}{m}$ is called the ($\ell\,{}^\centerdot m$-)\emph{collapsed model} associated to the Kato data $(\pi,\sigma)$.
When $\ell=-\infty$ and $m \in \nZ$, this corresponds to the collapsed model defined in \cite[Chapter I.3]{dl84}. We note that since for any $\ell \in \nZ$, $\colfan{A}{\ell}{}\subset \Sigma_A$, $\colvar{F}{\ell}{}$ naturally sits inside $\colvar{F}{}{}$ as an open subset.

Notice that $\colmap{\ell}{m}{m'}$ leaves $W_T(F)$ invariant, and hence defines a regular map $\colmap{\ell}{m}{m'}:\colvar{F}{\ell}{m'} \to \colvar{F}{\ell}{m}$, which is a proper modification of $0_m$. 

Notice also that for any $m \in \nZ_{\geq 1}$, the collapsing maps $\colmap{0}{0}{m}:\colvar{F}{0}{m} \to \colvar{F}{0}{0} =W^s(F)$ correspond to the toric modifications $\pi^m$ given by \ref{powerKD}.

\begin{remark}\label{rmk:collapsinginversesystem}
The natural diagrams involving collapsing maps commute, and in particular, for any $\ell \leq m \leq m' \leq m''$, one has $\colmap{\ell}{m}{m''} = \colmap{\ell}{m}{m'} \circ \colmap{\ell}{m'}{m''}$.
This allows to construct projective limits of the families $\Big(X\big(\colfan{A}{\ell}{m},N\big)\Big)_{m \in \nZ_{\geq \ell}}$ and $\Big(\colvar{F}{\ell}{m}\Big)_{m \in \nZ_{\geq \ell}}$.
The spaces $X\big(\colfan{A}{\ell}{},N\big)$ and $\colvar{F}{\ell}{}$ sit naturally inside such projective limits as open dense subsets.
Their complement can be interpreted in terms of valuations (see, e.g., \cite{teissier}), and the whole projective limits in terms of hybrid spaces (see, e.g., \cite{berkovich, boucksom-jonsson, favre}). 
\end{remark}

\subsection{Isomorphism classes}

We will say that two toric Kato manifolds $X$ and $X'$ are \textit{equivariantly isomorphic} if there exists a biholomorphism $\Phi:X\rightarrow X'$ and a group isomorphism $\nu\in\Aut_{gr.}(T_N)$ so that the lift $\tilde\Phi:\tilde X\rightarrow \tilde X'$ satisfies $\tilde\Phi(\un\la x)=\nu(\un\la)\tilde\Phi(x)$ for any $x\in X$ and $\un\la\in T_N$ for which this is defined. 

Also, recall that a holomorphic germ $\Phi$ satisfying $\Phi \circ F = G \circ \Phi$ is called a \emph{semi-conjugacy} between $F$ and $G$. When $\Phi$ is moreover an invertible germ,  
we call $\Phi$ a \emph{conjugacy}.

We start by describing the compact hypersurfaces in the universal covering of a toric Kato variety.

\begin{lemma}\label{lem:compacthypersurfunivcover}
Let $X$ be a toric Kato manifold, and $\wt{X}$ its universal covering. Then the compact hypersurfaces of $\wt{X}$ are exactly the toric hypersurfaces whose sum gives $D_+$: they correspond to the rays of the fan $\Sigma_A$ associated to $X$ not contained in $\tau_{P(A)}$.
\end{lemma}
\begin{proof}
Let $H$ be any compact hypersurface of $\wt{X}$.
Being compact, it must be contained in the open set $\left(\colmap{\ell}{m}{}\right)^{-1}\left(\colvar{F}{\ell}{m} \sminus \{0_m\}\right) \subset \colvar{F}{}{}$, for some $\ell \leq m \in \nZ$.
Notice that $\colmap{\ell}{\ell}{m}: \colvar{F}{\ell}{m} \to \colvar{F}{\ell}{\ell} \subseteq X(\colfan{A}{\ell}{\ell},N) \cong \nC^n$ is a toric modification above the origin of $\nC^n$, and we infer that $H$ must correspond to one of the exceptional primes of the toric modification.
This means that $H$ corresponds to one of the rays of $\colfan{A}{\ell}{m}$, and so to one of the rays of $\Sigma_A$. Finally, it is clear that a ray $\nu$ of $\Sigma_A$ gives rise to a compact hypersurface in $\tilde X$ if and only if $\nu$ is not contained in $\tau_{P(A)}$. 
\end{proof}

\begin{theorem}\label{thm:generaliso}
Let $X$ and $X'$ be two isomorphic toric Kato manifolds, with associated Kato germs $F=F_{\un\la,A}$ and $G=F_{\un\mu,B}$.
Set $H_\infty=\bigcup_{m>0}F^{-m}(0)$ and $H'_\infty=\bigcup_{m>0}G^{-m}(0)$.
Then there exists a holomorphic germ $\Phi:(\nC^n,0) \to (\nC^n,0)$ which is an isomorphism on $\nC^n \setminus H_\infty$, satisfying $\Phi^{-1}(\nC^n \setminus H_\infty') = \nC^n \setminus H_\infty$ and $\Phi \circ F = G \circ \Phi$.
\end{theorem}
\begin{proof}
Let $\Psi:X \to X'$ be an isomorphism, and let $\wt{\Psi}:\wt{X} \to \wt{X}'$ be a lift of $\Psi$ to the universal coverings.
	
To distinguish collapsing maps for $X$ and $X'$, we denote collapsing maps for $X$ with the letter $p$, and collapsing maps for $X'$ with the letter $q$.
To simplify notations, we also set $p_m:= \colmap{m}{m}{}$ and $q_m:= \colmap{m}{m}{}[q]$.
Finally, we denote by  $D_+$ (resp., $D'_+$) the union of the compact hypersurfaces of $\wt{X}$ (resp., $\wt{X}'$), and set $\wt{D}_+^{'+}=q_m^{-1}(0_m) \subset D'_+$.

Let us take $m \ll 0$ small enough so that the images through $\wt{\Psi}$ of the compact hypersurfaces in $p_0^{-1}(0)$ belong to $(q_{m})^{-1}(0_m)$. This can be done since all such (infinitely many) compact hypersurfaces project to finitely many hypersurfaces of $X$.

For any $\ell\in\ZZ$, let $\wt{B}^{-\ell}:\colvar{G}{\ell}{\ell}\rightarrow \CC^n$ be the isomorphism induced by $B^{-\ell}$, and let $\tilde q_\ell:=\wt{B}^{-\ell}\circ q_\ell:\colvar{G}{\ell}{}\rightarrow\CC^n$. Note that for $\ell<\ell'$, $\tilde q_{\ell}|_{\colvar{G}{\ell'}{}}=\tilde q_{\ell'}$. Define 
\begin{equation*}
\Phi:=\tilde q_m \circ \wt{\Psi} \circ p_0^{-1}:\CC^n\setminus\{0\}\rightarrow \CC^n.
\end{equation*}
By Hartogs' theorem, the map $\Phi$ extends to a holomorphic map on $\CC^n$. Moreover, it is a biholomorphism outside 
$p_0(\wt{\Psi}^{-1}(\wt{D}_+^{'+}))\subseteq p_0(D_+)$. Thanks to the description of compact hypersurfaces given by \ref{lem:compacthypersurfunivcover}, we have $\wt{\Psi}(D_+)=D'_+$, from which we infer that $\Phi(\CC^n\setminus H_\infty)=\CC^n\setminus H_\infty'$ and that $\Phi$ is an isomorphism outside $H_\infty$.

We are left with verifying that $\Phi\circ F=G\circ \Phi$. Let $\gamma_F$ be the deck group generator of $\tilde X$, which is induced by $F$, and let $\gamma_G$ be the deck group generator of $\tilde X'$. Being $\wt{\Psi}$ the lift of $\Psi$, we have $\wt{\Psi} \circ \gamma_F=\gamma_G \circ \wt{\Psi}$.

Note that $\gamma_F$ induce isomorphisms $\colvar{F}{\ell}{} \to \colvar{F}{\ell+1}{}$ for any $\ell \in \nZ$. Moreover, it is straightforward to check that $F \circ p_0 =p_0|_{\colvar{F}{1}{}}  \circ \gamma_F$, and similarly $G \circ \tilde q_m = \tilde q_m|_{\colvar{G}{m+1}{}}  \circ \gamma_G$.
Then we have
\begin{align*}
\Phi \circ F 
&= \tilde q_m \circ \wt{\Psi} \circ p_0^{-1}\circ F 
=\tilde q_m \circ \wt{\Psi} \circ \gamma_F \circ p_0^{-1}\\
&=\tilde q_m \circ \gamma_G \circ \wt{\Psi} \circ p_0^{-1}
= G \circ \tilde q_m \circ \wt{\Psi} \circ p_0^{-1}
= G \circ \Phi
\end{align*}
which concludes the proof.
\end{proof}

In general, the semi-conjugacy given by \ref{thm:generaliso} is not necessarily equivariant (see \ref{ex:notequivariantconj}).
The next proposition shows that this is the case, as soon as the union of the coordinate hyperplanes is $\Phi$-invariant.
Notice that this condition is automatically satisfied when $P(A)=\emptyset$.

\begin{proposition}\label{prop:equivariantiso}
Let $F=F_{\un\la,A}$ and $G=F_{\un\mu,B}$ be two Kato germs. Let $\Phi:(\nC^n,0) \to (\nC^n,0)$ be a semi-conjugacy between $F$ and $G$ which is a local isomorphism on $T_N$ and so that $\Phi^{-1}(T_N) = T_N$.
 Then $\Phi$ is $\nu$-equivariant for some $\nu\in\Aut_{gr.}(T_N)$.
\end{proposition}
\begin{proof}
Being  $K:=\{z_1 \cdots z_n=0\}$ totally invariant for the action of $\Phi$, we infer that $\Phi$ takes the form
$$
\Phi(z)=\un{\chi} z^{Q}(\one + u(z))\text,
$$
where $\one=(1, \ldots, 1)$, $u:(\nC^n,0)\to (\nC^n,0)$, and $Q$ is a suitable matrix in $\GL(n,\nZ)$, with non-negative entries since $\Phi$ is regular.

From the conjugacy relation we get
\begin{equation}\label{eqn:conjugacypowerseries}
\un{\chi} \un{\lambda}^Q z^{QA} \Big(\one+u\big(\un{\lambda} z^A\big)\Big) = 
\un{\mu} \un{\chi}^B z^{BQ} \big(\one+u(z)\big)^B\text.
\end{equation}

By checking the lowest degrees of \eqref{eqn:conjugacypowerseries}, we get $QA=BQ$ and $\un{\chi} \un{\lambda}^Q=\un{\mu}\un{\chi}^B$.
After simplifying the common factor in \eqref{eqn:conjugacypowerseries}, we get
\begin{equation}\label{eqn:conjugacysimplified}
\one+u\big(\un{\lambda} z^A\big) = 
\big(\one+u(z)\big)^B\text.
\end{equation}

Let $u_i z^i$, with $u_i \in \nC^n$ and $i \in \nN^n \sminus \{0\}$, be any ($n$-uple of) monomial(s) appearing in the formal power series expansion of $u$.

This monomial contributes with $\un{\lambda}^i u_i z^{iA}$ in the left-hand-side, and with $B u_i z^{i}+\textrm{h.o.t.}$ in the right-hand-side of \eqref{eqn:conjugacysimplified}, where $\textrm{h.o.t.}$ denotes a suitable formal power series with monomials $z^j$ with $\abs{j}> \abs{i}$.

Since $\sum_{k} a_{hk} \geq 1$ for all $h=1, \ldots, n$, we have $\abs{iA} \geq \abs{i}$ for any $i \in \nN^n$.

Assume that this inequality is strict for any $i \in \nN^n \sminus \{0\}$. This implies that $u\equiv 0$.
In fact, if this is not the case, there exists $i \in \nN^n \sminus \{0\}$ so that $u_i \neq 0$. We may assume that $\abs{i}$ is minimal among the multi-indices satisfying the above condition.
In this case we get a contradiction, since the left-hand-side of \eqref{eqn:conjugacysimplified} contains a non-trivial monomial of order $\abs{i}$, while the right-hand side does not.

Assume now that $\abs{iA} = \abs{i}$ for some $i \in \nN^n \sminus\{0\}$.
Notice that $\sum_{k} a_{jk} = 1$ for some $j\in\{1,\ldots, n\}$ if and only if $P(A)^c = \{j\}$.
We infer that in this case $P(A)^c = \{j\}$ and $i=m e_j$ for some $m > 0$. Then the coefficient of the monomial $z^i$ in \eqref{eqn:conjugacysimplified} gives $\lambda_{j} u_i= B u_i$, where $\un{\lambda}=(\lambda_1, \ldots, \lambda_n)$.
Since $\abs{\lambda_j} < 1$ and $B$ has only eigenvalues of modulus $1$, we deduce that $u_i=0$, and hence $u\equiv 0$.
\end{proof}

\begin{example}\label{ex:notequivariantconj}
Consider a Kato $4$-fold $X=X(\pi, \sigma)$ whose associated Kato germ is $F=F_{\un\lambda,A}$, with $$A=\begin{pmatrix}
1&1&1&1	\\
0&1&1&1 \\
0&1&2&1 \\
0&1&1&2
\end{pmatrix}\text.
$$
For any $\alpha \in \nC$, let $H_\alpha$ be the hyperplane of equation $z_1-\alpha z_4=0$.
Then $H_\alpha$ is $F$-invariant if and only if $\lambda_1 = \lambda_4$. 
In this case, for any $\alpha \in \nC$, the map $\Phi(z)=(z_1-\alpha z_4, z_2, z_3, z_4)$ commutes with $F$.

Suppose moreover that $\Phi$ lifts to an automorphism $\hat{\Phi}:\hat\CC^n \to \hat\CC^n$ on the total space $\hat\CC^n$ of the toric modification $\pi$.
Then $\Phi$ induces a (non-equivariant) automorphism of $X$.

This example also shows that the vector space $H^0(X,\Theta_X)$ is generally larger than the vector space $H^0(X,\Theta_X(-\log D_T))\cong\ker(A-\id)$ computed in \ref{cohDT}.
\end{example}

Since the semi-conjugacy $\Phi$ induced by an equivariant isomorphism $\Psi$ as in \ref{thm:generaliso} automatically satisfies the hypothesis of \ref{prop:equivariantiso}, we deduce the following classification of toric Kato manifolds, associated to the same toric modification, up to equivariant isomorphisms.

\begin{corollary}\label{cor:condiso}
Let $X=X(\pi,\sigma)$ be a toric Kato manifold with Kato germ $F=F_{\un\la,A}$.
Let $X'=X(\pi,\sigma')$ be another toric Kato manifold associated to the same toric modification $\pi$, and denote by $G=F_{\un\mu,B}$ its associated Kato germ.
Set $\un\lambda = e^{2 \pi \ui \un\ell}$ and $\un\mu = e^{2 \pi \ui \un m}$, where $\un\ell$ and $\un m$ belong to $\CC^n$.

Then $X'$ is equivariantly isomorphic to $X$ if and only if there exists $Q \in \GL(n,\nZ)$ such that $QA=BQ$, $Q(\Sigma_A)=\Sigma_B$, and so that
\begin{equation}\label{eqn:condiso}
Q\un{\ell}-\un{m} \in \operatorname{Im}(B-\id) + \nZ^n\text. 
\end{equation}
\end{corollary}
\begin{proof}
By \ref{prop:equivariantiso}, if $X$ and $X'$ are equivariantly isomorphic, then there exists a map $\Phi(z)=\un\chi z^Q$ conjugating $F$ and $G$ and so that $Q\Sigma_A=\Sigma_B$. Conversely, any such map $\Phi$ induces a natural equivariant biholomorphism $\tilde\Phi:\tilde X\rightarrow \tilde X'$ commuting with the deck group action, and so an equivariant biholomorphism between $X$ and $X'$.

By $\Phi \circ F= G \circ \Phi$, we infer $QA=BQ$, and $\un{\chi} \un{\lambda}^{Q} = \un{\mu} \un{\chi}^{B}$.
By applying $\frac{1}{2 \pi \ui} \log$ to both sides, we get
$$
Q\un{\ell} + \un{v} = \un{m} + B \un{v} \ (\textrm{mod}\ \nZ^n)\text, 
$$
where $\un{v}$ is such that $\un{\chi}=e^{2 \pi \ui \un{v}}$.
The statement follows.
\end{proof}

Notice that, as in \ref{thm:generaliso}, we may assume that the matrix $Q$ in the statement of \ref{cor:condiso} has non-negative entries.
As an immediate corollary, we get:

\begin{corollary}
Let $X$ and $X'$ be two toric Kato manifolds as in \ref{cor:condiso}, both with the same associated Kato matrix $A$.
If $1 \not \in \Spec(A)$, then $X$ and $X'$ are equivariantly isomorphic.
\end{corollary}
\begin{proof}
It suffices to apply \ref{cor:condiso} with $A=B$ and $Q=\id$.
\end{proof}

The non-triviality of the moduli space of toric Kato varieties with prescribed toric modification is related to the existence of (possibly non-equivariant) invariant families of subvarieties (or more generally foliations) for the associated Kato germs, which in turn corresponds to the existence of families of compact subvarieties (or foliations) on the toric Kato variety.

\begin{example}(\textbf{Parabolic Inoue surfaces})
For $n=2$, the only toric Kato matrices $A$ with $\ker(\id-A)\neq 0$ are
$
A_a:=\begin{pmatrix}
1 &a\\
0 &1\end{pmatrix}
$
and
$A'_a=\begin{pmatrix}
1 &0\\
a &1\end{pmatrix}$,
which is $\GL_2(\ZZ)$-conjugated to $A_a$, where $a\in\NN^*$.
Let $\pi_0:\Bl_0\CC^2\rightarrow\CC^2$ be the blow-up of the origin, and $\sigma_0:\CC^2\rightarrow\Bl_0\CC^2$ be given by the chart $\tau:=\langle e_1, e_1+e_2\rangle$. Let now $(\pi,\sigma)$ be the composition of $(\pi_0,\sigma_0)$ with itself $a$ times, as given by \ref{powerKD}.

Being $e_1$ invariant, the locus $\{x=0\}$ is invariant by $F=F_{\un\la, A_a}$, and induces an elliptic curve $\nC^*/\lambda_2$ on the parabolic Inoue surface $X_{\un\la, A_a}=X(\pi,\sigma)$.
In particular if $X_{\un\la, A_a}$ and $X_{\un\mu, A_a}$ are isomorphic, then we must have $\lambda_2=\mu_2$.

A direct computation shows that the matrices $Q$ satisfying the conditions of \ref{cor:condiso} for $A=B=A_a$ are of the form $Q=A_k$ for some $k \in \nN$.
Then condition \eqref{eqn:condiso} is satisfied if and only if there exists $t \in \nR$ so that
$$
\begin{cases}
&\ell_1 + k \ell_2 - m_1 - k m_2 -t \in \nZ,\\
&\ell_2-m_2 \in \nZ,
\end{cases}
$$
which holds if and only if $\ell_2 - m_2 \in \nZ$, i.e., if and only if $\lambda_2=\mu_2$.
\end{example}

\begin{example}\label{ex:nontrivialiso3d}
Let $\pi$ be a toric modification whose associated fan contains the cone generated by $A=\begin{pmatrix}
1&1&1\\
2&2&1\\
2&1&2
\end{pmatrix}
$
Let $X=X(\pi, \sigma)$ be a toric Kato data whose associated Kato germ takes the form $F=F_{\un\la,A}$.
In this example, $1 \in \Spec(A)$, and $\operatorname{Im}(A-\id)$ is generated by $e_1+e_2+e_3$ and $e_1$, hence it coincides with $\{v \in \nR^3\ |\ v_2-v_3=0\}$.
In fact, an eigenvector for $1$ is given by $(0,1,-1)$, and can be also interpreted by the fact that $z_2/z_3=\text{const}$ defines an invariant family of surfaces for $F$.
These families induce a regular foliation on $X$, whose geometrical properties depend on the value of $\lambda_2/\lambda_3=:\lambda$.

Let now $X'=X(\pi,\sigma')$ be another toric Kato data with associated germ $G=F_{\un\mu,A}$.
The condition \eqref{eqn:condiso} to have $X$ and $X'$ isomorphic, when $Q=\id$, gives
$$
\ell_2-\ell_3-(m_2-m_3) \in \nZ.
$$
Hence $X$ and $X'$ are isomorphic if $\lambda=\mu:=\mu_2/\mu_3$.

Notice that the matrix $A$ commutes with the permutation matrix
$
Q=\begin{pmatrix}
1&0&0\\
0&0&1\\
0&1&0
\end{pmatrix}.
$
If the fan associated to $\pi$ is invariant by the action of $Q$ (see for example Figure \ref{fig:example3}), then $\Phi=\un{\chi} z^Q$ induces an isomorphism between $X$ and $X'$, as long as
$$
\ell_3-\ell_2-(m_2-m_3) \in \nZ.
$$
Hence in this case $X$ and $X'$ are isomorphic also whenever $\lambda\mu=1$.

\begin{tiny}
\begin{figure}
\begin{center}
\begin{minipage}[htbp]{0.5\columnwidth}
\def\svgwidth{0.95\columnwidth}
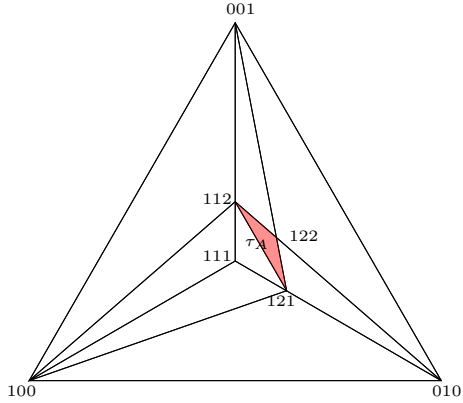
\end{minipage}
\end{center}
\caption{Fan associated to the toric modification $\pi$, invariant by the permutation $Q$ of \ref{ex:nontrivialiso3d}.}
\label{fig:example3}
\end{figure}
\end{tiny}
\end{example}

\subsection{Deformations of Kato manifolds of hyperbolic type}

\mbox{}

\begin{theorem}
Let $X=X(\pi,\sigma)$ be a toric Kato manifold of non-parabolic type, with toric divisor $D_T$. There exists a smooth versal family of deformations of the pair $p:(\mathcal X,\mathcal D)\rightarrow B$, where $B$ is a centered ball in $H^1(X,\Theta_X(-\log D_T))\cong\ker(A-\id)\subset N\otimes\CC$. Moreover, for each $u\in B$,  $(X_u,\mathcal D\cap X_u)$ is the toric Kato manifold $X(\pi,(\exp u).\sigma)$ together with its toric divisor. In particular, the family $p$ is versal for each of its fibers.
\end{theorem}
\begin{proof}
The existence of a smooth versal family of deformations of the pair follows simply from the fact that $H^2(X,\Theta_X(-\log D_T))=0$ and \cite{kns}, but we can in fact construct such a family explicitly. Let $A$ be a Kato matrix of $X$. Note first that since
\begin{equation*}
H^1(X,\Theta_X(-\log D_T))=\coker(A-\id)=H^1(\Gamma,H^0(\wt{X},\Theta_{\wt{X}}(-\log\wt{D}_T)))
\end{equation*}
the family $(\mathcal X,\mathcal D)$ must come from a trivial family 
\begin{equation*}
\wt{p}:(\wt{\mathcal X}=\wt{X}\times B,\wt{\mathcal D}=\wt{D}_T\times B)\rightarrow B
\end{equation*}
and a deformation of the action of $\Gamma$. Indeed, let $\gamma$ denote the positive generator of the deck group of $\wt{X}$, and consider an action of $\Gamma:=\langle \wt{\gamma}\rangle$ on $\wt{\mathcal X}$ by 
\begin{equation*}
\wt{\gamma}(z,u)=(\gamma(\exp u.z),u), \ \ (z,u)\in \wt{X}\times \ker(A-\id)
\end{equation*}
where $\exp:N\otimes\CC\rightarrow T_N=N\otimes\CC^*$ is induced by $\CC\ni t\mapsto \e^{2\pi it}\in\CC^*$. Since  $\gamma$ commutes with $\exp v$ for any $v\in\ker(A-\id)$, $\wt{X}\subset X(\Sigma_A,N)$ is invariant under the action of $\exp(\ker(A-\id))$, so $\wt{\gamma}$ is well-defined. Moreover, we take the ball $B\subset  \ker(A-\id)$ small enough so that $\sigma(\exp u.\ov{\BB})\subset\BB$ for  each $u\in B$. It follows that for each $u\in B$, 
\begin{equation*}
\wt{\gamma}|_{\wt{X}\times \{u\}}=\gamma(\exp u\cdot)=\exp u\gamma\in \Aut(\wt{X}).
\end{equation*}
Therefore $\Gamma$ acts freely and properly on $(\wt{\mathcal X},\wt{\mathcal D})$ and $p:(\mathcal X:=\wt{\mathcal X}/\Gamma, \mathcal D:=\wt{\mathcal D}/\Gamma)\rightarrow B$ gives the desired smooth versal family. 
\end{proof}

\section{Hermitian geometry of Kato manifolds}\label{metrici}

Kato manifolds cannot admit K\" ahler metrics since their fundamental group is isomorphic to $\mathbb{Z}$. However, a large class of them carry {\it locally conformally K\" ahler metrics (lcK)}, as shown in \cite[Theorem 2.2]{iop}. In this section, we give a characterization for the existence of lcK metrics on Kato manifolds. Furthermore, we investigate the possibility of endowing a Kato manifold with other special Hermitian metrics,  such as balanced, pluriclosed,  Hermitian symplectic or strongly Gauduchon. 

Recall the following definition:
\begin{definition}\label{deflck} A {\it locally conformally K\" ahler metric (lcK)} on a complex manifold $X$ is a Hermitian metric $\Omega$ satisfying $d\Omega=\theta\wedge\Omega$ for a closed one-form $\theta$. If $\theta$ is not exact, $\Omega$ is called strictly lcK. 
\end{definition}

This definition is equivalent to any of the following formulations:
\begin{enumerate}
\item There exists a covering with open sets $(U_i)_i$ of $X$ and smooth functions $f_i$ on every $U_i$ such that $e^{-f_i}\Omega$ is K\" ahler.
\item The universal cover $\wt{X}$ admits a K\" ahler metric $\wt{\Omega}$ such that $\gamma^*\wt{\Omega}=c_{\gamma}\wt{\Omega}$, $c_{\gamma}>0$, for any deck-transformation $\gamma$.
\end{enumerate}
 
For our characterization purpose, we start by giving a general principle used to modify K\" ahler metrics on modifications of balls, already used in \cite{bru} and \cite{iop}, which will be needed.

\begin{lemma}\label{extKA}
Let $\pi:\hat\CC^n\rightarrow \CC^n$ be a proper modification at $0$, let $\rho(z)=\rho(||z||)\in\ce(\CC^n,\RR)$ be a strictly plurisubharmonic function depending only on $||z||$, and let $\omega$ be a \Ka\ metric on $\hat{\ov\BB}:=\pi^{-1}(\ov\BB)$. Then for any $0<s<1$, there exists $\la>0$ and a \Ka\ metric $\wt \omega$ on $\hat\CC^n$ so that 
\begin{equation*}
\wt{\omega}|_{\hat\BB_s}=\omega, \ \ \  \wt{\omega}|_{\hat\CC^n\sminus \hat\BB}=\la\cdot\pi^*dd^c\rho.
\end{equation*}
\end{lemma}
\begin{proof}
There exists a pluriharmonic function $\psi$ on $\BB$, smooth and strictly plurisubharmonic on $\BB\sminus \{0\}$, so that $\pi_*\omega=dd^c\psi$. By the maximum principle, the functions
\begin{equation*}
u,v:(0,1]\rightarrow \RR, \ \ u(t)=\max_{\del\BB_t}\rho, \ \ v(t)=\max_{\del\BB_t}\psi
\end{equation*}
are strictly increasing. It follows that for any $0<s<1$, 
\begin{gather*}
r:=u(s)\equiv \rho|_{\del\BB_s}<R:=u(1)\equiv \rho|_{\del\BB}, \ \ m:=\min_{\del\BB_s}\psi\leq v(s)<M:=v(1).
\end{gather*}
Thus, for any $\al>1$, we can take 
\begin{equation*}
\la:=\al\frac{M-m}{R-r}>0, \ \ c\in\left(M-\la R,m-\la r\right)\neq\emptyset.
\end{equation*}
so that $\la\rho+c>\psi$ on $\del\BB$ and $\la\rho+c<\psi$ on $\del\BB_s$. Therefore, if $\wt{\psi}$ is the function on $\CC^n$ defined as the regularized maximum between $\psi$ and $\la\rho+c$ on $\BB\sminus \BB_s$, equal to $\psi$ on $\BB_s$ and equal to $\la\rho+c$ on $\CC^n\sminus \BB$,  then the metric $\pi^*(dd^c\wt{\psi})$ on $\hat\CC^n\sminus \pi^{-1}(0)$ glues to $\omega$ to define the desired metric $\wt{\omega}$ on $\hat\CC^n$.
\end{proof}

Now we are ready to characterize the existence of lcK metrics. For this, let us note that any proper modification $\pi:\hat\BB\rightarrow\BB$ induces naturally proper modifications $\pi:\hat\CC^n\rightarrow\CC^n$ and $\pi:\hat{\CC\PP}^n\rightarrow\CC\PP^n$, where $\hat{\CC\PP}^n$ is the compactification of $\hat\CC^n$ with a hyperplane at infinity.

 \begin{theorem}\label{lcK}	
 Let $(\pi:\hat\BB\rightarrow\BB,\sigma:\overline{\BB}\rightarrow\hat\BB)$ be a Kato data and let $X=X(\pi,\sigma)$ be the corresponding Kato manifold. 
	 The following are equivalent:
	\begin{enumerate}
		\item $X$ admits an lcK metric;
		\item $\wt{X}$ admits a \Ka\ metric;
		\item $\hat{\CC\PP}^n$ is a  projective manifold;
		\item $\hat{\mathbb{C}}^n$ admits a K\" ahler metric;
		
	\end{enumerate}
\end{theorem}

\begin{proof}
The implications $(1)\Rightarrow (2)$ and $(3)\Rightarrow (4)$ are clear, while the implication $(4) \Rightarrow (1)$ is precisely Brunella's theorem \cite{bru} adapted to all complex dimensions \cite[Theorem 2.2]{iop}. Indeed, in \cite[Theorem 2.2]{iop} we showed that if $\pi$ is a composition of smooth blow-ups, then $X$ is lcK. The same proof works if one merely supposes that $\hat\BB$, or equivalently $\hat\CC^n$, is K\" ahler. We are thus left with showing $(2) \Rightarrow (3)$. 

If $\wt{X}$ is K\" ahler, then \ref{embBall} implies that $\hat\BB\sminus \{\sigma(0)\}$ is also \Ka, so by Miyaoka's extension theorem \cite[Proposition A]{miyaoka} $\hat\BB$ admits a K\" ahler metric $\omega$. Now putting  $\rho:=\log(1+||z||^2)$ in \ref{extKA}, we find that $\hat\CC^n$ admits a \Ka\ metric $\wt{\omega}$ which glues to $\la\cdot\pi^*\omega_{FS}$ on $\hat{\CC\PP^n}\sminus \hat\BB$ for some $\la>0$, where $\omega_{FS}$ is the standard Fubini-Study metric on $\CC\PP^n$. Thus we have shown that $\hat{\CC\PP^n}$ is \Ka. Since on the other hand $\hat{\CC\PP^n}$ is a modification of $\CC\PP^n$, it is Moishezon, and therefore $\hat{\CC\PP^n}$ is projective by  Moishezon's theorem \cite[Chapter I, Theorem 11]{moi}. This concludes the proof.
	\end{proof}
	
\begin{remark}
We note here that a non-Hopf toric Kato manifold $X$ admitting an lcK metric can in no way be a toric lcK manifold in the sense of \cite[Definition 4.7]{is}, because in this case it would admit a Vaisman metric by \cite[Theorem A]{is}, which is impossible by \cite[Proposition~2.6]{iop}. At the same time, the compact torus $\TT$ acts effectively and holomorphically on the universal cover $\wt{X}$, and any lcK metric on $X$ gives rise to a \Ka\ metric $\omega$ on $\wt{X}$, which can be made $\TT$-invariant by averaging. Since $\wt{X}$ is moreover simply connected, it follows that $(\wt{X},\omega,\TT)$ is a toric \Ka\ manifold in the classical sense. In this manner, lcK toric Kato manifolds give interesting generalizations of the class of toric lcK manifolds. 
\end{remark}	


\hfill
Next, we recall the definitions of the other classes  of Hermitian metrics whose existence on Kato manifolds we discuss. 

\begin{definition} Let $\Omega$ be a Hermitian metric on a complex manifold $X$. $\Omega$ is called:
\begin{itemize}
\item {\it balanced} if $d\Omega^{n-1}=0$, or equivalently, if $\Omega$ is co-closed \cite{mic}.
\item  {\it pluriclosed}, or {\it strongly K\" ahler with torsion}, if $\partial \overline{\partial}\Omega=0$ \cite{bis}.
\item {\it Hermitian symplectic} if $\Omega$ is the $(1,1)$-component of a real $d$-closed two-form \cite{st}.
\item  {\it strongly Gauduchon} if $\partial\Omega^{n-1}$ is $\overline{\partial}$-exact \cite{P9}.

\end{itemize}
\end{definition}

There are certain inclusion relations between the classes of manifolds which admit one of the above metrics. For instance, it is easy to see that any balanced metric is strongly Gauduchon. Similarly, it is easy to see that any Hermitian symplectic metric is pluriclosed. Finally, \cite[Lemma~1]{yzz} shows that any manifold admitting a Hermitian symplectic metric also admits a strongly Gauduchon one.

We have the following non-existence results:
\begin{theorem}\label{balanced}
A Kato manifold $X$ admits no strongly Gauduchon metric, and in particular no balanced or Hermitian symplectic metric.
\end{theorem}

\begin{proof}
	Let $\mathcal{X}\xrightarrow{\pi} \mathbb{D}$ be the complex analytic deformation of $X$ given by  \cite[Theorem 1]{kato}, so that $\pi^{-1}(0)=X$ and $X_t:=\pi^{-1}(t)$ is a proper modification at a finite number of points of a Hopf manifold for $t\neq 0$. Assume that $X$ admits a strongly Gauduchon metric. Then by \cite[Theorem~3.1]{P14}, stating that the strongly Gauduchon property is open with respect to holomorphic deformations, $X_t$ also admits a strongly Gauduchon metric if $t$ is close to $0$. However, according to \cite[Theorem~2.2]{P10b}, if $\mu: \hat{Y} \rightarrow Y$ is a proper modification of complex manifolds and $\hat{Y}$ carries a strongly Gauduchon metric, then so does $Y$. Consequently, this would imply that Hopf manifolds admit strongly Gauduchon metrics, which is false by \cite[Section~2]{P14}. Therefore, $X$ cannot be endowed with strongly Gauduchon metrics, nor with balanced or Hermitian symplectic metrics \cite[Lemma~1]{yzz}. 
\end{proof}

\begin{remark}
The above result is related to the open question \cite[Question~1.7]{st}, asking whether there exist any non-\Ka\ Hermitian symplectic manifolds in dimension $\geq 3$. We answer it negatively for the class of Kato manifolds. 
\end{remark}

\begin{theorem}\label{pluri}
	If $X$ is a Kato manifold of dimension $n\geq 3$ satisfying $H^{1, 2}_{\overline{\partial}}(X)=0$, then it cannot be endowed with a pluriclosed metric. In particular, if $X$ is a toric Kato manifold of hyperbolic type, then it does not admit a pluriclosed metric unless it is a surface. 
\end{theorem}
\begin{proof}
Suppose that $X$ carries a pluriclosed metric. Let $\mathcal{X} \xrightarrow{\pi} \mathbb{D} \subset \mathbb{C}$ be the same deformation as above. Then by  \cite[Theorem 8.5]{c}, for any $t \in\DD$ close to $0$, $\pi^{-1}(t)$ also admits a pluriclosed metric. We recall that for each $t\neq 0$, $\pi^{-1}(t)$ is a modification at a finite number of points of a Hopf manifold $M_t$. According to the analogue of Miyaoka's theorem for pluriclosed metrics \cite[Theorem 4.3]{ft}, a complex manifold $M$ without a point admits a pluriclosed metric if and only if $M$ admits a pluriclosed metric. Therefore, the Hopf manifolds $M_t$ also admit pluriclosed metrics. However, by \cite[Example 5.17, Theorem 5.16]{c}, this can only happen if $n=2$. Finally, in the toric case the conclusion follows by \ref{cohOX2}. 
\end{proof}

\begin{remark} In \ref{balanced} and \ref{pluri} we referred to \cite{P14} and \cite{c} for the proofs that Hopf manifolds do not admit strongly Gauduchon/pluriclosed metrics. However, the proofs therein are given for diagonal Hopf manifolds. Nevertheless, they hold for any Hopf manifold, since they depend only on the existence of complex hypersurfaces for strongly Gauduchon metrics and the vanishing of the second Betti number and Hodge numbers $h^{3, 0}$ and $h^{2,1}$ for pluriclosed metrics. These are general facts that hold for any Hopf manifold, see \cite[Section 2]{kato79} and \cite[Theorem 3]{mall91}.
\end{remark}

\subsection*{Examples of non-lcK toric Kato manifolds}

We have shown in \cite{iop} that a Kato manifold of dimension $n\geq 4$ does not need to be lcK, however the examples we constructed there are not toric. We end this section by constructing examples of toric Kato manifolds in any dimension $n\geq 4$ which admit no lcK metric.

In order to do so, we will start from the simplest example of a complete smooth toric algebraic variety which admits no \Ka\ metric, cf. \cite{oda2}

Consider the toric manifold $\CC\PP^3$ with the standard action of $T_{N}=(\CC^*)^3$, $N=\ZZ^3$. Let $e_1,e_2, e_3$ denote the standard basis of $N$ and let $e_0=-e_1-e_2-e_3$. Then the fan of $\CC\PP^3$ is generated by the $3$-dimensional cones
\begin{gather*}
\sigma_0=\langle e_1, e_2, e_3  \rangle, \ \ \sigma_1=\langle e_0, e_2, e_3  \rangle, \ \ \sigma_2=\langle e_0, e_1,  e_3  \rangle, \ \ \sigma_3=\langle e_0, e_1, e_2 \rangle.
\end{gather*}
Consider the $T_N$-fixed points $P_j=\ov{\orb \sigma_j}$, $j=1,2,3$, and the one-dimensional orbits $d_{0j}=\ov{\orb\langle e_0, e_j\rangle}$, $j=1,2,3$, so that we have $d_{01}\cap d_{02}=P_3$ and so on.

Let $\mu:\hat{\CC\PP}^3\rightarrow \CC\PP^3$ be the toric modification of $\CC\PP^3$ constructed as follows. Around the point $P_1$, blow-up $d_{02}$ and then the strict transform of $d_{03}$, around $P_2$ blow-up $d_{03}$ and then the strict transform of $d_{01}$, and around the point $P_3$ blow up $d_{01}$ and then the strict transform of $d_{02}$. It is easy to check that this is a well defined operation globally, resulting in a new smooth toric variety $\hat{\CC\PP}^3$. The corresponding fan $\Sigma_3$ has $3$ new rays generated by $v_j=e_0+e_j$, $j=1,2,3$, and the following $3$-dimensional cones:
\begin{gather*}
\langle e_1, v_1, e_2  \rangle,  \ \ \langle v_1, e_2, v_2 \rangle,\ \ \langle v_1,v_2,e_0  \rangle \\
\langle e_2, e_3, v_2 \rangle, \ \ \langle v_2, e_3, v_3 \rangle, \ \ \langle v_2, v_3, e_0  \rangle\\
\langle e_1, e_3, v_3 \rangle, \ \ \langle e_1, v_3, v_1  \rangle, \ \ \langle v_1,v_3, e_0  \rangle\\
\langle e_1, e_2, e_3  \rangle. 
\end{gather*}

\begin{proposition}\label{nonproj3} The toric variety $\hat{\CC\PP}^3$ admits no \Ka\ metric.  
\end{proposition}
\begin{proof}
In \cite[Proposition 9.4]{oda} or \cite[Example on page 84]{oda2} it is shown that $\hat{\CC\PP}^3$ is not projective. However, being a toric variety, $\hat{\CC\PP}^3$ has $3$ algebraically independent meromorphic functions, i.e. it is a Moishezon manifold. Therefore, by Moishezon's theorem \cite[Chapter I, Theorem 11]{moi}, $\hat{\CC\PP}^3$ is non-\Ka.
\end{proof}

We now show:

\begin{proposition}
In any complex dimension $n\geq 4$, there exist toric Kato manifolds $X$ with holomorphic immersions $\hat{\CC\PP}^3\rightarrow X$. In particular, $X$ admits no lcK metric.
\end{proposition}
\begin{proof}
We wish first to construct a toric modification $\pi:\hat\CC^n\rightarrow\CC^n$ over $0$ which contains $\hat{\CC\PP}^3$ as a submanifold. Let $\pi'':\wt{\CC}^n\rightarrow\CC^n$ be the blow-up at $0$ of $\CC^n$, let $N=\ZZ^n$ of standard basis $f_1,\ldots, f_n$, let $f_0=\sum_{j=1}^nf_j$ and denote by $\wt{\Sigma}$ the fan of $\wt{\CC}^n$ with the standard $T_N$-action. Its $n$-dimensional cones are:
\begin{equation*}
\tau_j:=\langle f_0,\ldots, \widehat{f_j}, \ldots, f_n\rangle, \ \ j=0,1\ldots n.
\end{equation*}

Let $E=\ov{\orb\langle f_0\rangle}\cong \CC\PP^{n-1}$ be the exceptional divisor of $\wt{\CC}^n$ and identify 
\begin{equation*}
\CC\PP^3=\ov{\orb\langle f_0,f_5\ldots, f_n\rangle}. 
\end{equation*}
For $j=1,2,3$ denote by $L_j=\ov{\orb\langle f_0,f_j,f_4,f_5\ldots, f_n\rangle}$ the $T_N$-invariant one-dimensional submanifolds. Again, put $P_j=\ov{\orb\tau_j}$, $j=1,2,3$. Define the toric modification 
$\pi':\hat\CC^n\rightarrow \wt{\CC}^n$ similarly as before. Namely, around $P_1$, blow up $L_{2}$ and then the strict transform of $L_{3}$, around $P_2$ blow-up $L_{3}$ and then the strict transform of $L_{1}$ and around the point $P_3$ blow up $L_{1}$ and then the strict transform of $L_{2}$.  Equivalently, define $\hat\CC^n=X(\hat{\Sigma},N)$, where the fan $\hat{\Sigma}$ has the rays of $\wt{\Sigma}$ to which we add $\nu_j=f_0+f_j+\sum_{k=4}^nf_k$, $j=1,2,3$,  and is generated by the $n$-dimensional cones:
\begin{gather*}
\tau_0, \ \ \tau_j, \ \ 4\leq j\leq n\\
\tau_{1k}=\langle \nu_2, f_s\mid s\in\{0,\ldots, n\}\sminus \{1,k\}\rangle, \ \ k\in\{0,4,\ldots, n\}\\
\tau_{12l}=\langle\nu_2,\nu_3, f_s\mid s\in\{0,\ldots, n\}\sminus \{1,2,l\}\rangle, \ \ l\in\{0,3,4,\ldots, n\}\\
\tau_{2k}=\langle \nu_3, f_s\mid s\in\{0,\ldots, n\}\sminus \{2,k\}\rangle, \ \ k\in\{0,4,\ldots, n\}\\
\tau_{23l}=\langle\nu_3,\nu_1, f_s\mid s\in\{0,\ldots, n\}\sminus \{2,3,l\}\rangle, \ \ l\in\{0,1,4,\ldots, n\}\\
\tau_{3k}=\langle \nu_1, f_s\mid s\in\{0,\ldots, n\}\sminus \{3,k\}\rangle, \ \ k\in\{0,4,\ldots, n\}\\
\tau_{31l}=\langle \nu_1,\nu_2, f_s\mid s\in\{0,\ldots, n\}\sminus \{1,3,l\}\rangle, \ \ l\in\{0,2,4,\ldots, n\}.
\end{gather*}

It is clear that the strict transform of $\CC\PP^3$ under $\pi'$ is isomorphic to $\hat{\CC\PP}^3$, so $\pi=\pi''\circ \pi':\hat\CC^n\rightarrow\CC^n$ has the desired properties. Now it suffices to take any toric chart $\sigma$ with $\sigma(\BB)\subset(\hat\BB)$ and $\sigma(0)\notin\hat{\CC\PP}^3$, for instance one corresponding to the cone $\tau_{10}$, in order to obtain a toric Kato manifold $X=X(\pi,\sigma)$. Clearly $\hat{\CC\PP}^3\subset \wt{X}$, therefore $\wt{X}$ is not \Ka\ and so $X$ admits no lcK metric.
\end{proof} 

{\noindent{\bf Aknowledgements:}} N. Istrati would like to thank Tel Aviv University for the hospitality, where part of this project was carried out. A. Otiman and M. Pontecorvo are partially supported by GNSAGA and PRIN 2017 Real and Complex Manifolds: Topology, geometry and holomorphic dynamics.
M. Ruggiero is partially supported by the ANR grant Fatou ANR-17-CE40-0002-01.
The first three named authors would also like to thank Victor Vuletescu for many stimulating discussions.

\end{document}

%% file: example4a.pdf_tex
\begingroup%
  \makeatletter%
  \providecommand\color[2][]{%
    \errmessage{(Inkscape) Color is used for the text in Inkscape, but the package 'color.sty' is not loaded}%
    \renewcommand\color[2][]{}%
  }%
  \providecommand\transparent[1]{%
    \errmessage{(Inkscape) Transparency is used (non-zero) for the text in Inkscape, but the package 'transparent.sty' is not loaded}%
    \renewcommand\transparent[1]{}%
  }%
  \providecommand\rotatebox[2]{#2}%
  \newcommand*\fsize{\dimexpr\f@size pt\relax}%
  \newcommand*\lineheight[1]{\fontsize{\fsize}{#1\fsize}\selectfont}%
  \ifx\svgwidth\undefined%
    \setlength{\unitlength}{341.6678772bp}%
    \ifx\svgscale\undefined%
      \relax%
    \else%
      \setlength{\unitlength}{\unitlength * \real{\svgscale}}%
    \fi%
  \else%
    \setlength{\unitlength}{\svgwidth}%
  \fi%
  \global\let\svgwidth\undefined%
  \global\let\svgscale\undefined%
  \makeatother%
  \begin{picture}(1,0.89059585)%
    \lineheight{1}%
    \setlength\tabcolsep{0pt}%
    \put(-0.2230795,-2.14191836){\color[rgb]{0,0,0}\makebox(0,0)[lt]{\begin{minipage}{0.24889485\unitlength}\raggedright \end{minipage}}}%
    \put(0,0){\includegraphics[width=\unitlength,page=1]{example4a.pdf}}%
    \put(0.37587807,0.27802812){\color[rgb]{0,0,0}\makebox(0,0)[lt]{\begin{minipage}{0.21467115\unitlength}\raggedright $\tau_A$\\ \end{minipage}}}%
    \put(0.51197516,0.35266194){\color[rgb]{0,0,0}\makebox(0,0)[lt]{\begin{minipage}{0.12073128\unitlength}\raggedright $111$\end{minipage}}}%
    \put(0.32758553,0.26924761){\color[rgb]{0,0,0}\makebox(0,0)[lt]{\begin{minipage}{0.12073128\unitlength}\raggedright $211$\end{minipage}}}%
    \put(0.41539023,0.1682723){\color[rgb]{0,0,0}\makebox(0,0)[lt]{\begin{minipage}{0.12073128\unitlength}\raggedright $321$\end{minipage}}}%
    \put(0,0){\includegraphics[width=\unitlength,page=2]{example4a.pdf}}%
    \put(0.34075619,0.39656423){\color[rgb]{0,0,0}\makebox(0,0)[lt]{\begin{minipage}{0.12073128\unitlength}\raggedright $212$\end{minipage}}}%
    \put(0.01587916,0.04973606){\color[rgb]{0,0,0}\makebox(0,0)[lt]{\begin{minipage}{0.16198862\unitlength}\raggedright $1 0 0$\end{minipage}}}%
    \put(0.9202671,0.04973606){\color[rgb]{0,0,0}\makebox(0,0)[lt]{\begin{minipage}{0.12073128\unitlength}\raggedright $010$\end{minipage}}}%
    \put(0.48124359,0.87070923){\color[rgb]{0,0,0}\makebox(0,0)[lt]{\begin{minipage}{0.12073128\unitlength}\raggedright $001$\end{minipage}}}%
    \put(0.46368266,0.29558897){\color[rgb]{0,0,0}\makebox(0,0)[lt]{\begin{minipage}{0.12073128\unitlength}\raggedright $322$\end{minipage}}}%
    \put(0.49880457,0.19900389){\color[rgb]{0,0,0}\makebox(0,0)[lt]{\begin{minipage}{0.12073128\unitlength}\raggedright $221$\end{minipage}}}%
    \put(0.44173159,0.22534527){\color[rgb]{0,0,0}\makebox(0,0)[lt]{\begin{minipage}{0.12073128\unitlength}\raggedright $432$\end{minipage}}}%
  \end{picture}%
\endgroup%

%% file: example4b.pdf_tex
\begingroup%
  \makeatletter%
  \providecommand\color[2][]{%
    \errmessage{(Inkscape) Color is used for the text in Inkscape, but the package 'color.sty' is not loaded}%
    \renewcommand\color[2][]{}%
  }%
  \providecommand\transparent[1]{%
    \errmessage{(Inkscape) Transparency is used (non-zero) for the text in Inkscape, but the package 'transparent.sty' is not loaded}%
    \renewcommand\transparent[1]{}%
  }%
  \providecommand\rotatebox[2]{#2}%
  \newcommand*\fsize{\dimexpr\f@size pt\relax}%
  \newcommand*\lineheight[1]{\fontsize{\fsize}{#1\fsize}\selectfont}%
  \ifx\svgwidth\undefined%
    \setlength{\unitlength}{341.6678772bp}%
    \ifx\svgscale\undefined%
      \relax%
    \else%
      \setlength{\unitlength}{\unitlength * \real{\svgscale}}%
    \fi%
  \else%
    \setlength{\unitlength}{\svgwidth}%
  \fi%
  \global\let\svgwidth\undefined%
  \global\let\svgscale\undefined%
  \makeatother%
  \begin{picture}(1,0.89059585)%
    \lineheight{1}%
    \setlength\tabcolsep{0pt}%
    \put(-0.2230795,-2.14191836){\color[rgb]{0,0,0}\makebox(0,0)[lt]{\begin{minipage}{0.24889485\unitlength}\raggedright \end{minipage}}}%
    \put(0,0){\includegraphics[width=\unitlength,page=1]{example4b.pdf}}%
    \put(0.37587807,0.27802812){\color[rgb]{0,0,0}\makebox(0,0)[lt]{\begin{minipage}{0.21467115\unitlength}\raggedright $\tau_A$\\ \end{minipage}}}%
    \put(0.51197516,0.35266194){\color[rgb]{0,0,0}\makebox(0,0)[lt]{\begin{minipage}{0.12073128\unitlength}\raggedright $111$\end{minipage}}}%
    \put(0.32758553,0.26924761){\color[rgb]{0,0,0}\makebox(0,0)[lt]{\begin{minipage}{0.12073128\unitlength}\raggedright $211$\end{minipage}}}%
    \put(0.41539023,0.1682723){\color[rgb]{0,0,0}\makebox(0,0)[lt]{\begin{minipage}{0.12073128\unitlength}\raggedright $321$\end{minipage}}}%
    \put(0,0){\includegraphics[width=\unitlength,page=2]{example4b.pdf}}%
    \put(0.34075619,0.39656423){\color[rgb]{0,0,0}\makebox(0,0)[lt]{\begin{minipage}{0.12073128\unitlength}\raggedright $212$\end{minipage}}}%
    \put(0.01587916,0.04973606){\color[rgb]{0,0,0}\makebox(0,0)[lt]{\begin{minipage}{0.16198862\unitlength}\raggedright $1 0 0$\end{minipage}}}%
    \put(0.9202671,0.04973606){\color[rgb]{0,0,0}\makebox(0,0)[lt]{\begin{minipage}{0.12073128\unitlength}\raggedright $010$\end{minipage}}}%
    \put(0.48124359,0.87070923){\color[rgb]{0,0,0}\makebox(0,0)[lt]{\begin{minipage}{0.12073128\unitlength}\raggedright $001$\end{minipage}}}%
    \put(0.44612174,0.29119874){\color[rgb]{0,0,0}\makebox(0,0)[lt]{\begin{minipage}{0.12073128\unitlength}\raggedright $322$\end{minipage}}}%
    \put(0.51197526,0.23851594){\color[rgb]{0,0,0}\makebox(0,0)[lt]{\begin{minipage}{0.12073128\unitlength}\raggedright $221$\end{minipage}}}%
    \put(0,0){\includegraphics[width=\unitlength,page=3]{example4b.pdf}}%
  \end{picture}%
\endgroup%

%% file: example1a.pdf_tex
\begingroup%
  \makeatletter%
  \providecommand\color[2][]{%
    \errmessage{(Inkscape) Color is used for the text in Inkscape, but the package 'color.sty' is not loaded}%
    \renewcommand\color[2][]{}%
  }%
  \providecommand\transparent[1]{%
    \errmessage{(Inkscape) Transparency is used (non-zero) for the text in Inkscape, but the package 'transparent.sty' is not loaded}%
    \renewcommand\transparent[1]{}%
  }%
  \providecommand\rotatebox[2]{#2}%
  \newcommand*\fsize{\dimexpr\f@size pt\relax}%
  \newcommand*\lineheight[1]{\fontsize{\fsize}{#1\fsize}\selectfont}%
  \ifx\svgwidth\undefined%
    \setlength{\unitlength}{341.6678772bp}%
    \ifx\svgscale\undefined%
      \relax%
    \else%
      \setlength{\unitlength}{\unitlength * \real{\svgscale}}%
    \fi%
  \else%
    \setlength{\unitlength}{\svgwidth}%
  \fi%
  \global\let\svgwidth\undefined%
  \global\let\svgscale\undefined%
  \makeatother%
  \begin{picture}(1,0.89059585)%
    \lineheight{1}%
    \setlength\tabcolsep{0pt}%
    \put(0,0){\includegraphics[width=\unitlength,page=1]{example1a.pdf}}%
    \put(-0.2230795,-2.14191836){\color[rgb]{0,0,0}\makebox(0,0)[lt]{\begin{minipage}{0.24889485\unitlength}\raggedright \end{minipage}}}%
    \put(0,0){\includegraphics[width=\unitlength,page=1]{example1a.pdf}}%
    \put(0.01587916,0.04973605){\color[rgb]{0,0,0}\makebox(0,0)[lt]{\begin{minipage}{0.16198862\unitlength}\raggedright $1 0 0$\end{minipage}}}%
    \put(0.92026711,0.04973605){\color[rgb]{0,0,0}\makebox(0,0)[lt]{\begin{minipage}{0.12073128\unitlength}\raggedright $010$\end{minipage}}}%
    \put(0.4812436,0.87070922){\color[rgb]{0,0,0}\makebox(0,0)[lt]{\begin{minipage}{0.12073128\unitlength}\raggedright $001$\end{minipage}}}%
    \put(0.43295105,0.33949125){\color[rgb]{0,0,0}\makebox(0,0)[lt]{\begin{minipage}{0.12073128\unitlength}\raggedright $111$\end{minipage}}}%
    \put(0.56904822,0.24290617){\color[rgb]{0,0,0}\makebox(0,0)[lt]{\begin{minipage}{0.12073128\unitlength}\raggedright $121$\end{minipage}}}%
    \put(0.46807285,0.27802812){\color[rgb]{0,0,0}\makebox(0,0)[lt]{\begin{minipage}{0.21467115\unitlength}\raggedright $\tau_A$\\ \end{minipage}}}%
  \end{picture}%
\endgroup%

%% file: example1b.pdf_tex
\begingroup%
  \makeatletter%
  \providecommand\color[2][]{%
    \errmessage{(Inkscape) Color is used for the text in Inkscape, but the package 'color.sty' is not loaded}%
    \renewcommand\color[2][]{}%
  }%
  \providecommand\transparent[1]{%
    \errmessage{(Inkscape) Transparency is used (non-zero) for the text in Inkscape, but the package 'transparent.sty' is not loaded}%
    \renewcommand\transparent[1]{}%
  }%
  \providecommand\rotatebox[2]{#2}%
  \newcommand*\fsize{\dimexpr\f@size pt\relax}%
  \newcommand*\lineheight[1]{\fontsize{\fsize}{#1\fsize}\selectfont}%
  \ifx\svgwidth\undefined%
    \setlength{\unitlength}{341.6678772bp}%
    \ifx\svgscale\undefined%
      \relax%
    \else%
      \setlength{\unitlength}{\unitlength * \real{\svgscale}}%
    \fi%
  \else%
    \setlength{\unitlength}{\svgwidth}%
  \fi%
  \global\let\svgwidth\undefined%
  \global\let\svgscale\undefined%
  \makeatother%
  \begin{picture}(1,0.89059585)%
    \lineheight{1}%
    \setlength\tabcolsep{0pt}%
    \put(0,0){\includegraphics[width=\unitlength,page=1]{example1b.pdf}}%
    \put(-0.2230795,-2.14191836){\color[rgb]{0,0,0}\makebox(0,0)[lt]{\begin{minipage}{0.24889485\unitlength}\raggedright \end{minipage}}}%
    \put(0,0){\includegraphics[width=\unitlength,page=1]{example1b.pdf}}%
    \put(0.01587916,0.04973605){\color[rgb]{0,0,0}\makebox(0,0)[lt]{\begin{minipage}{0.16198862\unitlength}\raggedright $1 0 0$\end{minipage}}}%
    \put(0.92026711,0.04973605){\color[rgb]{0,0,0}\makebox(0,0)[lt]{\begin{minipage}{0.12073128\unitlength}\raggedright $010$\end{minipage}}}%
    \put(0.4812436,0.87070922){\color[rgb]{0,0,0}\makebox(0,0)[lt]{\begin{minipage}{0.12073128\unitlength}\raggedright $001$\end{minipage}}}%
    \put(0.43734129,0.35266195){\color[rgb]{0,0,0}\makebox(0,0)[lt]{\begin{minipage}{0.12073128\unitlength}\raggedright $111$\end{minipage}}}%
    \put(0.61734076,0.28241825){\color[rgb]{0,0,0}\makebox(0,0)[lt]{\begin{minipage}{0.12073128\unitlength}\raggedright $121$\end{minipage}}}%
    \put(0.48124355,0.25607696){\color[rgb]{0,0,0}\makebox(0,0)[lt]{\begin{minipage}{0.21467115\unitlength}\raggedright $\tau_{A^2}$\\ \end{minipage}}}%
    \put(0.43734129,0.3043694){\color[rgb]{0,0,0}\makebox(0,0)[lt]{\begin{minipage}{0.12073128\unitlength}\raggedright $332$\end{minipage}}}%
    \put(0.53831661,0.29558894){\color[rgb]{0,0,0}\makebox(0,0)[lt]{\begin{minipage}{0.12073128\unitlength}\raggedright $453$\end{minipage}}}%
  \end{picture}%
\endgroup%

%% file: example3a.pdf_tex
\begingroup%
  \makeatletter%
  \providecommand\color[2][]{%
    \errmessage{(Inkscape) Color is used for the text in Inkscape, but the package 'color.sty' is not loaded}%
    \renewcommand\color[2][]{}%
  }%
  \providecommand\transparent[1]{%
    \errmessage{(Inkscape) Transparency is used (non-zero) for the text in Inkscape, but the package 'transparent.sty' is not loaded}%
    \renewcommand\transparent[1]{}%
  }%
  \providecommand\rotatebox[2]{#2}%
  \newcommand*\fsize{\dimexpr\f@size pt\relax}%
  \newcommand*\lineheight[1]{\fontsize{\fsize}{#1\fsize}\selectfont}%
  \ifx\svgwidth\undefined%
    \setlength{\unitlength}{341.6678772bp}%
    \ifx\svgscale\undefined%
      \relax%
    \else%
      \setlength{\unitlength}{\unitlength * \real{\svgscale}}%
    \fi%
  \else%
    \setlength{\unitlength}{\svgwidth}%
  \fi%
  \global\let\svgwidth\undefined%
  \global\let\svgscale\undefined%
  \makeatother%
  \begin{picture}(1,0.89059585)%
    \lineheight{1}%
    \setlength\tabcolsep{0pt}%
    \put(0,0){\includegraphics[width=\unitlength,page=1]{example3a.pdf}}%
    \put(-0.2230795,-2.14191836){\color[rgb]{0,0,0}\makebox(0,0)[lt]{\begin{minipage}{0.24889485\unitlength}\raggedright \end{minipage}}}%
    \put(0,0){\includegraphics[width=\unitlength,page=1]{example3a.pdf}}%
    \put(0.01587916,0.04973605){\color[rgb]{0,0,0}\makebox(0,0)[lt]{\begin{minipage}{0.16198862\unitlength}\raggedright $1 0 0$\end{minipage}}}%
    \put(0.92026711,0.04973605){\color[rgb]{0,0,0}\makebox(0,0)[lt]{\begin{minipage}{0.12073128\unitlength}\raggedright $010$\end{minipage}}}%
    \put(0.4812436,0.87070922){\color[rgb]{0,0,0}\makebox(0,0)[lt]{\begin{minipage}{0.12073128\unitlength}\raggedright $001$\end{minipage}}}%
    \put(0.43295105,0.33949125){\color[rgb]{0,0,0}\makebox(0,0)[lt]{\begin{minipage}{0.12073128\unitlength}\raggedright $111$\end{minipage}}}%
    \put(0.43295105,0.46241772){\color[rgb]{0,0,0}\makebox(0,0)[lt]{\begin{minipage}{0.18613465\unitlength}\raggedright $112$\\ \end{minipage}}}%
    \put(0.56904822,0.24290617){\color[rgb]{0,0,0}\makebox(0,0)[lt]{\begin{minipage}{0.12073128\unitlength}\raggedright $121$\end{minipage}}}%
    \put(0.61734076,0.38339356){\color[rgb]{0,0,0}\makebox(0,0)[lt]{\begin{minipage}{0.12073128\unitlength}\raggedright $122$\end{minipage}}}%
    \put(0.52075563,0.36583269){\color[rgb]{0,0,0}\makebox(0,0)[lt]{\begin{minipage}{0.21467115\unitlength}\raggedright $\tau_A$\\ \end{minipage}}}%
  \end{picture}%
\endgroup%